\documentclass{asl}
\usepackage{lineno,hyperref}
 
\usepackage{amsfonts,amssymb}
\usepackage{verbatim}
\usepackage{placeins}

\newtheorem{Theorem}{Theorem}[section] 
 
\newtheorem{Lemma}[Theorem]{Lemma}
\newtheorem{Corollary}[Theorem]{Corollary}
\newtheorem{Definition}[Theorem]{Definition}

\def\axioms{\begin{tabular*}{0.8\textwidth}{@{\extracolsep{\fill}}lr}}
\def\newaxioms#1{\begin{tabular*}{0.8\textwidth}{@{}l@{\hskip-#1cm}r}}
\def\endaxioms{\end{tabular*}\smallskip}

\def\seq{\Rightarrow}

\def\B{{\bf B}}
\def\T{{\bf T}}   

\def\implies{\ \rightarrow\ }
\def\defined{\downarrow}
\def\Left{\mbox{\it Left\,}}
\def\Right{\mbox{\it Right\,}}

\def\MakePoint{\mbox{\it MakePoint\,}}

\def\Line{\mbox{\it Line\,}}

\def\Circle{\mbox{\it Circle\,}}

\def\il{i\ell}
\def\ilcone{i\ell c_1}
\def\ilctwo{i\ell c_2}
\def\icone{ic_1}
\def\ictwo{ic_2}
 
\def\Reflect{\mbox{\it Reflect\,}}

\def\SameSide{\mbox{\it SameSide\,}}
\def\OppositeSide{\mbox{\it OppositeSide\,}}

\def\HilbertMultiply{\mbox{\it HilbertMultiply\,}}

\def\Add{\mbox{\it Add\,}}

\def\R{{\mathbb R}}

\def\K{{\mathbb K}}

\def\EF{{\bf EF}}

\def\LPT{{\bf LPT}}

\def\F{{\mathbb F}}

\usepackage{graphicx}
\usepackage{pstricks}
\psset{unit=3cm}

\def\EuclidParallelRawFigure{%
\pspicture(3.2, 1.2)
\pspolygon[fillstyle=solid,fillcolor=lightgray](0.67,0.7)(0.78,0.425)(1.49,0.7)
\pspolygon[fillstyle=solid,fillcolor=lightgray](0.07,0.15)(0.78,0.425)(0.9,0.15)
\qline(-0.15,0.15)(2.47,0.15)
\qline(2.53,0.15)(2.7,0.15)
\qline(-0.15,0.7)(2.7,0.7)
\qline(0.0,0.9)(2.473,0.163)
\qline(2.5275,0.14)(2.7,0.085)
\pscircle(2.5,0.15){0.03}
\psdot(0.67,0.7)
\put(0.67,0.78){$p$}
\psdot(1.29,0.515)
\put(1.28,0.39) {$a$}
\psdot(0.9,0.15)
\put(0.88,0.04) {$q$}
\psdot(0.07,0.15)
\put(0.07,0.04){$s$}
\qline(0.67,0.7)(0.9,0.15)
\psdot(1.49,0.7)
\put(1.5,0.78) {$r$}
\qline(1.49,0.7)(0.07,0.15)  
\psdot(0.78,0.425)
\put(0.65,0.3){$t$}
\qline(0.9,0.15)(1.49,0.7)  
\put(-0.3,0.12) {$L$}
\put(-0.3,0.68) {$K$}
\put(-0.15,0.88) {$M$}
\endpspicture}

\def\AlternateInteriorAnglesFigure{%
\pspolygon[fillstyle=solid,fillcolor=lightgray](0.7,0.7)(1.3,0.7)(0.8,0.425)
\pspolygon[fillstyle=solid,fillcolor=lightgray](0.9,0.15)(0.3,0.15)(0.8,0.425)
\pspicture(3.2, 1.2)
\qline(0.0,0.15)(2,0.15)
\qline(0.0,0.7)(2,0.7)
\psdot(0.7,0.7)
\put(0.67,0.78){$p$}
\psdot(0.9,0.15)
\put(0.88,0.04) {$q$}
\psdot(0.3,0.15)
\put(0.3,0.04){$s$}
\psdot(0.8,0.425)
\put(0.85,0.36){$t$}
\qline(0.3,0.15)(1.3,0.7)
\qline(0.7,0.7)(0.9,0.15)  
\psdot(1.3,0.7)
\put(1.3,0.78) {$r$}
\put(-0.15,0.12) {$L$}
\put(-0.15,0.68) {$K$}
\endpspicture}

\def\TarskiFiveSegmentFigure{%
\psset{unit=2.25cm}
\pspicture(2.5,0)(2.2, 1.2)  
\qline(0.0,0.15)(2.0,0.15)
\qline(0.0,0.15)(1.0,0.85)
\psline[linestyle=dashed](1.0,0.85)(2.0,0.15)
\qline(1.0,0.85)(0.7,0.15)
\put(1,0.9){$d$}
\put(0.0,0){$a$}
\put(0.7,0){$b$}
\put(2.0,0){$c$}
\qline(2.5,0.15)(4.5,0.15)
\qline(2.5,0.15)(3.5,0.85)
\psline[linestyle=dashed](3.5,0.85)(4.5,0.15)
\qline(3.5,0.85)(3.2,0.15)
\put(3.5,0.9){$D$}
\put(2.5,0){$A$}
\put(3.2,0){$B$}
\put(4.5,0){$C$}
\psset{unit=3cm}
\endpspicture}

\def\InteriorFiveFigure{%
 
\psset{unit=2.25cm}
\pspicture(2.5,0)(2.2, 1.2)  
\qline(0.0,0.15)(2.0,0.15)
\qline(0.0,0.15)(1.0,0.85)
\qline(1.0,0.85)(2.0,0.15)
\psline[linestyle=dashed](1.0,0.85)(0.7,0.15)
\put(1,0.9){$d$}
\put(0.0,0){$a$}
\put(0.7,0){$c$}
\put(2.0,0){$b$}
\qline(2.5,0.15)(4.5,0.15)
\qline(2.5,0.15)(3.5,0.85)
\qline(3.5,0.85)(4.5,0.15)
\psline[linestyle=dashed](3.5,0.85)(3.2,0.15)
\put(3.5,0.9){$D$}
\put(2.5,0){$A$}
\put(3.2,0){$C$}
\put(4.5,0){$B$}
\psset{unit=3cm}
\endpspicture}

\def\LineCircleContinuityFigureTwo{%
\hskip2cm   
\pspicture(3,2.2)
\psdot(1.45,1.05)
\put(1.45,1.12){$a$}
\pscircle(1.45,1.05){1.05}
\qline(0.2,1.25)(2.4,1.8)
\put(2.45,1.7){$L$}
\psdot(0.5,1.475)
\put(0.39,1.45){$b$}
\qline(0.5,1.475)(1.45,1.05)  
\psdot(0.72, 1.38)
\put(0.74,1.25){$p$}
\pscircle(2.22,1.755){0.035}
\pscircle(0.436,1.31){0.035}
\endpspicture}

\def\SegmentCircleContinuityFigure{%
\hskip2cm   
\pspicture(3,2.2)
\psdot(1.45,1.05)
\put(1.45,0.92){$a$}
\pscircle(1.45,1.05){1.05}
\qline(0,1.20)(2.4,1.8)
\psdot(0.2,1.25)
\put(0.13,1.33){$q$}
\psdot(0.4,1.05)
\put(0.31,0.9){$b$}
\qline(0,1.05)(1.45,1.05)
\psdot(0.65,1.05)
\put(0.7,0.92){$x$}
\psdot(0.18,1.05)
\put(0.08,0.92){$y$}
\psarcn(1.45,1.05){0.80}{190}{125}
\psarcn(1.45,1.05){1.27}{200}{145}
\psdot(0.72, 1.38)
\put(0.74,1.25){$p$}
\pscircle(0.436,1.31){0.035}
\endpspicture}

\def\CircleCircleContinuityFigure{%
\hskip2cm  
\pspicture(3,2.4)
\psdot(2.0,1.0)
\put(2,0.87){$c$}
\psarc(1.95,1.0){1.045}{-146.7}{146.7}
\psarc(1.95,1.0){1.045}{149.5}{210.5}
\psdot(0.5,1.0)
\put(0.5,0.87){$a$}
\qline(0.5,1.0)(-0.3,1.8)
\psdot(-0.3,1.8)
\put(-0.33,1.7){$z$}
\psdot(-0.069,1.564)
\put(-0.195,1.55){$b$}
\psdot(0.1,1.4)
\put(0.05,1.29){$x$}
\psarc(0.5,1.0){1.133}{46.5}{135}
\psarc(0.5,1.0){0.57}{38}{135}
\psdot(1.288,1.811)
\put(1.28,1.7){$q$}
\psdot(2.612,1.811)
\put(2.54,1.67){$d$}
\psdot(0.962,1.338)
\put(0.975,1.25){$p$}
\psline{->}(0.5,1)(2,1)
\psarc(0.5,1.0){0.8}{46.5}{313.5}
\psarc(0.5,1.0){0.8}{317}{43}
\pscircle(1.066,1.557){0.03}
\put(1.03,1.65){$2$}
\qline (1.086,1.545)(2,1)
\pscircle(1.066,0.443){0.03}
\put(1.03,0.3){$1$}
\put(2.7,0.1){$K$}
\put(0,0.1){$C$}
\endpspicture}

\def\InnerOuterPaschFigure{%
\hskip0.9cm   
\pspicture(3.6,1.8)
\psdot(0,0.6)
\put(0,0.47){$a$}
\pscircle(1,0.6){0.03}
\put(1,0.47){$x$}
\qline(0,0.6)(0.97,0.6)
\psdot(1.3,0.6)
\put(1.35,0.55){$q$}
\qline(1.03,0.6)(1.3,0.6)
\psdot(1.1,1.4)
\put(1.1,1.45){$c$}
\qline(0,0.6)(1.1,1.4)  
\qline(1.1,1.4)(1.375,0.3) 
\psdot(1.375,0.3)
\put(1.43,0.25){$b$}
\psdot(0.52,0.98)
\put(0.43,1.03){$p$}
\qline(1.375,0.3)(1.015,0.584)  
\qline(0.52,0.98)(0.979,0.62)  
\psdot(2,0.3)
\put(1.95,0.19){$a$}
\psdot(3.5,0.3)
\put(3.45,0.19){$q$}
\qline(2,0.3)(2.47,0.3)
\qline(2.53,0.3)(3.5,0.3)
\pscircle(2.5,0.3){0.03}
\put(2.45,0.19){$x$}
\psdot(2.7,1.5)
\put(2.65,1.57){$b$}
\qline(2.7,1.5)(3.5,0.3)  
\psdot(3.0,1.05)
\put(3.07,1.02) {$c$}
\qline(2.7,1.5)(2.5,0.33)  
\qline(2,0.3)(3.0,1.05)   
\psdot(2.566,0.722)
\put(2.6,0.6){$p$}
\endpspicture}

\def\ContinuousPaschFigure{%
\hskip0.9cm   
\pspicture(-0.5,0.2)(3.0,2.0)
\psdot(0,0.6)
\put(0,0.47){$a$}
\pscircle(1,0.6){0.03}
\put(0.98,0.47){$x$}
\qline(0,0.6)(0.97,0.6)
\psdot(1.3,0.6)
\put(1.23,0.49){$q$}
\qline(1.03,0.6)(1.39,0.6)
\psdot(1.1,1.4)
\put(1.07,1.45){$c$}
\psdot(1.6,1.764)
\put(1.63,1.68){$d$}
\qline(0,0.6)(1.1,1.4)  
\qline(0,0.6)(1.6,1.764) 
\pscircle(1.421,0.6){0.03} 
\qline(1.375,0.3)(1.42,0.57) 
\qline(1.422,0.63)(1.6,1.764) 
\qline(1.1,1.4)(1.375,0.3) 
\psdot(1.375,0.3)
\put(1.43,0.25){$b$}
\psdot(0.52,0.98)
\put(0.43,1.03){$p$}
\qline(1.375,0.3)(1.015,0.584)  
\qline(0.52,0.98)(0.979,0.62)  
\psdot(1.35,1.582)  
\put(1.32,1.65){$p?$}
\psline[linestyle=dashed](1.375,0.3)(1.35,1.582)
\endpspicture}

\def\SameSideFigure{%
\pspicture(-0.6,0.4)(3.0,1.8)
\qline(0.7,1)(2,1)
\put(0.53,0.97){$L$}
\psdot(0.9,0.5)
\put(0.83,0.37){$a$}
\psdot(1.2,0.5)
\put(1.17,0.37){$b$}
\psdot(1.8,1.5)
\put(1.85,1.5){$c$}
\psline(1.8,1.5)(1.2,0.5)
\psline(1.8,1.5)(0.9,0.5)
\psdot(1.5,1)
\put(1.54,0.9){$y$}
\psdot(1.346,1)
\put(1.27,1.05){$x$}
\endpspicture
}

\def\TarskiParallelFigure{%
\hskip2cm  
\pspicture(2.5,1.2)  
\qline(0.03,0)(2.47,0)  
\pscircle(0,0){0.03}
\put(-0.13,-0.02){$x$}
\pscircle(2.5,0){0.03}
\put(2.55,-0.02){$y$}
\psdot(0.5,1)
\put(0.47,1.05){$a$}
\psdot(0.167,0.333)
\put(0.06,0.33){$b$}
\qline(0.5,1)(2.48,0.022)  
\qline(0.01,0.025)(0.5,1)  
\psdot(1.1,0.703)
\put(1.1,0.74){$c$}
\qline(0.167,0.333)(1.1,0.703)  
\qline(0.5,1)(1.0,0)  
\psdot(1.0,0)
\put(1.02,0.05){$t$}
\psdot(0.723,0.553)
\put(0.78,0.48){$d$}
\endpspicture
}

\def\TarskiParallelProofFigure{%
\pspicture(2.5,2.3)(-1,0.6)
\pscircle(-0.2,0.6){0.03}
\put(-0.32,0.6){$x$}
\psdot(0.5,2)
\put(0.47,2.05){$a$}
\psdot(0.1667,1.333)
\put(0.06,1.33){$b$}
\qline(0.5,2)(1.881,1.314)  
\qline(-0.19,0.625)(0.5,2)  
\psdot(1.17,1.667)
\put(1.18,1.72){$c$}
\qline(0.5,2)(1.0,1)  
\psdot(1,1)
\put(0.99,0.88){$t$}
\psdot(0.74,1.52)
\put(0.79,1.45){$d$}
\pscircle(1.8969,1.30154){0.03}
\put(1.9,1.2){$y$}
\qline(-0.172,0.61)(1.87,1.295)  
\qline(1.92,1.31)(2.393,1.46)  
\put(2.393,1.35){$M$}
\psdot(1.60,1.81)
\put(1.635, 1.81){$e$}
\psdot(-0.4055,1.1427)
\put(-0.4055,1.044){$g$}
\qline(-0.4055,1.1427)(1.60,1.81) 
\qline(-0.4055,1.1427)(1,1)   
\qline(1.60,1.81)(1,1)   
\psdot(1.405,1.55)
\put(1.47,1.56){$f$}
\psdot(0.045,1.093)
\put(0.05,0.97){$h$}
\endpspicture
}

\def\TarskiImpliesEuclidFigure{%
\pspicture(3.2, 1.1)(-0.8,-0.45)
\pspolygon[fillstyle=solid,fillcolor=lightgray](0.67,0.7)(0.78,0.425)(1.49,0.7)
\pspolygon[fillstyle=solid,fillcolor=lightgray](0.07,0.15)(0.78,0.425)(0.9,0.15)
\qline(-0.15,0.15)(1.22,0.15)  
\qline(1.28,0.15)(2.7,0.15)    
\qline(-0.15,0.7)(2.7,0.7)     
\pscircle(1.25,0.15){0.03}    
\put(1.22,0.04){$e$}
\qline(0.52,0.85)(1.225,0.17)    
\qline(1.27,0.13)(1.93,-0.5)    
\qline(2,-0.5)(0,0.7)    
\psdot(0,0.7)
\put(0,0.75){$x$}
\psdot(1.83,-0.4)
\put(1.7,-0.44){$y$}
\psdot(0.67,0.7)
\put(0.67,0.75){$p$}
\psdot(0.5,0.7)
\put(0.45,0.75){$u$}
\psdot(1.075,0.315)
\put(1.14,0.29) {$a$}
\psline[linestyle=dotted](0.5,0.7)(1.075,0.315)
\psdot(0.725,0.55)
\put(0.62,0.5){$v$}
\psdot(0.9,0.15)
\put(0.88,0.04) {$q$}
\psdot(0.07,0.15)
\put(0.07,0.04){$s$}
\qline(0.67,0.7)(0.9,0.15)
\psdot(1.49,0.7)
\put(1.5,0.75) {$r$}
\qline(1.49,0.7)(0.07,0.15)  
\psdot(0.78,0.425)
\put(0.73,0.32){$t$}
\qline(0.9,0.15)(1.49,0.7)  
\put(-0.3,0.12) {$L$}
\put(-0.3,0.68) {$K$}
\put(0.4,0.88) {$M$}
\endpspicture
}

\def\PerpFigure{%
\pspicture(-0.5,-0.6)(2.5,2.6)
\qline(-0.5,0.6)(1.75,0.6)
\put(1.8,0.57){$L$}
\psdot(0.625,1.0)
\put(0.635,1.05){$x$}
\psdot(0.25,0.6)
\put(0.205,0.48){$a$}
\psdot(0.5,0.6)
\put(0.5,0.48){$b$}
\psdot(1.2185,0.6)
\put(1.2,0.45){$c$}
\pscircle(0.625,1.0){1.0185}
\put(1.7,1.2){$C$}
\psdot(-0.305,0.6)
\put(-0.38,0.48){$p$}
\psdot(1.556,0.6)
\put(1.556,0.48){$q$}
\put(0.66,2.2){$K$}
\qline(0.625,2.3)(0.625,-0.5)
\endpspicture}

\def\GuptaMidpointFigure{%
\psset{unit=2cm}
\pspicture(0,0)(2,1.6)  
\psdot(0,0)
\put(-0.15,-0.15){$x$}
\psdot(1.414,0)
\put(1.4,-0.15){$z$}
\psdot(1.414,1.414)
\put(1.48,1.4){$t$}
\psdot(2,0)
\put(2,-0.15){$u$}
\psdot(1,1)
\put(0.9,1.1){$y$}
\qline(0,0)(2,0)  
\qline(0,0)(1.414,1.414) 
\qline(1.414,0)(1.414,0.545) 
\qline(1.414,0.62)(1.414,1.414)  
\qline(1,1)(1.397,0.60)  
\qline(2,0)(1.43,0.57)  
\qline(1,1)(1.20,0.53) 
\qline(1.225,0.47)(1.414,0) 
\qline(0,0)(1.187,0.491)
\qline(1.222,0.512)(1.39,0.575)
\qline(1.45,0.594)(2,0.8284)
\pscircle(1.207,0.5){0.04}
\put(1.1,0.34){$w$}
\pscircle(1.414,0.5857){0.04}
\put(1.52, 0.525){$v$}
\psset{unit=3cm}
\endpspicture}

\def\GuptaPerpendicularFigure{%
\pspicture(-0.2,-1.2)(4.2,1.3)
\psset{unit=2cm}
\qline(0.8,0)(2.57,0)
\pscircle(2.6,0){0.03}
\put(2.65,-0.16){$x$}
\qline(2.63,0)(5.1,0)
\psdot(4.2,0)
\put(4.15,-0.16){$a$}
\psdot(5.0,0)
\put(4.95,-0.18){$b$}
\psdot(2,0)
\put(2.12,-0.16){$y$}
\psdot(2.6,1.51)
\put(2.7,1.5){$c$}
\psdot(1.3,0)
\put(1.35,-0.16){$z$}
\psdot(2.6,-1.51)
\put(2.7,-1.55){$c^\prime$}
\qline(2.6,-1.51)(2.6,-0.03)  
\qline(2.6,0.03)(2.6,1.51)  
\qline(4.2,0)(2.6,1.51)   
\psdot(1.3,1.762)
\put(1.38,1.75){$q^\prime$}
\psdot(1.3,-1.762)
\put(1.4,-1.8){$q$}
\qline(1.3,1.762)(2.6,-1.51)  
\qline(1.3,-1.762)(2.6,1.51)  
\qline(1.3,-1.762)(1.3,1.762) 
\psdot(2.3,0.755)
\put(2.15,0.75){$p$}
\qline(2.3,0.755)(4.2,0)   
\psset{unit=3cm}
\endpspicture}

\def\ErectPerpFigure{%
\pspicture(-1.2,0)(2.5,1.5) 
\qline(-0.5,0.6)(1.75,0.6)
\put(1.8,0.57){$L$}
\psdot(0.6,0)
\put(0.65,-0.07){$c$}
\psdot(0.6,1.2)
\put(0.65,1.25){$d$}
\psdot(0.6,0.6)
\put(0.65,0.5){$x$}
\psdot(1.1,0.6)
\put(1.17,0.525){$a$}
\psdot(1.1,1.2)
\put(1.13,1.25){$p$}
\psdot(-0.3,0.6)
\put(-0.38,0.5){$b$}
\psline[linecolor=blue](0.6,0)(1.6,0)
\psline[linecolor=blue](0.6,0)(1.1,1.2)
\psline[linecolor=blue](1.6,0)(1.1,1.2)
\psdot(1.1,0)
\put(1.12,-0.07){$q$}
\psdot(1.6,0.6)
\put(1.625,0.5){$y$}
\psdot(1.6,1.2)
\put(1.556,1.25){$e$}
\psdot(1.6,0)
\put(1.625,-0.08){$f$}
\psdot(-0.1,0.47)
\put(-0.13,0.37){$s$}
 \qline(0.6,1.2)(0.6,0)
\qline(1.6,1.2)(1.6,0)
\qline(0.6,1.2)(1.6,0)
\qline(0.6,0)(1.6,1.2)
\qline(0.6,1.2)(1.6,1.2)
\qline(-0.3,0.6)(0.6,0)
\psline[linestyle=dashed](1.1,1.2)(0.6,0.6)
\psline[linestyle=dashed](1.1,1.2)(1.6,0.6)
\psline[linecolor=red](1.1,1.2)(1.1,0)
\endpspicture}

\def\ErectPerpFinishFigure{%
\pspicture(-1.3,0)(2.5,1.5)  
\qline(-0.5,0.6)(0.08,0.6)
\qline(0.14,0.6)(0.82,0.6)
\qline(0.88,0.6)(1.75,0.6)
\put(1.8,0.57){$L$}
\psdot(0.6,0)
\put(0.65,-0.07){$c$}
\psdot(0.6,1.2)
\put(0.65,1.25){$d$}
\psdot(0.6,0.6)
\put(0.65,0.48){$x$}
\psdot(1.1,0.6)
\put(1.08,0.48){$a$}
\psdot(1.1,1.2)
\put(1.13,1.25){$p$}
\psdot(-0.3,0.6)
\put(-0.38,0.48){$b$}
\psdot(1.6,1.2)
\put(1.556,1.25){$e$}
\psdot(-0.1,0.47)
\put(-0.13,0.35){$s$}
\qline(0.6,1.2)(0.6,0)
\qline(0.6,0)(1.6,1.2)
\qline(0.6,1.2)(1.6,1.2)
\qline(-0.3,0.6)(0.6,0)
\psline[linestyle=dashed](1.1,1.2)(0.865,0.63) 
\psline[linestyle=dashed](0.6,0)(0.835,0.57)
\pscircle(0.85,0.6){0.03}
\put(0.86,0.48){$t$}
\pscircle(0.11,0.6){0.03}
\put(0.11,0.48){$r$}
\psline[linestyle=dashed](1.1,1.2)(0.125,0.605)
\psline[linestyle=dashed](0.085,0.58) (-0.1,0.475)

\endpspicture}

\def\CenterImpliesStrongParallelFigure{%
\pspicture(-0.5,-0.5)(3,1.5)
\qline(0.2,0)(2.97,0)
\put(0.1,-0.05){$L$}
\qline(0,1)(2.9,1)
\put(2.93,0.95){$K$}
\psdot(1,1)
\put(1.02,0.89){$p$}
\pscircle(3,0){0.03}
\put(3,-0.1){$e$}
\psline(0.75,1.125)(2.98,0.015)
\put(0.6,1.125){$M$}
\psdot(2,0.5)
\put(2.05,0.53){$r$}
\psdot(1,0)
\put(1,-0.1){$w$}
\qline(1,0)(1,1)
\qline(0.4,-0.2)(1.225,1.45)
\psdot(0.4,-0.2)
\put(0.45,-0.22){$j$}
\psdot(0.775,0.55)
\put(0.64,0.5){$x$}
\psdot(1.225,1.45)
\put(1.14,1.45){$y$}
\psdot(1,0.5)
\put(1.03,0.57){$m$}
\pscircle[linestyle=dotted](1,1){0.5} 
\psdot(0.78,-0.55)
\put(0.66,-0.55){$z$}
\qline(0.78,-0.55)(0.78,0.55)
\psline[linestyle=dashed](0.78,-0.55)(2.97,-0.01)
\psline[linestyle=dashed](0.78,0.55)(2.97,0.01)
\psline[linestyle=dashed](1.225,1.45)(2.97,0.03)
\endpspicture}

\def\PointnotonLFigure{%
\pspicture(0,-1)(4,1.7)
\qline(0.2,0)(3.8,0)
\put(0.1,-0.05){$L$}
\psdot(1,1)
\put(1.02,0.89){$q$}
\pscircle(3.5,-0.225){0.03}
\put(3.5,-0.325){$c$}
\psline(0.75,1.125)(3.49,-0.22) 
\put(0.6,1.125){$K$}
\psdot(2,0.51)
\put(2.05,0.53){$r$}
\qline(0.4,-0.2)(1.225,1.45)
\put(0.45,-0.22){$M$}
\psdot(0.775,0.55)
\put(0.64,0.5){$u$}
\psdot(1.225,1.45)
\put(1.14,1.45){$v$}
\pscircle[linestyle=dotted](1,1){0.5} 
\psdot(0.78,-0.55)
\put(0.66,-0.55){$w$}
\psdot(0.78,-0.9)
\put(0.66,-0.9){$d$}
\qline(0.78,-1)(0.78,0.55)
\psline[linestyle=dashed](0.78,-0.9)(3.5,-0.235)
\psline[linestyle=dashed](0.78,0.55)(3.5,-0.225)
\psline[linestyle=dashed](1.225,1.45)(3.5,-0.21)
\psline[linestyle=dotted](0.78,-0.23)(3.5,-0.23)
\endpspicture}

\def\PerpFromCircleFigure{%
\pspicture(-1,-0.2)(2,1.5)
\psdot(0,0)
\put(-0.1,-0.14){$a$}
\psdot(1,0)
\put(1,-0.14){$b$}
\psdot(0,1)
\put(-0.1,1.03){$e$}
\psdot(0.5,0.5)
\put(0.57, 0.5){$c$}
\pscircle (0.5,0.5){0.707}
\qline(-0.5,0)(1.5,0)
\put(1.55,-0.03){$L$}
\qline(1,0)(0,1)
\psline[linecolor=red](0,-0.2)(0,1.2)
\endpspicture}

\def\MidpointButterflyFigure{%
\pspicture(-2,-0.9)(2,1)
\psset{unit=2.5cm}
\psdot(-2,0)
\put(-2.13,0.03){$a$}
\psdot(2,0)
\put(2.05,-0.02){$b$}
\qline(-2,0)(2,0)
\qline(-2,-1)(-2,1)
\qline(2,-1)(2,1)
\psdot(-2,-1)
\put(-2.13,-0.97){$p$}
\psdot(-2,1)
\put(-2.13,1){$q$}
\qline(-2,-1)(2,1)
\psdot(0,0)
\put(0,-0.15){$m$}
\psdot(2,1)
\put(2.05,1){$r$}
\psdot(2,-1)
\put(2.05,-1){$s$}
\psline[linestyle=dashed](-2,1)(2,1)
\psline[linestyle=dashed](-2,-1)(2,-1)
\psline[linestyle=dashed](-2,1)(2,-1)

\psset{unit=3cm}
\endpspicture}

\def\MidpointPaschFigure{%
\pspicture(-1.8,-0.7)(2,1)
\psset{unit=2cm}
\psdot(-2,0)
\put(-2.15,0){$a$}
\psdot(2,0)
\put(1.97,-0.17){$b$}
\psdot(-0.265,0)
\put(-0.3,0.07){$t$}
\qline(-2,0)(-0.03,0)
\qline(0.03,0)(2,0)
\qline(-2,-1)(-2,0)
\qline(2,0)(2,1.3)
\psdot(-2,-1)
\put(-2.15,-0.98){$p$}
\psdot(2,1.3)
\put(1.86,1.34){$q$}
\qline(-2,-1)(-0.03,-0.02)
\pscircle(0,0){0.03}
\qline(0.03,0.02)(2,1)
\put(0.03,-0.14){$m$}
\qline(-2,-1)(2,1.3)
\psdot(2,1)
\put(2.06,0.95){$r$}
\psdot(2,0.3)
\put(2.06,0.25){$w$}
\qline(-2,-1)(2,0.3)
\psdot(1.092,0)
\put(1.092,-0.14){$v$}
\psset{unit=3cm}
\endpspicture}

\def\RadicalAxisFigure{%
\pspicture(-1,-0.7)(2,2.5)
\psset{unit=0.6cm}
\psdot(-0.94,5.26)
\put(-1.2,4.8){$s$}
\psdot(8.22,4.64)
\put(8.3,4.1){$t$}
\psdot(-1.8,-0.99)
\put(-1.9,-1.6){$a$}
\psdot(8.37,0.21)
\put(8.3,-0.45){$b$}
\qline(-0.94,5.26)(8.22,4.64)   
\psdot(3.28,-0.39)
\put(3.28,-1){$m$}
\psdot(4.64,3.84)
\put(4.05,3.8){$p$}
\qline(-1.8,-0.99)(8.37,0.21)  
\qline(-0.94,5.26)(3.28,-0.39)  
\qline(8.24,4.64)(3.28,-0.39)   
\pscircle(-0.94,5.26){6.3}  
\pscircle(8.22,4.64){4.45} 
\psline[linecolor=blue](-4,-2.651)(10,7.85)  
\psline[linecolor=blue](0.35,8)(9.64,-1)  
\psline[linecolor=red](4,-4)(5.1,10)
\put(-6.5,7){$C$}
\put(11.5,5){$K$}
\psset{unit=3cm}
\endpspicture}

\def\RadicalAxisProofFigure{%
\pspicture(-1,-1)(2,2.5)
\psset{unit=0.6cm}
\pscircle[linecolor=green](3.28,-0.39){5.15}
\psdot(-0.94,5.26)
\put(-1.2,4.8){$s$}
\psdot(8.22,4.64)
\put(8.3,4.1){$t$}
\psdot(-1.8,-0.99)
\put(-1.9,-1.6){$a$}
\psdot(8.37,0.21)
\put(8.3,-0.45){$b$}
\qline(-0.94,5.26)(8.22,4.64)   
\psdot(3.28,-0.39)
\put(3.28,-1){$m$}
\psdot(4.64,3.84)
\put(4.05,3.8){$p$}
\qline(-1.8,-0.99)(8.37,0.21)  
\qline(-0.94,5.26)(3.28,-0.39)  
\qline(8.24,4.64)(3.28,-0.39)   
\pscircle(-0.94,5.26){6.3}  
\pscircle(8.22,4.64){4.45} 
\psline[linecolor=blue](-4,-2.651)(10,7.85)  
\psline[linecolor=blue](0.35,8)(9.64,-1)  
\psline[linecolor=red](4,-4)(5.1,10)
\put(-6.5,7){$C$}
\put(11.5,5){$K$}
\put(3,-5){$M$}
\psset{unit=3cm}
\endpspicture}

\def\RadicalAxisProofFigureTwo{%
\pspicture(-1,-1)(2,2.5)
\psset{unit=0.6cm}
\pscircle[linecolor=green](3.28,-0.39){5.15}
\psdot(-0.94,5.26)
\put(-1.2,4.8){$s$}
\psdot(8.22,4.64)
\put(8.3,4.1){$t$}
\psdot(-1.8,-0.99)
\put(-1.9,-1.6){$a$}
\psdot(8.37,0.21)
\put(8.3,-0.45){$b$}
\qline(-0.94,5.26)(8.22,4.64)   
\psdot(3.28,-0.39)
\put(3.28,-1){$m$}
\psdot(4.64,3.84)
\put(4.05,3.8){$p$}
\psdot(5.08,3.42)
\put(5.3,3.43){$y$}
\qline(-1.8,-0.99)(8.37,0.21)  
\qline(-0.94,5.26)(3.28,-0.39)  
\qline(8.24,4.64)(3.28,-0.39)   
\pscircle(-0.94,5.26){6.3}  
\pscircle(8.22,4.64){4.45} 
\psline[linecolor=blue](-4,-2.651)(10,7.85)  
\psline[linecolor=blue](-2.95,11.2)(9.64,-1)  
\psdot(-2.95,11.2)
\put(-3.1,10.6){$v$}
\psline[linecolor=red](4,-4)(5.1,10)
\put(-6.5,7){$C$}
\put(11.5,5){$K$}
\put(3,-5){$M$}
\psdot(3.77,4.7)
\put(3.1,4.4){$x$}
\psline[linestyle=dashed](3.77,4.7)(-1.8,-0.99)
\psset{unit=3cm}
\endpspicture}

\def\PowerFigure{%
\pspicture(0,0)(2,1)   
\psset{unit=0.3cm}
\pscircle(4.22,4.64){4.45} 
\qline(-2,0)(10,7.85)  
\qline(-1,8)(9.64,-0.5)  
\psdot(4.06,3.96)
\put(3.9,4.5){$b$}
\psdot(0.4 ,6.9)
\put(1,6.8){$a$}
\psdot(0.75,1.9)
\put(0.7,2.5){$d$}
\psdot(8.1,6.65)
\put(8.1,7.5){$e$}
\psdot(7.25,1.4)
\put(8,1.4){$c$}
\pscircle(20.22,4.64){4.45} 
\qline(12,6)(26,7.85)  
\qline(12,6)(25.64,-0.5) 
\psdot(12,6)
\put(11.8,4.6){$b$}
\psdot(16.2,6.55)
\put(15.3,7){$a$}
\psdot(23.5,7.5)
\put(23,6.3){$c$}
\psdot(15.8,4.2)
\put(15, 3){$d$}
\psdot(22.55,1.0)
\put(22.2,1.5){$e$}
\psline[linecolor=blue](16.2,6.55)(22.55,1.0)
\psline[linecolor=blue](23.5,7.5)(15.8,4.2)
\psset{unit=3cm}
\endpspicture}

\def\ParallelogramsFigure{%
\pspicture(-0.5,-0.1)(2,0.8)   
\psset{unit=2cm}
\psdot(0,0)
\put(0,-0.15){$a$}
\psdot(0.2,1)
\put(0.2,1.05){$b$}
\psdot(1,0)
\put(1,-0.15){$\ell$}
\psdot(2,0)
\put(2,-0.15){$d$}
\psdot(1.2,1)
\put(1.2,1.05){$n$}
\psdot(2.2,1)
\put(2.2,1.05){$c$}
\qline(0,0)(2,0)
\qline(0.2,1)(2.2,1)
\qline(0,0)(0.2,1)
\qline(2,0)(2.2,1)
\qline(1,0)(1.2,1)
\psdot(0.1,0.5)
\put(0.1,0.35){$m$}
\psdot(1.1,0.5)
\put(1,0.55){$x$}
\psdot(2.1,0.5)
\put(2.15,0.5){$k$}
\qline(0.1,0.5)(2.1,0.5)
\psline[linestyle=dashed](0,0)(2.2,1)
\psline[linestyle=dashed](0.2,1)(2,0)
\psline[linecolor=red](0,0)(0,0.5)
\psline[linecolor=red](0,0.5)(0.1,0.5)
\psline[linecolor=red](0.2,0.5)(0.2,1)
\psline[linecolor=red](2,0)(2,0.5)
\psdot(2,0.5)
\put(1.9,0.58){$f$}
\psdot(0.2,0.5)
\put(0.25,0.58){$g$}
\psdot(0,0.5)
\put(-0.2,0.5){$h$}
\psset{unit=3cm}
\endpspicture}

\def\UniformPerpFigureTwo{%
\pspicture(-0.6,-0.5)(2,1.5)  
\psset{unit=2cm}
\qline(-1,0)(2.5,0)
\put(2.55,0){$L$}
\psdot(2,0)
\put(2.05,-0.14){$p$}
\psdot(2,2)
\put(2.05,2){$q$}
\psdot(0,-0.5)
\put(-0.13,-0.5){$s$}
\psdot(0,1.5)
\put(-0.15,1.5){$t$}
\qline(0,-0.5)(2,0)
\qline(0,-0.5)(0,1.5)
\qline(0,1.5)(2,2)
\qline(2,0)(2,2)
\psdot(1,0.75)
\put(1.08,0.75){$x$}
\psdot(1,-0.25)
\put(1.05,-0.37){$k$}
\psdot(1,0)
\put(1.05,0.15){$f$}
\psdot(1,1.75)
\put(1.05,1.63){$m$}
\psline[linestyle=dashed](0,-0.5)(2,2)
\psline[linestyle=dashed](0,1.5)(2,0)
\qline(1,-0.5)(1,2)
\put(1.05,2){$J$}
\psset{unit=3cm}
\endpspicture}

\def\TwoPerpsFigure{%
\pspicture(-0.8,-0.6)(2,0.85) 
\psset{unit=1.5cm}
\qline(-1,0)(2.5,0)
\put(2.55,0){$L$}
 \psdot(1,0)
\put(0.91,0.15){$m$}
\psdot(0,0)
\put(-0.17,0.09){$x$}
\psdot(2,0)
\put(2.1,-0.17){$p$}
\psdot(2,-1)
\put(2.1,-1.1){$r$}
\psdot(2,1)
\put(2.1,1){$q$}
\psdot(0,-1)
\put(-0.18,-1.03){$s$}
\psdot(0,1)
 \put(-0.17,0.97){$t$}
\qline(0,-1.2)(0,1.3)
\put(0.05,1.3){$J$}
\qline(2,-1.2)(2,1.3)
\put(2.05,1.3){$M$}
\psline[linestyle=dashed](0,1)(2,-1)
\psline[linestyle=dashed](0,-1)(2,1)
\psset{unit=3cm}
\endpspicture}

\begin{document}

\title{A Constructive Version of Tarski's Geometry} 
\author{Michael Beeson}        

\maketitle 

\begin{abstract}
Constructivity, in this context, refers to a theory of geometry whose
axioms and language are closely related to ruler and compass constructions.
It may also refer to the use of intuitionistic (or constructive) logic, but the reader
who is interested in ruler and compass geometry but not in constructive logic,
will still find this work of interest.  
Euclid's reasoning is essentially constructive (in both senses).
Tarski's elegant and concise first-order theory of Euclidean geometry,
on the other hand, is essentially non-constructive (in both senses), even if we restrict
attention (as we do here) to the theory with line-circle   
continuity in place of first-order Dedekind completeness.  Hilbert's
axiomatization has a much more elaborate language and many more axioms,
but it contains no essential non-constructivities.   Here we exhibit three constructive
versions of Tarski's theory.  One, like Tarski's theory, has existential axioms and 
no function symbols.   We then consider a version in which function symbols are used
instead of existential quantifiers.  This theory is quantifier-free and 
proves the continuous dependence on parameters of the terms giving the intersections
of lines and circles, and of circles and circles.  The third version has a function 
symbol for the intersection point of two non-parallel, non-coincident lines, 
instead of only for intersection points produced by Pasch's axiom and the parallel axiom;
this choice of function symbols connects directly to ruler and compass constructions.
All three versions have this in common: the axioms have been modified so that the
points they assert to exist are unique and depend continuously on parameters.  This
modification of Tarski's axioms, with classical logic, has the same theorems as
Tarski's theory, but we obtain results connecting it with ruler and compass 
constructions as well.  In particular, we show that constructions involving the 
intersection points of two circles are justified, even though only line-circle 
continuity is included as an axiom.
  We obtain metamathematical results based on the G\"odel double-negation interpretation,
which permit the wholesale importation of proofs of negative theorems from classical to 
constructive geometry, and of proofs of existential theorems where the object asserted
to exist is constructed by a single construction (as opposed to several constructions 
applying in different cases).  In particular, this enables us to import the proofs 
of correctness of the geometric definitions of addition and multiplication, once
these can be given by a uniform construction.  

We also show, using cut-elimination, that objects proved
to exist can be constructed by ruler and compass.  (This was proved in \cite{beeson-kobe}
for a version of constructive geometry based on Hilbert's axioms.)
  Since these theories are 
interpretable in the theory of Euclidean fields, the independence
results about different versions of the parallel postulate given in \cite{beeson-bsl}
apply to them; and since addition and multiplication can be defined geometrically,
their models are exactly the planes over (constructive) Euclidean fields.
 \end{abstract}

\tableofcontents

\section{Introduction } 

 Euclidean geometry, as presented by Euclid,  consists of straightedge-and-compass constructions 
and rigorous reasoning about the results of those constructions.  Tarski's twentieth-century
axiomatization of geometry does not bear any direct relation to ruler and compass constructions.
Here we present modifications of Tarski's theory whose axioms correspond more closely to 
straightedge-and-compass constructions.  These theories can be considered either with 
intuitionistic (constructive) logic,  or with ordinary (``classical'') logic. 
Both versions are of interest.
 
In \cite{beeson-kobe},
we gave an axiomatization of constructive geometry based on a version of 
Hilbert's axioms (which contain no essential non-constructivities).  
In \cite{beeson-bsl}, we obtained metamathematical results about constructive 
geometry, and showed that those results do not depend on the details of the 
axiomatization.  
In this paper, we focus on formulating constructive geometry in the 
language and style that Tarski used for his well-known axiomatization of geometry.
What is striking about Tarski's theory is its use of only one sort of variables,
for points, and the small number of axioms.  Here  
we give what may be the shortest possible axiomatization
of constructive geometry, following Tarski's example.%
\footnote{Readers unfamiliar with
Tarski's geometry may want to read \cite{tarski-givant}, which summarizes 
the axioms of Tarski's geometry and gives some of their history; but we do give
a basic explanation of Tarski's axioms in this paper.}

In \cite{beeson-bsl}, we discussed Euclidean constructive geometry in general terms, and 
worked informally with a theory that had three sorts of variables for points, lines, and circles.
Here, in the spirit of Tarski, we work with a one-sorted theory, with variables for points only.
In order to provide terms for points proved to exist, we need some function symbols.  Tarski's
axioms have existential quantifiers; we are interested (both classically and constructively)
in extensions of the language that provide function symbols to construct points.  Three 
of these symbols are Skolem symbols that correspond immediately to ruler and compass 
constructions: one for extending a segment $ab$ by another segment $cd$,  and
 two for the intersection points of a 
line and circle.  (In our previous work we also had function symbols for the 
intersection points of two circles; here we will prove that these are not needed, 
as the  intersection points of two circles 
are already constructible.)   Then we need a way
to construct certain intersection points of two lines.  Such points are proved to exist
 by versions of Pasch's axiom; so one obvious approach is just to provide a Skolem 
symbol for a suitable version of Pasch's axiom. (This has been done for decades by people
using theorem-provers with Tarski's axioms.)  

However, Tarski's version of Pasch's axiom allows ``degenerate cases'' in which the 
``triangle'' collapses to three points on a line, or the line through the triangle coincides
with a side of the triangle.  In these cases, the point asserted to exist is not really 
constructed by intersecting two lines and does not correspond to a ruler and compass construction.
Therefore, even with classical logic, Tarski's axioms need some modifications before they
really correspond to ruler and compass constructions.  To start with, we
 require that the points in Pasch's axiom be not collinear.  Then we have 
 to ``put back'' the two fundamental axioms about betweenness that Tarski originally had,
 but which were eliminated when Tarski and his students realized that they followed from 
 the degenerate cases of Pasch.  Finally, we have to restrict the segment-extension axiom 
 to extending non-null segments, i.e.,$ab$ with $a \neq b$,  since extending a null segment
 is not done by laying a straightedge between two points.  More formally, the extension 
 of segment $ab$ by a non-null segment $cd$ will not depend continuously on $a$ as $a$ 
 approaches $b$,  while ruler and compass constructions should depend continuously on 
 parameters.   The resulting modification of Tarski's classical axioms we call 
 ``continuous Tarski geometry''.  If we add the function symbols mentioned above,
 then all those function symbols correspond to ruler and compass constructions, 
 and Herbrand's theorem then tells us that if we can prove $\forall x \exists y\, A(x,y)$, 
 and $A$ is quantifier-free, 
 then there are finitely many ruler and compass constructions $t_1,\ldots,t_n$ such 
 that for each $x$, one of the $t_i(x)$ constructs $y$ such that $A(x,y)$.
 
 We said that ruler and compass constructions should depend continuously on parameters,
 but there is a problem about that: we need to distinguish axiomatically between
 the two intersection points of a line and a circle.  Since lines are given by two 
 distinct points, our solution to this problem is to require that the two intersection points
 of $\Line(a,b)$ and circle $C$ occur in the same order on $L$ as $a$ and $b$.  
 Thus if $a$ and $b$ are interchanged, the intersection points given by the two 
 function symbols also are interchanged.

 All the changes discussed above make sense and are desirable even with classical logic.
 They connect the axioms of geometry with ruler and compass constructions and, in the 
 case of Pasch's axiom, with its intuitive justification.  The degenerate cases of Pasch
 have nothing to do with triangles and lines;  they are really about betweenness relations
 between points on a single line, so it is philosophically better to formulate the 
 axioms as in continuous Tarski geometry.  Having the smallest possible number of axioms
 is not necessarily the criterion for the best version of a theory.
 
 There is also an issue regarding the best form of the parallel axiom.  Historically,
 several versions have been considered for use with Tarski's theories.  Two in particular
 are of interest:  the axiom (A10) that Tarski eventually settled upon, and the 
 ``triangle circumscription principle'', which says that given three non-collinear points, 
 there is a point $e$ equidistant from all three (which is then the 
 center of a circle containing the three points).  Classically, these two formulations
  are equivalent, so it is just a matter
 of personal preference which to take as an axiom.   Constructively, the two versions
 mentioned  are 
 also equivalent, as follows from the results of \cite{beeson-bsl} and this paper, but the proof is much lengthier than
 with classical logic.  Euclid's own formulation of the parallel postulate, ``Euclid 5'',
 mentions angles, so it requires a reformulation to be expressed in the ``points only'' language
 of Tarski's theory; a points-only version of Euclid 5 is given in \cite{beeson-bsl} and 
 repeated below.  In \cite{beeson-bsl} it is proved that Euclid 5 is equivalent to the 
 triangle circumscription principle,  which is considerably shorter than Euclid 5.
  We follow Szmielew in adopting the triangle circumscription
 principle as our parallel axiom, although our results show that we {\em could} have 
 retained Tarski's version. 
 
 There is also ``Playfair's axiom'', which is the version of the parallel axiom adopted 
 by Hilbert in \cite{hilbert1899}.   That version, unlike all the other versions, makes 
 no existence assertion at all, but only asserts that there cannot exist two different lines
 parallel to a given line through a given point.  This version, making no existence assertion,
 appears to be constructively weaker than the others, and in \cite{beeson-bsl}, it is proved
 that this is indeed the case.  
 
 Our aim in this paper is a constructive version of Tarski's geometry.  The changes
 described above, however, make sense with classical logic and are the primary changes
 that allow a connection between proofs from Tarski's axioms and ruler and compass 
 constructions.   If we still use classical logic, proofs in this theory yield a finite 
 number of ruler and compass constructions,  to be used in the different cases required
 in the proof.    To make the theory constructive,   we  do just
 two things more:  (1) we use intuitionistic logic instead of classical logic, and 
 (2) we add ``stability axioms'', allowing us to prove equality or inequality of 
 points by contradiction. The reasons for accepting the stability axioms are 
 discussed in \S\ref{section:stability} below. 
  
 It turns out that no more changes are needed.  This theory is called
``intuitionistic Tarski geometry''.  As in classical geometry, we can consider 
it with or without function symbols.  

Even though this theory is constructively acceptable, one might not like the 
fact that the Skolem symbols are total, i.e.,they have {\em some} (undetermined)
value even in ``undefined'' cases, where 
they do not actually correspond to ruler and compass constructions.  Therefore 
we also consider a version of Tarski geometry in which the logic is further modified
to use the ``logic of partial terms'' LPT, permitting the use of undefined terms.
In this theory, we replace the Skolem function for Pasch's axiom by a more natural 
term $\il(a,b,c,d)$ for the intersection point of  $\Line(a,b)$ and $\Line(c,d)$.

The main difference between constructive and classical geometry is that, in 
constructive geometry, one is not allowed to make a case distinction when proving 
that something exists.   For example,
to prove that there always exists a perpendicular to line $L$ through point $x$,
we may classically use two different constructions, one of which works when $x$ is not
on $L$ (a ``dropped perpendicular''), and a different construction that works when 
$x$ is on $L$ (an ``erected perpendicular'').   But constructively, we need a single 
construction (a ``uniform perpendicular'') that handles both cases with one construction.
In this paper we show that such uniform constructions can be found, using the Tarski
axioms,  for perpendiculars, reflections, and rotations.  Then the methods of \cite{beeson-bsl}
can be used to define addition and multiplication geometrically,  as was done classically 
by Descartes and Hilbert.  This shows that every model of the theory is a plane over
a Euclidean ordered field that can be explicitly constructed.  

Having formulated intuitionistic Tarski geometry, we then study its metamathematics.
using two logical tools:  the G\"odel double-negation interpretation, and cut-elimination.
The  double-negation interpretation
 is just a formal way of saying that, by pushing double negations inwards,
we can convert a classical proof of a basic statement like equality of two points,
or incidence of a point on a line, or a betweenness statement,  to a constructive proof.
(The same is of course not true for statements asserting that something exists.)
This provides us with tools for the wholesale importation of certain types of theorems
from the long and careful formal development from Tarski's classical axioms in \cite{schwabhauser}.
But since we modified Tarski's axioms, to make them correspond better to ruler and compass,
some care is required in this metatheorem.

Cut-elimination provides us with the theorem that things proved to exist in intuitionistic
Tarski geometry can be constructed by ruler and compass.  The point here is that they can
be constructed by a {\em uniform} construction, i.e.,a single construction that works
for all cases.  We already mentioned the example of dropped and erected 
perpendiculars in classical geometry, versus a uniform perpendicular construction 
in constructive geometry.  Using  cut-elimination we prove that this feature of 
constructive proofs,  so evident from examples,  is a necessary feature of any existence
proof in intuitionistic Tarski geometry:  an existence proof {\em always} provides
a uniform construction.

On the other hand, our version of Tarski geometry with classical logic, which we call 
``continuous Tarski geometry'',   supports a similar theorem.   If it proves 
$\forall x \exists y\, A(x,y)$,  with $A$ quantifier-free,  then there are a finite 
number (not just one) of ruler and compass constructions,  given terms of the theory,
such that for every $x$,  one of those constructions produces $y$ such that $A(x,y)$.
 
Readers not familiar with intuitionistic or constructive mathematics will find 
additional introductory
material in \S~\ref{section:constructivity}.  

I would like to thank Marvin J.~Greenberg and Victor Pambuccian for conversations and emails
on this subject, and the anonymous referees for their detailed comments.

\section{Hilbert and Tarski}
It is not our purpose here to review in detail the (long, complicated, and interesting)
history of axiomatic geometry, but some history is helpful in understanding the variety of geometrical
axiom systems. We begin by mentioning the standard English translation of Euclid \cite{euclid1956}
and the beautiful commentary-free edition \cite{euclid2007}.
Euclid is the touchstone against which axiomatizations are measured.
 We restrict our attention to the two most famous axiomatizations, those of Hilbert and 
Tarski.  Previous work on constructive geometry is discussed in \cite{beeson-bsl}.
 
\subsection{Hilbert}
Hilbert's influential book \cite{hilbert1899} used the notion of betweenness and the axioms for 
betweenness studied by Pasch \cite{pasch1882}. 
Hilbert's theory was what would today be called ``second-order'', in that sets were freely used in the axioms.
Segments, for example, were defined as sets of two points, so by definition $AB = BA$ since the set $\{A,B\}$
does not depend on the order.  Of course, this is a trivial departure from first-order language; but Hilbert's
last two axioms,  Archimedes's axiom and the continuity axiom, are not expressible in a first-order geometrical 
theory.  On the other hand, lines and planes were regarded not as sets of points, but as (what today would be called)
first-order objects, so incidence was an undefined relation, not set-theoretic membership. 
At the time (1899) the concept of first-order language had not yet been developed, and set theory 
was still fairly new.   Congruence was treated by Hilbert as a binary relation on sets of two points, not as 
a 4-ary relation on points.   

Early geometers thought that the purpose of axioms was to set down the truth about space, so as to ensure
accurate and correct reasoning about the one true (or as we now would say, ``intended'')  model of those axioms.
Hilbert's book promoted the idea that axioms may have many models;  the axioms and deductions from them should make 
sense if we read ``tables, chairs, and beer mugs''  instead of ``points, lines, and planes.''   This is evident 
from the very first sentence of Hilbert's book:  

\begin{quote}{   Consider three distinct sets of objects. Let the objects of the first set be called {\em points}
$\ldots$; let the objects of the second set be called {\em lines} $\ldots$; let the objects of the third set be callled {\em planes}.}
\end{quote}

Hilbert defines segments as pairs of points (the endpoints), although lines are primitive objects.  On the other hand,
  a ray is the set of all points on the ray, and angles are sets consisting of two rays.  So an angle is a set of 
  sets of points.  
Hence technically Hilbert's theory, which is often described as second order, is at least third order.
(We say ``technically'', because it would not be difficult to reduce Hilbert's theory to an 
equivalent theory that would really be second order.)  

Hilbert's language has a congruence relation for segments, and a separate congruence relation for 
angles.   Hilbert's congruence axioms involve the concept of angles:  his fourth congruence axiom involves
``angle transport'' (constructing an angle on a given base equal to a given angle), and his fifth congruence axiom 
is the SAS triangle congruence principle.

Hilbert's Chapter VII discusses geometric constructions with a limited set of tools, a ``segment transporter''
and an ``angle transporter''.  These correspond to the betweenness and congruence axioms.  Hilbert does not discuss
the special cases of line-circle continuity and circle-circle continuity axioms that correspond to ruler and compass
constructions, despite the mention of ``compass'' in the section titles of Chapter VII.

Hilbert's geometry contained two axioms that go beyond first-order logic.   First, the axiom of Archimedes (which
requires the notion of natural number),
and second, an axiom of continuity, which (rephrased in modern terms) says that Dedekind cuts are filled on 
any line.   This axiom requires mentioning a set of points,  so Hilbert's theory with this axiom included is 
not a ``first-order theory''  in a language with variables only over points, lines, and circles.

\subsection{Tarski}
Later in 
the 20th century, when the concept of ``first-order theory'' was widely understood, Tarski formulated his 
theory of elementary geometry, in which Hilbert's axiom of continuity was replaced with an axiom schemata.
The set variable in the continuity axiom was replaced by any first-order formula.   Tarski proved that this 
theory (unlike number theory) is complete: every statement in the first-order language can be proved or refuted
from Tarski's axioms.   In addition to being a first-order theory, Tarski also made  other 
simplifications.  He realized that lines, angles, circles, segments, and rays could all be 
treated as auxiliary objects, merely enabling the construction of some new points from some given points.
Tarski's axioms are  stated using only variables for points.   He has only two relations: 
(non-strict) betweenness, ``$b$ is between $a$ and $c$'',  and ``equidistance'',  $E(a,b,c,d)$,
which means what Euclid expressed as  ``$ab$ is congruent to segment $cd$''.   We abbreviate
$E(a,b,c,d)$ informally as $ab=cd$.
 We have listed Tarski's axioms
for reference in \S\ref{section:axioms} of this paper, along with the axioms of our constructive version of Tarski geometry, 
adhering to the numbering of \cite{tarski-givant}, which has become standard.

\begin{figure}[h]
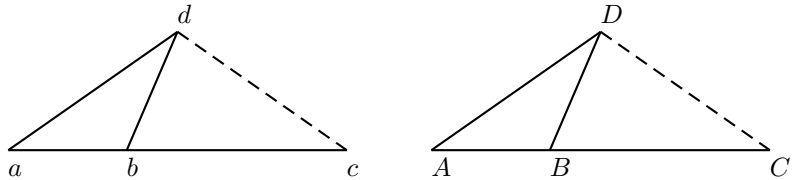
  
\caption{Tarski's 5-segment axiom.  $cd = CD$. 
\label{figure:TarskiFiveSegmentFigure}}
\TarskiFiveSegmentFigure
\end{figure}

Tarski replaced Hilbert's fourth and 
fifth congruence axioms (angle transport and SAS)  with an elegant axiom, known as the five-segment axiom.
This axiom is best understood not through its formal statement, but through   
Fig.~\ref{figure:TarskiFiveSegmentFigure}.  
The 5-segment 
axiom says that in Fig.~\ref{figure:TarskiFiveSegmentFigure}, the length of the dashed segment $cd$ is determined by the lengths 
of the other four segments in the left-hand triangle.  Formally, if the four solid 
segments in the first triangle are pairwise congruent to the corresponding segments in the second triangle, then the dashed segments
are also congruent.
  
\noindent
Tarski's 5-segment axiom is a thinly-disguised variant of the SAS criterion for triangle congruence.  To see 
this, refer to the figure.  The triangles we are to prove congruent are $dbc$ and
$DBC$.  We are given that $bc$ is congruent to $BC$ and $db$ is congruent to $DB$.   The congruence of angles $dbc$ and $DBC$ is expressed in
Tarski's axiom by the congruence of triangles $abd$ and $ABD$, whose sides are pairwise equal.   The conclusion,
that $cd$ is congruent to $CD$,  give the congruence of triangles $dbc$ and $DBC$.  In Chapter 11 of \cite{schwabhauser}, one can find a formal proof of the 
SAS criterion from the 5-segment axiom.  Borsuk-Szmielew also took this as an axiom (see \cite{borsuk-szmielew}, p. 81, Axiom C-5). 

An earlier version of Tarski's theory included as an axiom the 
``triangle construction theorem'',  which says that if we are given triangle $abc$, and segment $AB$ congruent to $ab$,
and a point $x$ not on $\Line(A,B)$, then we can construct a point $C$ on the same side of $\Line(A,B)$ as $x$ such that 
triangle $ABC$ is congruent to triangle $abc$.
It was later realized%
\footnote{
 Acording to \cite{tarski-givant}, Tarski included this principle  
 as an axiom in his first two published axiom sets, but then 
discovered in 1956-57 with the aid of Eva Kallin and Scott Taylor, that it was derivable; so he did not include it 
in \cite{tarski1959}. (See the footnote, p. 20 of \cite{tarski1959}.)  But Tarski did not publish the 
proof, and Borsuk-Szmielew take the principle as their Axiom C-7 \cite{borsuk-szmielew}. 
}
that this axiom is provable.  For example, one can drop a perpendicular from $c$ to $\Line(a,b)$, whose foot is the point $d$
on $\Line(a,b)$, and then find a corresponding point $D$ on $\Line(A,B)$, and then lay off $dc$ on the perpendicular to $\Line(A,B)$
at $D$ on the same side of $\Line(A,B)$ as $x$, ending in the desired point $C$.  Of course one must check that this construction 
can be done and proved correct on the basis of the other axioms.  But as it stands, this construction demands a case 
distinction about the order and possible identity of the points $d$, $a$, and $c$ on $\Line(a,b)$.  Hence, at least this
 proof of the triangle construction theorem from the axioms of Tarski's theory is non-constructive.

Tarski's early axiom systems also included axioms about betweenness and congruence that were later shown \cite{gupta1965}
to be superfluous.  The final version of this theory appeared in \cite{schwabhauser}; for the full history 
see \cite{tarski-givant}.%
\footnote{Note that the version mentioned in 
\cite{avigad2009} is not the final version used in \cite{schwabhauser}; 
inner transitivity for betweenness was eliminated in \cite{schwabhauser}.
} 
  The achievement of Szmielew and Gupta (who are mainly responsible for Part I of \cite{schwabhauser})
   is to develop a really minimal set of axioms for 
betweenness and congruence.%
\footnote{We would like to emphasize the important contributions of Gupta, 
which are important to the development in \cite{schwabhauser}, and are credited appropriately there, 
but without a careful study one might not realize how central Gupta's results were.   These results 
were apparently never published under Gupta's own name, and still languish in the Berkeley math library 
in his doctoral dissertation \cite{gupta1965}.  However, you can get that thesis and others from the 
ProQuest database, accessible from most university libraries.
}
  Hilbert's intuitive axioms about betweenness disappeared, leaving only 
the axiom $\neg \B(a,b,a)$ and the Pasch axiom and axioms to guarantee that congruence is an equivalence relation.  

\subsection{Strict {\em vs.}  non-strict betweenness and collinearity}
The (strict) betweenness relation is written
$\B(a,b,c)$.  We read this ``$b$ is between $a$ and $c$''.  The intended meaning is that 
that the three points are collinear and distinct, and $b$ is the middle one of the three.

Hilbert \cite{hilbert1899} and Greenberg \cite{greenberg} use strict betweenness, as we do.
Tarski \cite{tarski-givant} used non-strict betweenness.  They all used the same letter $\B$
for the betweenness relation, which is confusing.  For clarity we always use
$\B$ for strict betweenness, and introduce $\T(a,b,c)$ for non-strict betweenness.  Since $\T$
is Tarski's initial, and he used non-strict betweenness, that should be a memory aid.
 The two notions are interdefinable 
(even constructively): 

\begin{Definition} \label{defn:T} Non-strict betweenness is defined by 
$$ \T(a,b,c) := \neg(a \neq b \land b \neq c \land \neg \B(a,b,c))$$
\end{Definition}

\noindent
In the other direction,  $\B(a,b,c)$ can be defined as $$\T(a,b,c) \land a \neq b \land a \neq c.$$
If we express $\B(a,b,c)$ in terms of $\T(a,b,c)$, and then again express $\T$ in terms of $\B$,
we obtain a formula that is equivalent to the original using only axioms of classical propositional
logic. 
We mention this point to emphasize
that using these 
definitions, Tarski could have taken either strict or non-strict betweenness as primitive.  
In fact, to show that these definitions are ``inverses'' in the sense mentioned, we need
only intuitionistic logic plus the stability of equality and betweenness.  Therefore
(since we do accept stability of equality and betweenness), neither $\B$ nor $\T$
  is inherently more constructive than the other.

Why then did Tarski choose to use non-strict betweenness, when Hilbert had used strict betweenness?
Possibly, as suggested by \cite{tarski-givant},
 because this allowed him to both simplify the axioms, and reduce their number.  By using 
$\T$ instead of $\B$, the axioms cover various ``degenerate cases'',  when diagrams collapse onto 
lines, etc.  Some of these degenerate cases were useful.   From the point of view of 
constructivity, however, this is not desirable.  It renders Tarski's axioms {\em prima facie}
non-constructive (as we will show below).  Therefore the inclusion of degenerate cases in the axioms is something that 
will need to be eliminated in making a constructive version of Tarski's theories.  The 
same is true even if our only aim is to connect the axioms with ruler and compass constructions,
while retaining classical logic.

We next want to give a constructive definition of collinearity.
Classically we would define this as 
$\T(p,a,b) \lor \T(a,p,b) \lor \T(a,b,p)$.   That wouldn't work as a constructive definition 
of collinearity, because we have no way to decide in general which alternative might hold,
and the constructive meaning of disjunction would require it.  In other words,
 we can know that $p$ lies on the line determined by distinct points $a$ and $b$ 
without knowing its order relations with $a$ and $b$. But we can find a classically
equivalent (yet constructively valid) form by using the law that 
$ \neg\, \neg (P \lor Q)$
is equivalent to $\neg(\neg P \land \neg Q)$.   By that method we arrive at

\begin{Definition} \label{defn:collinearity}
$Col(a,b,p)$ is the formula expressing that $a$, $b$, and $p$ lie on a line.
$$   \neg ( \neg \T(p,a,b) \land \neg \T(a,p,b) \land \neg \T(a,b,p))$$
or equivalently, in terms of $\B$, 
$$   \neg ( \neg \B(p,a,b) \land \neg \B(a,p,b) \land \neg \B(a,b,p) \land a \neq p \land b \neq p \land a \neq b)$$ 
\end{Definition}

\noindent
Informally, we use the notation $\Line(a,b)$ to stand for the line determined by distinct
points $a$ and $b$, so ``$c$ lies on $\Line(a,b)$'' means $a \neq b \land Col(a,b,p)$.
Thus $Col(a,b,p)$ 
only expresses that $p$ lies on $\Line(a,b)$ if we also specify $a \neq b$.  We do not put 
the condition $a\neq b$ into the definition of $Col(a,b,p)$ for two reasons:  it would destroy
the symmetry between the three arguments, and more important, it would cause confusion in 
comparing our work with the standard reference for Tarski's theories, namely \cite{schwabhauser}.

\subsection{Pasch's axiom}

\begin{figure}[h]
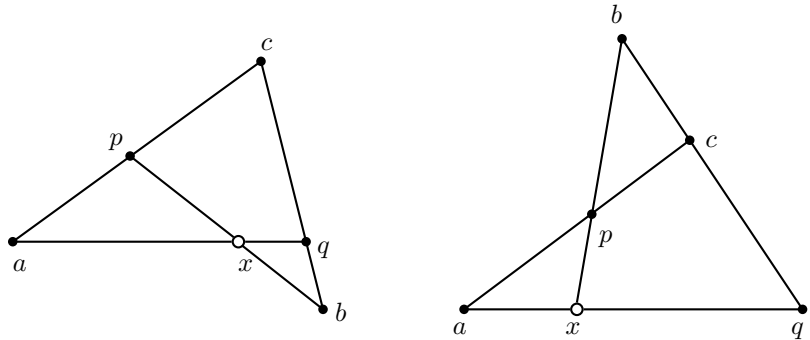
  
\caption{ Inner Pasch (left) and outer Pasch (right).  Line $pb$ meets triangle $acq$ in one side.
 The open circles show the points asserted to exist on the other side. }
\label{figure:InnerOuterPaschFigure}
\InnerOuterPaschFigure
\end{figure}

  Hilbert's fourth
betweenness axiom is often known as Pasch's axiom, because it was first studied
by Pasch in 1882 \cite{pasch1882}.  It says that if line $L$ meets (the interior of) side $AB$ of triangle $ABC$
then it meets (the interior of) side $AC$ or side $BC$ as well. 
But Tarski considered instead, two restricted versions of Pasch's axioms known as ``inner Pasch''
and ``outer Pasch'', illustrated in   Fig.~\ref{figure:InnerOuterPaschFigure}.

Outer Pasch was an axiom (instead of, not in addition to, inner Pasch) 
in versions of Tarski's theories until 1965, when it was proved from inner Pasch in Gupta's thesis \cite{gupta1965}, Theorem 3.70, or Satz 9.6 in  \cite{schwabhauser}.%
\footnote{But apparently, judging from footnote 4 on p. 191 of \cite{tarski-givant}, Tarski knew as early as 
1956-57 that outer Pasch implies inner Pasch; in that footnote Tarski argues against replacing outer Pasch with 
inner Pasch as an axiom, as Szmielew and Schwabh\"auser chose to do. Also on p. 196 of \cite{tarski-givant},
Tarski attributes the idea for the proof of inner Pasch from outer Pasch to specific other people; the history 
is too long to review here, but he credits only Gupta with the derivation of outer Pasch from inner Pasch.
}  
Outer Pasch appears as Satz 9.6 in \cite{schwabhauser}.
The proof  given in \cite{schwabhauser}, applied to the formulation of outer Pasch
with strict betweenness,  is constructive.  The proof is complicated, however, in that 
it depends on the ability to drop a perpendicular to a line from a point not on the line.
As we shall discuss extensively, proving the existence of such perpendiculars is problematic:
easy constructions require either line-circle continuity or the parallel axiom, and only 
with Gupta's complicated proof do we have dropped perpendiculars without either of those 
assumptions. In \S\ref{section:doublenegation} we shall show how to use the double-negation 
interpretation (a metamathematical tool developed by G\"odel) to ensure the constructivity
of the proof of outer Pasch without checking it line by line.  The unsatisfied reader 
has the choice to either 
  check the proof directly, or look ahead to \S\ref{section:doublenegation}, which 
 does not depend on the intervening material.
 We will use outer Pasch as needed.

After Gupta proved outer Pasch from inner, Szmielew chose to take inner Pasch as an axiom instead of outer Pasch, although
a footnote in \cite{tarski-givant} shows that Tarski disagreed with that choice (on grounds 
of how easy or difficult it is to deduce other things).
 Gupta's thesis also contains a proof that outer Pasch implies inner Pasch.

It is not completely clear why Tarski wanted to restrict Pasch's axiom in the first place,
but two good reasons come to mind.  First,  
the restricted forms are valid even in three-dimensional space, so they 
do not make an implicit dimensional assertion, as the unrestricted Pasch axiom does (it fails in three-space).  Second, there is the simpler logical form of inner (or outer) Pasch: unrestricted
Pasch needs either a disjunction, or a universal quantifier in the hypothesis, so the condition 
to be satisfied by the point whose existence is asserted is not quantifier-free and disjunction-free,
as it is with inner and outer Pasch.   This simplicity of logical form is important for our purposes
in constructive geometry, but for Tarski it may just have been a matter of ``elegance.''

\subsection{Sides of a line} \label{section:sameside} 
The notions of ``same side'' and ``opposite side'' of a line will be needed below, 
and are also of interest in comparing Hilbert's and Tarski's geometries.  
One of Hilbert's axioms was the {\em plane separation axiom}, according to which a line 
separates a plane into (exactly) two regions.  
Two points $a$ and $b$ not on line $L$ are on opposite sides of $L$ if $a \neq b$ and 
there is a point of $L$ between $a$ and $b$, i.e., the segment $ab$ meets $L$.

\begin{Definition} \label{definition:oppositesides}
\begin{eqnarray*}
\hskip -0.5cm &&\OppositeSide(a,b,L) :=  \exists x\,( on(x,L) \land \B(a,x,b) )
\end{eqnarray*}
\end{Definition}
Of course, in Tarski geometry, we cannot mention lines directly, so $L$ has to be 
replaced by two distinct points, yielding a 4-ary relation.

The definition of being on the same side is less straightforward.  Hilbert's definition of 
$\SameSide(a,b,L)$ was that segment $ab$ does not meet $L$.   That involves a universal 
quantifier:
$$ \forall x \neg (\B(a,x,b) \land on(x,L)).$$
One can get an existential quantifier instead of a universal quantifier by using Tarski's definition,
illustrated in Fig.~\ref{figure:SameSideFigure}:

\begin{figure}[h]
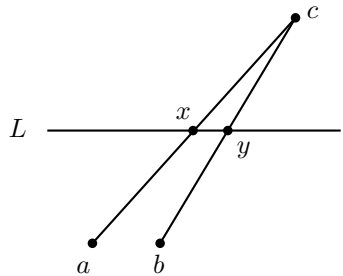
  
\caption{Tarski's definition: $a$ and $b$ are on the same side of line $L$, as witnessed by point $c$
on the other side.} 
\label{figure:SameSideFigure}
\SameSideFigure
\end{figure} 

\begin{Definition} \label{definition:sameside}
$a$ and $b$ are on the same side of $L$ if there is 
some $c$ such that both $a$ and $b$ are on the opposite side of $L$ from $c$.  Formally:
\begin{eqnarray*}
 \SameSide(a,b,L) :=\exists c,x,y\,(\B(a,x,c) \land \B(b,y,c) \land on(x,L) \land on(y,L))
\end{eqnarray*}
\end{Definition}
 Another advantage of 
this definition is that it works in more than two dimensions.
It can be proved equivalent to  
 Hilbert's definition above (as is discussed in section \ref{section:TarskiProvesHilbert} below),
if all points are restricted to lie in the same plane; but in Tarski's geometry, planes are 
a defined concept (using $SameSide$) rather than a primitive concept.
 
Hilbert took it as axiomatic that a line divides a plane into two regions.  In Tarski's 
system this becomes a fairly difficult theorem:

\begin{Theorem}[Plane separation theorem] \label{theorem:planeseparation}
If $p$ and $q$ are on the same side of line $L$,
and $p$ and $r$ are on opposite sides of $L$, then $q$ and $r$ are also on opposite sides
of $L$.  Formally,
$$\SameSide(a,b,L) \land \OppositeSide(a,c,L) \implies \OppositeSide(b,c,L)$$
is provable in neutral constructive geometry (i.e., without using the parallel axiom).
\end{Theorem}

\noindent{\em Proof}.  This is proved in Gupta \cite{gupta1965}, and also as Satz 9.8 of \cite{schwabhauser}.  The proof follows fairly easily 
from outer Pasch and the definition of $\SameSide$, and occurs in \cite{schwabhauser}
right after the proof of outer Pasch.  The proof (from outer Pasch) 
is completely and unproblematically 
constructive. 

 \subsection{The parallel axiom according to Hilbert and Tarski}
As is well-known, there are many propositions equivalent to the parallel postulate in classical
geometry.  The main point of \cite{beeson-bsl} is to establish which of these
 versions of the parallel postulate
are   equivalent in constructive geometry, and which are not.
 Hilbert's parallel axiom (Axiom IV, p. 25 of \cite{hilbert1899}) is the version we call Playfair's Axiom, introduced by Playfair in 1729:
There cannot be more than one parallel to a given line through a point not on the line.  
Tarski's axiom A10 as published in \cite{schwabhauser} is a more complicated statement, classically equivalent. Specifically, it says that if $p$ is in the (closed) interior of angle $\alpha$, then there 
exist points $x$ and $y$ on the sides of $\alpha$ such that $\T(x,p,yt)$.  Of course, one cannot 
mention ``interior of angle $\alpha$'' directly, so the formulation in Tarski's language is a bit 
more complex.   Szmielew's manuscript, on which Part I of \cite{schwabhauser}  is based, took instead
the ``triangle circumscription principle'', which says that for every three non-collinear  
points $a$, $b$, $c$,  there exists a point $d$ equidistant from all three (thus $d$ is the 
center of a circle passing through $a$, $b$, and $c$, thus circumscribing triangle $abc$).%
\footnote{The triangle circumscription principle is equivalent (with classical logic)
to Euclid's parallel axiom.  Euclid IV.5 proves the triangle circumscription principle;
the converse implication was first proved by Farkas Bolyai, father of Janos Bolyai,
who thought he had proved Euclid's parallel postulate, but had actually assumed 
the triangle circumscription principle. See \cite{greenberg}, pp. 229--30 and p. 240.
}

In \cite{beeson-bsl}, we considered the parallel axiom from the constructive point of view,
and gave a points-only version of Euclid's parallel postulate, called ``Euclid 5'', as well as 
a stronger version called the ``strong parallel postulate.''  These turned out to be constructively
equivalent, though the proof requires the prior development of considerable ``machinery'' based
on Euclid 5.  We also showed that the triangle circumscription principle is equivalent to 
the strong parallel postulate, and hence to Euclid 5.  In this paper (Theorems~\ref{theorem:8.2}
and \ref{theorem:8.3}), we show that Tarski's 
parallel axiom is equivalent to Euclid 5, too.  Hence all the versions of the parallel postulate
that make an existential assertion turn out to be equivalent.    

For the reason of simplicity, we follow Szmielew
in using the triangle circumscription principle as the parallel axiom in Tarski's theories.%
\footnote{The change in the parallel axiom was apparently one of the ``inessential changes''
Schwabh\"auser introduced in publishing Szmielew's work.  I have not seen Szmielew's manuscript,
but base what I say about it here on \cite{tarski-givant}, page 190.} 
The center of the circumscribed circle $abc$ can be constructed with ruler and compass as the 
intersection point of the perpendicular bisectors of $ab$ and $bc$; the point of the axiom 
is that these lines do indeed meet (which without some form of the parallel axiom, they 
cannot be proved to do). The axiom lends itself well to a points-only theory, since it does 
not actually mention circles.  It merely says there is a point equidistant from the three 
given points. 

Tarski and Givant wrote a letter to Schwabh\"auser ``around 1978'', which was published 
in 1998 \cite{tarski-givant} and has served, in the absence of an English translation of 
\cite{schwabhauser}, as a common reference for Tarski's axioms and their history.  
The letter mentions   equivalent versions of the parallel axiom:  the two mentioned above and 
a ``Third version of the parallel axiom'',  which says that if one connects the midpoints of 
two sides of a triangle, the connecting segment is congruent to half the third side.   In spite
of the name ``Third version of the parallel axiom'', the letter makes no claim that the different
versions are equivalent (in any theory at all).  One has to be careful when speaking about 
``versions of the parallel postulate.''  According to \cite{szmielew1959}, p.~51, any statement that 
holds in Euclidean geometry but not in the standard hyperbolic plane is (classically) 
equivalent to Euclid's parallel postulate in Tarki's geometry with full first-order continuity axioms
(Axiom (A11) of \cite{tarski-givant}).
In other words, there are only two complete extensions of neutral geometry with full continuity.
But no such thing is true in the theories considered here, which have only line-circle and 
circle-circle continuity (one as an axiom, and one as a theorem).

Indeed, the ``third version'' mentioned above (which we here call $M$, for ``midline'')  
is not equivalent to the parallel postulate  (in neutral geometry with line-circle and circle-circle 
continuity), but instead to the weaker assertion that the sum 
of the angles of every triangle is equal to two right angles.    The non-equivalence 
with the parallel axiom is proved as follows:

\begin{Theorem}\label{theorem:Munprovable}
 No quantifier-free statement can be equivalent to the parallel axiom
in neutral geometry  with  circle-circle and line-circle continuity.
\end{Theorem} 

\noindent{\em Proof}.  We give a model of neutral geometry in which $M$ (or any quantifier-free formula
that is provable with the aid of the parallel axiom) holds, but the parallel axiom fails.
Let $\F$ be a non-Archimedean Euclidean field, and let $\K$ 
be the finitely bounded elements of $\F$,  i.e.,elements
between $-n$ and $n$ for some integer $n$.  (Then $\K$ is Archimedean but is not a field,
because it contains ``infinitesimal'' elements whose inverses are in $\F$ but not in $\K$.)
 The model is $\K^2$.  
This model is due to Max Dehn, and is  described in Example 18.4.3  and Exercise 18.4
of \cite{hartshorne},
where it is stated that $\K^2$ is a Hilbert plane, and also satisfies line-circle and 
circle-circle continuity, since the intersection points with finitely bounded circles have
finitely bounded coordinates.

Since $\F^2$ is a 
model of geometry including the parallel axiom, $M$ holds there, and since $M$ is quantifier 
free, it holds also in $\K$.  Yet,  $\K$ is not a Euclidean plane; let $L$ be the $x$-axis and let $t$ be an infinitesimal. There are many lines through $(0,1)$  that are  parallel to $L$  in $\K$ (all but one of them are restrictions to $P$ of lines in $F^2$  that meet the $x$-axis at some non finitely bounded point). 
That completes the proof.
\medskip

\noindent{\em Discussion}. As remarked above, it follows from Szmielew's work \cite{szmielew1959}, p.~51, that $M$ is equivalent to the parallel axiom in Tarki's geometry with classical 
logic and the full first-order continuity axiom (A11).   The question then arises, how exactly can we 
use elementary continuity to prove Euclid 5 from $M$?  Here is a proof:  Assume, 
for proof by contradiction, the negation of Euclid 5.  Then, by elementary continuity, limiting 
parallels exist (see \cite{greenberg}, p. 261). Then Aristotle's axiom holds, as proved in \cite{hartshorne}, Prop.~40.8, p.~380.
But $M$ plus Aristotle's
axiom implies Euclid 5 (see \cite{greenberg}, p.~220), contradiction, QED. 

 This proof is interesting 
because it uses quite a bit of machinery from hyperbolic geometry to prove a result that, on the face 
of it, has nothing to do with hyperbolic geometry.   That is, of course, also true of the proof
via Szmielew's metamathematics.   Note that a non-quantifier-free instance of 
elementary continuity is needed to get the existence of limiting parallels directly; in the presence
of Aristotle's axiom, line-circle continuity suffices (see \cite{greenberg}, p. 258), but Aristotle's
axiom does not hold in $P$.  Finally, the proof of Theorem~\ref{theorem:Munprovable}
shows that the use of 
a non-quantifier-free instance of continuity is essential, since quantifier-free instances will 
hold in Dehn's model (just like line-circle and circle-circle continuity).

\subsection{Interpreting Hilbert in Tarski} \label{section:TarskiProvesHilbert}
The fundamental results about betweenness discussed in section~\ref{section:betweenness},
 along with many pages of further work, enabled
Szmielew to prove (interpretations of) Hilbert's axioms in Tarski's theory.  Neither she 
nor her (posthumous) co-authors pointed this out explicitly in \cite{schwabhauser}, but
it is not difficult to find each of Hilbert's axioms among the theorems of \cite{schwabhauser}
(this has been done explicitly, with computer-checked proofs, in \cite{narboux2012}).  Here we illustrate by 
comparing Hilbert's betweenness axioms to Tarski's:  Both have symmetry.
Hilbert's II,3  is ``Of any three points on a line there exists no more than one that 
lies between the other two.''  We render that formally as
$$ a\neq b \land b\neq  c \land a\neq c \land \B(a,b,c) \implies \neg \B(b,a,c) \land \neg \B(a,c,b).$$ 
This can be proved from Tarski's axioms as follows:  suppose $a$, $b$, and $c$
are distinct, and $\B(a,b,c)$.  
Then $\neg\B(b,a,c)$, since if $\B(b,a,c)$ then $\B(a,b,a)$, by inner transitivity and symmetry.
(See \S\ref{section:axioms} for the formulas mentioned by name here.)
Also, $\neg \B(a,c,b)$, since if $\B(a,b,c)$ and $\B(a,c,b)$, then $\B(a,b,a)$ by inner transitivity
and symmetry.  

Hilbert has a ``density'' axiom
   (between two distinct points there is a third).  This is listed as (A22) in \cite{tarski-givant}, but was never an axiom of Tarski's theory.  
Density can be proved classically even without line-circle or circle-circle 
continuity:  Gupta (\cite{gupta1965}, or \cite{schwabhauser}, Satz 8.22) showed that the midpoint of a segment can be constructed without continuity.  It is also possible to give a very short
direct proof of the density lemma:

\begin{Lemma}[Density]\label{lemma:density} 
Given distinct points $a$ and $c$, and point $p$ not collinear with $a$ and $c$,
there exists a point $b$ with $\B(a,b,c)$.
\end{Lemma}

\noindent{\em Remark}.   
The proof requires only an axiom stating that the order is unending (here 
the extension axiom supplies that need),
the inner form of the Pasch axiom, and the existence of a point not on a given line.
 This theorem has been known for sixty years \cite{lumiste1964,lumiste2005,pambuccian-pasch},
as it turns out, and rediscovered by Ben Richert in 2014, whose proof is given here.
 \medskip

\noindent{\em Proof}.  Extend $ap$ by $ac$ to point $r$, and extend $rc$ by $ac$ to point $s$.
Apply inner Pasch to $scrap$. The result is point $b$ with $\B(a,b,c)$  (and $\B(p,b,s)$,
but that is irrelevant).  That completes the proof.
\medskip

As discussed above, 
one can prove in Tarski's system (using the dimension axioms) that Hilbert's and 
Tarski's definitions of $\SameSide$ coincide; and Hilbert's plane separation axiom 
becomes a theorem in Tarski's system.

Hilbert's theory has variables for angles; but in Tarski's theory, angles are given by 
ordered triples of non-collinear points, and the theory of congruence and ordering of angles has to be 
developed, somewhat laboriously, but along quite predictable lines, carried out in \cite{schwabhauser}.
Some details of the Tarskian theory of angles are discussed in \S\ref{section:angles} below;
the upshot is that (even constructively) one can 
construct a 
conservative extension of Tarski geometry that has variables for angles and directly supports the 
kind of arguments one finds in Euclid.

\begin{figure}[ht]
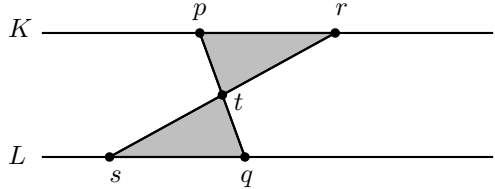
   
\caption{Transversal $pq$ makes alternate interior angles equal with $L$ and $K$,  if $pt=tq$ and 
$rt=st$.}
\label{figure:AlternateInteriorAnglesFigure}
\hskip 2.5cm
\AlternateInteriorAnglesFigure
\end{figure}

It is sometimes possible to reduce theorems about angles directly; in particular it is not 
necessary to develop the theory of angle ordering to state Euclid's parallel postulate.
In Fig.~\ref{figure:AlternateInteriorAnglesFigure}, we show how to 
translate the concept ``equal alternating interior angles'' into Tarski's language.

\smallskip    

\section{Tarski's theory of straightedge and compass geometry}
Tarski's theory is ``elementary'' only in the sense that it is first-order.  It still goes far beyond Euclid.%
\footnote{It is confusing that in axiomatic geometry,  ``elementary''
sometimes refers to the elementary constructions, and sometimes to the full first-order theory of Tarski.   In this 
paper we shall not refer again to the full first-order theory.  
}
To capture Euclid's geometry, Tarski considered the subtheory in which the continuity axiom is replaced by 
``segment-circle continuity''.   This axiom asserts the existence of the intersection 
points of a line segment and a circle, if  some point on the segment lies inside the circle and
some point on the segment lies outside the circle.

It is this theory that we refer to in the title of this paper as ``Tarski's geometry''.
In the section title, we mention ``straightedge and compass'';  but henceforth we use the 
more common terminology ``ruler and compass'', with the same meaning.

 \subsection{Line-circle continuity} 
We now formulate the axiom of line-circle continuity.
This tells us when a line and a circle intersect--namely, when there is a point 
on the line closer (or equally close) to the center than the radius of the circle.%
\footnote{Note that in spite of the use of the word ``circle'' the axiom, in the form that 
only asserts the existence of an intersection point, is valid in $n$-dimensional
Euclidean space, where it refers to the intersections of lines and spheres.}
  But we have not 
defined inequalities for segments yet, so the formal statement is a bit more complex.  Moreover, we 
have to include the case of a degenerate circle or a line tangent to a circle, without making a case distinction.%
\footnote{Avigad {\em et. al.} count only {\em transverse} intersection, not tangential intersection, as \\ ``intersection.''}
Therefore we must find a way to express ``$p$ is inside the closed circle with center $a$ passing 
through $y$''.  For that it suffices that there should 
be some $x$  non-strictly between $a$ and $y$ such that $ax=ap$.   
Since this will appear in the antecedent of the 
axiom, the ``some $x$'' will not involve an existential quantifier.  

\begin{Definition}  \label{defn:lessthan}
$ab < cd$ (or $cd > ab$)  means $\exists x(\B(c,x,d) \land ax = ab$)).

\noindent
$ab \le cd$ (or $cd \ge ab$)  means $\exists x(\T(c,x,d) \land ax = ab)$, where $\T$ is 
non-strict betweenness.
\end{Definition}
 
\begin{Definition} \label{defn:inside}
 Let $C$ be a circle with center $a$.  Then point $p$ is {\bf strictly inside} $C$
means there exists a point $b$ on $C$ such that $ap < ab$, and $p$ is {\bf inside} $C$, or {\bf non-strictly inside} $C$, means $ap \le ab$.  
\end{Definition}

\begin{figure}[h]
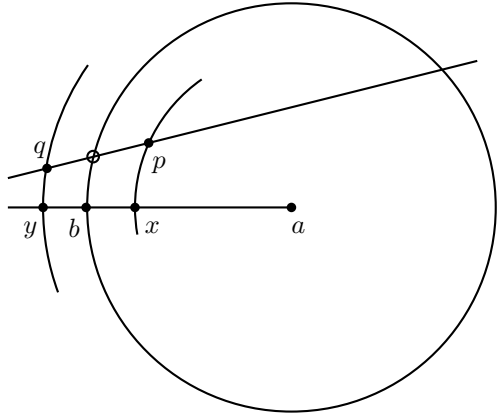
  
\caption{Segment-circle continuity. $p$ is inside the circle, $q$ is outside,
 so $L$ meets the segment $pq$. \label{figure:SegmentCircleContinuityFigure}}
\psset{unit=2.6cm}
\SegmentCircleContinuityFigure
\psset{unit=3cm}
\end{figure}

Replacing `$<$' by `$>$', we obtain the definition of {\bf outside}. 
The version of line-circle continuity given in \cite{tarski1959} is better
described as ``segment-circle'' continuity:
$$ax = ap \land \T(a,x,b)\land \T(a,b,y) \land ay=aq   \implies
 \exists z\,( \T(p,z,q) \land az = ab)$$
 This axiom says that if $p$ is inside circle $C$ and $q$ is outside $C$, then segment $pq$ meets
circle $C$.  See Fig.~\ref{figure:SegmentCircleContinuityFigure}.

\smallskip
\medskip

\begin{figure}[h]
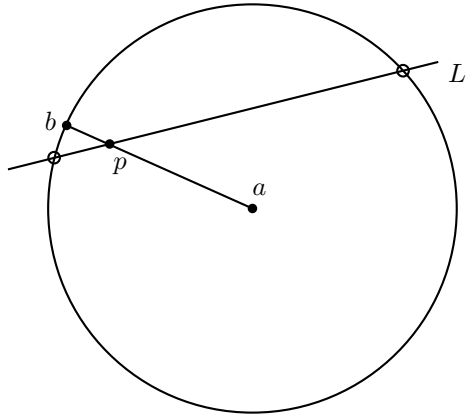
  
\caption{Line-circle continuity.  $p$ is inside the circle, so $L$ meets the circle. \label{figure:LineCircleContinuityFigureTwo}}
\psset{unit=2.6cm}
\LineCircleContinuityFigureTwo
\psset{unit=3cm}
\end{figure}

One may also consider a geometrically simpler formulation
of line-circle continuity:  if line $L = \Line(u,v)$ has a point $p$ inside circle $C$,  then there is a point that lies on both $L$ and $C$.
   See Fig.~\ref{figure:LineCircleContinuityFigureTwo}.

\smallskip
\medskip

We consider two versions of this axiom.   The weaker version (one-point line-circle)
only asserts the existence of one intersection point.  The stronger version (two-point
line-circle) adds the extra assertion that if $p\neq b$ (i.e., $p$ is strictly inside 
the circle) then there are two distinct intersection points.

Here are the formal expressions of these axioms.  (The formula $Col$, for collinearity, 
is given in 
Definition~\ref{defn:collinearity}.)
\begin{eqnarray*}
Col(u,v,p) \land u \neq v \land  \T(a,p,b)   \implies  &  \mbox{\qquad (one-point line-circle)}\\
   \exists  z\,(Col(u,v,z) \land az = ab) &
\end{eqnarray*}

$$Col(u,v,p) \land u \neq v \land  \T(a,p,b)    \implies   \mbox{\qquad (two-point line-circle)}$$
$$  \exists  y,z\,( az = ab \land ay = ab \land \T(y,p,z) \land (p \neq a \implies y \neq z)) $$ 
 
Classically, we could take a shorter version of two-point line-circle:
\begin{eqnarray*}
Col(u,v,p) \land u \neq v \land  \T(a,p,b) \land p\neq a  \implies  &  \mbox{\quad (classical two-point line-circle)}\\
   \exists  y,z\,( ay= ab \land az = ab \land \B(y,p,z)) &
\end{eqnarray*}
This is classically equivalent to two-point line circle, since the case when $p=a$ is trivial; 
but constructively, we cannot make a case distinction whether $p=a$ or not.  The longer form 
is necessary for a constructive version.

The equivalence of these three continuity axioms, relative 
to the other axioms of Tarski geometry, is not at all obvious (even with classical logic),
because
\smallskip

(i) in order to show line-circle implies segment-circle, we need to 
construct points on the line outside the circle, 
which requires the triangle inequality.  In turn the triangle inequality 
requires perpendiculars.
\smallskip

(ii) in order to show one-point line-circle implies two-point line-circle,
we need to construct the second point somehow.  To do that, we need to be able 
to construct a perpendicular to the line through the center.  Classically this 
requires a dropped perpendicular from the center to the line (as the case when 
the center is on the line is trivial);  constructively it requires even more,
a ``uniform perpendicular'' construction that works without a case distinction.
But even the former is difficult.

Since two-point line-circle continuity corresponds directly to the uses made 
(implicitly) of line-circle continuity in Euclid, we adopt it as an axiom of 
our constructive version(s) of Tarski's theory.  We shall show eventually that 
all three versions are in fact equivalent,  using the other axioms of Tarski's 
theory (and not even using any form of the parallel axiom).  But this proof
rests on the work of Gupta \cite{gupta1965}, which we will also discuss below.

\subsection{Intersections of circles}
\noindent
We next give the principle known as circle--circle continuity.  
It should say that if point $p$ on circle $K$ lies (non-strictly) inside circle $C$, and point $q$
on $K$ lies (non-strictly) outside $C$, then both 
intersection points of the circles are defined.    This principle would be taken as an axiom,
except that it turns out to be derivable from line-circle continuity, so it is not necessary
as an axiom.  This implication will be 
proved and discussed fully in \S~\ref{section:linecircle}, where it will be shown 
to be true also with intuitionistic logic.%
\footnote{It is also true that circle-circle continuity implies line-circle
continuity.  See for example 
\cite{greenberg}, p.~201.
Proofs of the equivalence of line-circle and circle-circle
continuity using Hilbert's axioms (with no continuity and without even the 
parallel axiom) were found by Strommer \cite{strommer1973}.  Since these axioms
are derivable from (A1)-(A9), as shown by Gupta and Szmielew \cite{schwabhauser,narboux2012},
the equivalence can be proved in (A1)-(A9) (with classical logic).  We have not 
studied the constructivity of Strommer's proof.}

\begin{figure}[h]
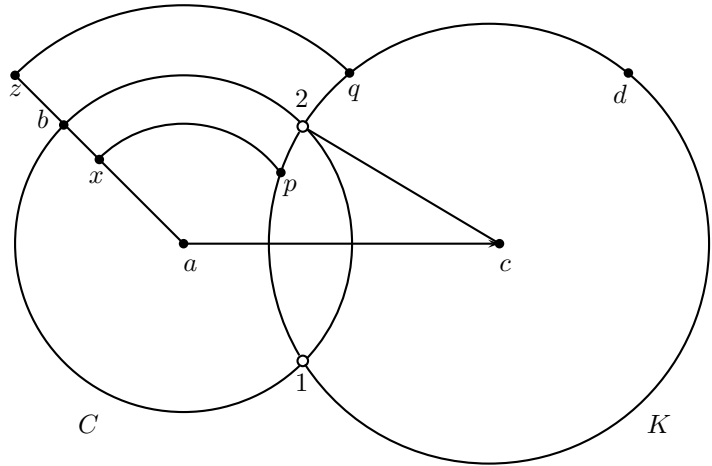
  
\caption{Circle-circle continuity. $p$ is inside $C$ and $q$ is outside $C$, as witnessed by $x$  and $z$, so the 
intersection points 1 and 2 exist. }
\label{figure:CircleCircleContinuityFigure}
\psset{unit=2.8cm}
\CircleCircleContinuityFigure
\psset{unit=3cm}
\end{figure}

In Tarski's points-only language, circles are given by specifying the center and a point
through which the circle passes. Informally, we write $Circle(a,b)$ for the circle with 
center $a$ passing through point $b$.  The case $a=b$ (a ``degenerate circle'') is allowed.
To express that $p$ is non-strictly inside (or outside) a circle $C$,
we use the same technique as just above.  Namely, $p$ is inside $C$ if $ap=ax$ for some $x$
non-strictly between $a$ and a point $b$ on $C$.   The situation is illustrated in Fig.~\ref{figure:CircleCircleContinuityFigure}, where circle $C$ is given by center $a$ and point $b$,
and circle $K$ is given by center $c$ and point $d$.
\medskip

We want this principle
 to apply even to degenerate circles, and to points that are on $C$ rather than strictly inside, so we must use $\T$ rather than $\B$ to allow
 $x=y$ or $y=z$, and we must even allow $a=x=b=z$. But we do not want it to apply when $C$ and 
 $K$ are the same circle, for then,  although classically plenty of intersection points exist,
 no special point is singled out by a ``construction.''
 
 \smallskip

In order to express this axiom using point variables only, we think of $K$ as $\Circle(c,d)$
and $C$ as $\Circle(a,b)$.  Then the axiom becomes
\medskip

\axioms
$  \ ap=ax \land aq=az \land\ cp = cd \ \land cq=cd\ \land \T(a,x,b)\  \land $ & \hskip-0.3cm$\T(a,b,z) \land a \neq c$    \\
$  \implies \exists z_1,z_2\,(cz_1=cd \land az_1 = ab \land  cz_2=cd \land az_2 = ab) $  & (circle-circle) 
\endaxioms
\smallskip

\noindent
The use of non-strict betweenness $\T$ allows for the cases when the circles are tangent (either 
exterior or interior tangency).

It is not necessary to assert the existence of two distinct intersection points when $p$ is 
strictly inside $C$, since the second intersection point can be constructed as the 
reflection of the first in the line connecting the two centers.  Then, using the 
plane separation theorem, one can prove the existence of an intersection point on a 
given side of the line connecting the centers.

\section{Listing of axioms for reference} \label{section:axioms}
In this section, we list the axioms of all theories used in this paper.
It is intended for reference, rather than as part of a first reading,
but we follow Tarski's example and a referee's advice in putting the axioms, in the 
referee's phrase, ``front and center.''

In the following, $ab=cd$ abbreviates $E(a,b,c,d)$, and $\T(a,b,c)$ is 
non-strict betweenness, while $\B(a,b,c)$ is strict betweenness.  

\subsection{Classical two-dimensional Tarski geometry} 
We give the version preferred by 
Szmielew.  The version in \cite{schwabhauser} has (the classically equivalent) 
(A10) instead of 
(A10${}_3$).   We also give the Skolemized versions here.  $Col(a,b,c)$ (collinearity)
is an abbreviation for $\T(a,b,c) \lor \T(b,c,a) \lor \T(c,a,b)$.
 
\smallskip
\newaxioms{2}
$ab = ba$ & (A1) Reflexivity of equidistance \\
$ab = pq \land ab = rs \implies pq=rs$ & (A2) Transitivity of equidistance \\
$ab = cc \implies a=b$ & (A3) Identity of equidistance \\
$\exists x\, (\T(q,a,x) \land ax = bc)$ & (A4) Segment extension \\
$ \T(q,a,ext(q,a,b,c)) \land E(q,ext(q,a,b,c),b,c)$ & (A4), Skolemized  \\
$(a \neq b \land \T(a,b,c) \land \T(A,B,C) \land ab=AB \land bc = BC$ & \\
\qquad $ ad=AD \land bd=BD) \implies cd = CD $ & (A5) Five-segment axiom \\
$ \T(a,b,a) \implies a=b$ & (A6) Identity for betweenness \\
$ \T(a,p,c) \land \T(b,q,c) \implies \exists x\,(\T(p,x,b) \land \T(q,x,a))$ & (A7) inner Pasch \\
$ \T(a,p,c) \land \T(b,q,c) \implies$ &\\
\quad $ \T(p,ip(a,p,c,b,q),b) \land  \T(q,ip(a,p,c,b,q),a)$ & (A7), Skolemized \\
$ \exists a,b,c\, \neg Col(a,b,c)  $ & (A8), lower dimension \\
$ \neg Col(\alpha,\beta,\gamma) $ & (A8), Skolemized \\
$ pa=pb \land qa=qb \land ra=rb \implies Col(a,b,c)$ &   (A9), upper dimension\\
$\neg\, Col(a,b,c)\implies  \exists x\,(ax = bx \land ax = cx)$ & (A10${}_3$), triangle circumscription \\
$ax = ap \land \T(a,x,b)\land \T(a,b,y) \land ay=aq   \implies $ & segment-circle continuity\\
$ \exists z\,( \T(p,z,q) \land az = ab)$&
 \endaxioms
\medskip

In the Skolemized version of the triangle circumscription principle, $x$ is given by $center(a,b,c)$.
We make no use of a Skolemized version of segment-circle, so we do not give one.
\medskip

\subsection{Intuitionistic Tarski geometry} \label{section:axioms2}
  This theory takes $\B$ as 
primitive rather than $\T$, so $\T(a,b,c)$ is an abbreviation for 
$\neg(a \neq b \land b \neq c \land \neg \B(a,b,c))$, and $Col(a,b,c)$ is 
an abbreviation for 
$$ a \neq b \land \neg ( \neg \B(p,a,b) \land \neg \B(a,p,b) \land \neg \B(a,b,p) \land a \neq p \land b \neq p),$$ 
which is equivalent to the double negation of the classical definition of $Col(a,b,c)$ together 
with $a \neq b$.  In other words, $Col(a,b,c)$ says $c$ lies on $\Line(a,b)$.
The axioms (A1)-(A3) and (A5) are unchanged, except that now $\T$ is defined in terms of $\B$.
It is inessential whether $\T$ or $\B$ is taken as primitive.

The differences between classical and intuitionistic Tarski geometry are
\begin{itemize}
\item (A4): Only non-null segments can be extended.
\item Axiom (A6) becomes $\neg \B(a,b,a)$.  
\item inner Pasch (A7): The hypothesis $\T(a,p,c)$ is changed to $\B(a,p,c)$,
 and the hypothesis $\neg Col(a,b,c)$ is added, and 
 although $\T$ is strengthened to $\B$ in the conclusion.  But the hypothesis
 $\T(b,q,c)$ is not changed.
\item Symmetry and inner transitivity of betweenness (A14) and (A15) are added.
\item A negative
formula is used for collinearity in the dimension axioms and the triangle circumscription principle.
\item In line-circle continuity, the two points $p$ and $q$ determining the line are assumed 
to be unequal, and we use two-point line-circle continuity instead of segment-circle.
\item We use intuitionistic logic and add the stability axioms.
\end{itemize}
\medskip

Intuitionistic Tarski geometry plus classical logic is called ``continuous Tarski geometry'';
we can have continuous Tarski geometry with or without Skolem functions.
The changed axioms are as follows:
\medskip

\newaxioms{2}
$ q \neq a \implies \exists x\, (\T(q,a,x) \land ax = bc)$ & (A4-i) Segment extension \\
$ q\neq a \implies \T(q,a,ext(q,a,b,c)) \land E(q,ext(q,a,b,c),b,c)$ & (A4-i), Skolemized  \\
$ \neg \B(a,b,a)$ & (A6-i) \\
$ \B(a,p,c) \land \T(b,q,c)   \land \neg\, Col(a,b,c) \implies$&\\
\qquad $\exists x\,(\B(p,x,b) \land \B(q,x,a))$  & (A7-i)  strict  inner Pasch \\
$ \B(a,p,c) \land \T(b,q,c)   \land \neg\, Col(a,b,c) \implies$& \\
\qquad $  \B(p,ip(a,p,c,b,q),b) \land  \B(q,ip(a,p,c,b,q),a)$ & (A7-i), Skolemized \\
$ \B(a,b,c) \implies \B(c,b,a) $ &(A14-i), symmetry of betweenness \\
$ \B(a,b,d) \land \B(b,c,d) \implies \B(a,b,c)$&(A15-i), inner transitivity\\
$Col(u,v,p) \land u \neq v \land  \T(a,p,b)    \implies $&   two-point line-circle \\
$   \exists  y,z\,( az = ab \land ay = ab \land \T(y,p,z) \land (p \neq a \implies y \neq z))$& \\ 
\endaxioms
\medskip

The axioms of stability are as follows:
\medskip

\axioms
$\neg \neg \B(a,b,c) \implies \B(a,b,c)$ & \\
$\neg \neg E(a,b,c,d) \implies E(a,b,c,d)$ & \\
$\neg a \neq b \implies a = b$ &
\endaxioms
\medskip

For reference we also state the circle-circle continuity principle, which is not 
an axiom but a theorem.  The circles must have distinct centers but one of them 
could be a null circle (zero radius).  See Fig.~\ref{figure:CircleCircleContinuityFigure}.
\medskip

\axioms
$  \ ap=ax \land aq=az \land\ cp = cd \ \land cq=cd\ \land \T(a,x,b)\  \land $ & \hskip-0.3cm$\T(a,b,z) \land a \neq c$    \\
$  \implies \exists z_1,z_2\,(cz_1=cd \land az_1 = ab \land  cz_2=cd \land az_2 = ab) $  & (circle-circle) 
\endaxioms

\subsection{Ruler-and-compass Tarski geometry}
This theory uses  \LPT\ (logic of partial terms) as given in \cite{beeson-book}, p.~97, which 
allows a formal treatment of ``undefined terms''.  
Its axioms are similar to intuitionistic Tarski geometry with Skolem functions, except that 
there is an additional 4-ary function symbol $\il$ with the axioms
\medskip

\axioms
$ Col(a,b,x) \land Col(p,q,x) \land \neg\, (Col(a,b,p)\land Col(a,b,q))\implies$& \\
\qquad $  x = \il(a,b,p,q)$ & Axiom $\il$-i \\
$ \il(a,b,p,q) \defined \implies   Col(a,b,\il(a,b,p,q)) \land Col(p,q,\il(a,b,p,q)) $ & Axiom $\il$-ii 
\endaxioms
\medskip

The Skolem term $ip(a,p,c,b,q)$ is replaced in the Skolemized inner Pasch axiom by
$\il(a,q,b,p)$.  Point $c$ does not occur in this term.  The term $center(a,b,c)$ in the 
triangle circumscription axiom is not changed.

\section{Tarski's axioms, continuity, and ruler and compass}
Two of Tarski's axioms have ``degenerate cases'', in the sense that they introduce points
that do not depend continuously on the parameters of the axiom.  (The two axioms are 
segment extension, which permits extending a null segment, and inner Pasch, 
which allows the diagram to collapse to a line.)  Even using classical logic,
we consider this undesirable.  We would like to have a formulation of Tarski's theory that 
would permit us to use Herbrand's theorem to show that if $\exists y\, A(x,y)$ is provable
(where $x$ stands for several variables, not just one),
then there are finitely many ruler and compass constructions $t_1(x),\ldots,t_n(x)$ such 
that for each $x$, one of the $t_i$ constructs the desired $y$, i.e.,$A(x,t_i(x))$. In this 
section, we discuss how Tarski's axioms can be slightly modified to eliminate 
discontinuities. It may be worth pointing out that a reformulation is {\em necessary},
as Tarski's formulation definitely does {\em not} have this property:  all ruler and 
compass constructions produce points that depend continuously on parameters, but as remarked
above, the existential theorems of Tarski's theory can produce points that do not depend
continuously on parameters.

\subsection{Segment extension and Euclid I.2}
(A4) is the segment construction axiom.  Tarski's version is  $\exists x\, (\T(q,a,x) \land ax = bc)$.
The degenerate case  is extending a null segment, i.e.,when $q=a$; then the point
$x$ is not uniquely determined, and moreover, $x$ 
does not depend continuously on $q$ as $q$ approaches $a$.
One might wonder if $x=a$,  or in other words $b=c$ (extending by a null segment) 
is also a degenerate case, but we do not consider it as degenerate, 
since there is no discontinuous dependence in that case.    Then to avoid
degenerate cases, we could  consider
\medskip

\axioms
$q \neq a \implies \exists x\,(\T(q,a,x) \land ax = bc)$ & (A4-i) 
\endaxioms
\medskip

Classically, disallowing $q=a$ costs nothing, since to extend a null segment $aa$ 
by $bc$, we just pick any point $d \neq a$ and extend the non-null segment $da$ by $bc$.
Of course, this introduces a discontinuous dependence. 

\subsection{Degenerate cases of inner Pasch}
(A7) is inner Pasch; please refer to Fig.~\ref{figure:InnerOuterPaschFigure}.
 This has a degenerate case when $p=a$ and $q=b$,  for as $(p,q)$ 
approaches ($a,b$), the intersection point $x$ does not have a unique limit,  but could approach
any point on $ab$ or not have a limit at all, depending on how $(p,q)$ approaches $(a,b)$.  
If $p=c$ or $q=c$, or if $p=a$ but $q \neq b$, or if $q=b$ but $p \neq a$,
then there is an obvious choice of $x$, so this degenerate case can be 
removed simply by replacing $\T$ by $\B$ in inner Pasch.  In fact, it is enough 
to make this replacement on one side of the triangle, leaving $\T$ on the other side.
 
Tarski's version of inner Pasch allows the points $a$, $b$, and $c$ to be collinear,
and this case is technically important, because it allows a number of fundamental theorems about 
betweenness to be derived that originally were taken as axioms.%
\footnote{Tarski viewed it as a good thing when the number of axioms could be 
reduced by using degenerate cases of remaining axioms.   We note that in 2013, a further possible 
reduction in the number of axioms was proved possible by Makarios \cite{makarios2013}:  
interchanging two variables in the conclusion 
of the five-segment allows the elimination of the symmetry axiom of congruence, $ab = ba$.
} 
The point asserted
to exist is unique when $a$, $b$, and $c$ are not collinear;  the technical 
question arises, whether the point can be chosen continuously in the five parameters
$a$, $b$, $c$, $p$, and $q$,  in case collinearity is allowed, but the five points are
required to be distinct.  Some computations (not provided here) show that indeed the 
point can be continuously chosen. 

Nevertheless, we consider the case when $a$, $b$, and $c$ are collinear to be objectionable,
on philosophical grounds.  Pasch's axiom is supposed to justify the construction of certain 
points by labeling the intersections of lines drawn with a straightedge as actually ``existing''
points.  In the case when the lines coincide, the axiom has no conceptual connection with 
the idea of intersecting lines, and hence would need some other justification to be accepted
as an axiom.  If the justification is just that it provides a single axiom from which several
intuitively evident propositions about betweenness can be deduced,  that is a distortion 
of the meaning of the word ``axiom.''  

 Whether or not one gives weight to this philosophical argument,
there is a related technical point:  we consider below a version of geometry with terms for the 
intersection points of lines, and we want to be able to use those terms to construct the points
shown to exist by Pasch's axiom.  In other words, the problem with Tarski's too-general version 
of inner Pasch is that it asserts the existence of points for which there is no ruler and compass
construction.   In that respect, it is unlike any of the other axioms (A1) to (A10), and also 
unlike the line-circle and circle-circle continuity axioms.
 This issue reflects in a precise mathematical way the 
philosophical issue about the collinear case of Pasch's axiom.  
 
Therefore, we reformulate inner Pasch for continuity, and for constructivity in the 
sense of ruler and compass constructions of the points asserted to exist, 
as follows:
\begin{itemize}
\item We change $\T(a,p,c)$ to $\B(a,p,c)$.
\item  We add the hypothesis, $\neg Col(a,b,c)$.
\end{itemize}

The resulting axiom is 
\smallskip

\newaxioms{1}
$ \T(a,p,c) \land \T(b,q,c) \land p \neq a   \land \neg\, Col(a,b,c) \implies$&\\
$\exists x\,(\B(p,x,b) \land \B(q,x,a))$ & (A7-i) strict inner Pasch \\
\endaxioms
\medskip

Note that we did not require both $\B(a,p,c)$ and $\B(b,q,c)$.  Changing just one of those
from $\T$ to $\B$ is  
sufficient to allow a ruler and compass construction.  We do not need two versions,
one with $\B(a,p,c)$ and one with $\B(b,q,c)$, by symmetry. 
As it turns out, we could use $\B$ instead of $\T$ in all the parts of this axiom and prove the same 
theorems, as is shown in Section~\ref{section:meta} below.
\medskip

\subsection{Inner Pasch and betweenness}\label{section:betweenness}
Tarski's final theory \cite{tarski-givant} had only one betweenness axiom, known as (A6) or ``the identity
axiom for betweenness'':
 $$\T(a,b,a) \implies a=b.$$
In terms of strict betweenness, that becomes $\neg \B(a,x,a)$,  or otherwise expressed,
$\B(a,b,c) \implies a \neq c$.   We also refer to this axiom as (A6).   
The original version of Tarski's theory had  more betweenness axioms (see \cite{tarski-givant}, p. 188).
These were all shown eventually to be superfluous in classical Tarski geometry, through the work of Eva Kallin, Scott Taylor,
Tarski himself, and especially Tarski's student H.~N.~Gupta \cite{gupta1965}.   These proofs 
appear in \cite{schwabhauser}. 
Here we give the axiom numbers from \cite{tarski-givant}, names by which they are known,
 and also the theorem numbers of their 
proofs in \cite{schwabhauser}:
\smallskip

\axioms
$\T(a,b,c) \implies \T(c,b,a)$ &  (A14), symmetry, Satz 3.2 \\
$\T(a,b,d) \land \T(b,c,d) \implies \T(a,b,c)$ & (A15), inner transitivity, Satz 3.5a\\
$\T(a,b,c) \land \T(b,c,d) \land b \neq c \implies $ & \\
\qquad $\T(a,b,d)$ & (A16), outer transitivity, Satz 3.7b \\
$\T(a,b,d) \land \T(a,c,d) \implies$ & \\
\qquad $\T(a,b,c) \lor \T(a,c,b)$ & (A17), inner connectivity, Satz 5.3 \\
$\T(a,b,c) \land \T(a,b,d) \land a \neq b \implies$ & \\
\qquad $\T(a,c,d) \lor \T(a,d,c)$ & (A18), outer connectivity, Satz 5.1
\endaxioms

\noindent

The first of these (A14), is a consequence of inner Pasch, formulated with 
$\T$, but the proof uses a degenerate case of inner Pasch, so if we replace inner Pasch by
the non-degenerate form (strict inner Pasch), we will (apparently) have to reinstate
(A14) as an axiom.  The question arises as to whether this is also true of the others.
Certainly these cases suffice:

\begin{Lemma} \label{lemma:collinearA7}
(A14) and (A15) suffice to prove the collinear case of Tarski's inner Pasch, using 
(A4-i) and (A7-i) instead of (A4) and (A7).
That is, 
$$ Col(a,b,c) \land a \neq b \land \T(a,p,c) \land \T(b,q,c) \implies \exists x\,(\T(p,x,b) \land \T(q,x,a)).$$
\end{Lemma}

\noindent{\em Proof}.   We first note that $\T(a,b,b)$ follows immediately from the 
definition of $\T(a,b,c)$ in terms of $\B$.  

Since we checked above that the degenerate cases of (A7) are provable, 
we can assume that all five of the given points are distinct. Since $Col(a,b,c)$, we have 
$\B(a,b,c) \lor \B(a,c,b) \lor \B(c,a,b)$.    

Case 1, $\B(a,b,c)$.  Then we take $x=b$.  We have to prove $\T(p,b,b) \land \T(q,b,a)$.
From $\T(a,b,c) \land \T(b,q,c)$ we have 
 $\T(a,b,q)$ by (A15).  Then $\T(q,b,a)$ by (A14).  Since $p \neq b$
we have
$\T(p,b,b)$ as shown above.  That completes Case 1.

Case 2, $\B(c,a,b)$.  Then we take $x=a$.  We have to prove $\T(p,a,b) \land \T(q,a,a)$.
Since $q \neq a$ we have $\T(q,a,a)$ as shown above.  By symmetry (A14) we have 
$\T(a,p,c)$ and $\T(b,a,c)$, so by (A15) we have $\T(b,a,p)$, so by (A14) again
we have $\T(p,a,b)$ as desired. That completes Case 2.

Case 3, $\B(a,c,b)$.  Then we take $x=c$.  We have to prove $\T(p,c,b) \land \T(q,c,a)$.
From $\T(a,c,b)$ and $\T(c,q,b)$ we have by (A15) $\T(a,c,q)$, whence by (A14), $\T(q,c,a)$.
From $\T(a,c,b)$ by (A14), we have $\T(b,c,a)$.  From $\T(a,p,c)$ by (A14), we have $\T(c,p,a)$.
From that and $\T(b,c,a)$ we have by (A15) $\T(b,c,p)$.  By (A14) we have $\T(p,c,b)$ as desired.
That completes Case 3, and the proof of the lemma.

\section{Alternate formulations of Tarski's theory}
In this section we consider some reformulations of Tarski's theories (still using 
classical logic) that (i) isolate and remove ``degenerate cases'' of the axioms, and 
(ii) introduce Skolem functions to achieve a quantifier-free axiomatization, and 
(iii) introduce additional axioms to make the intersection points of lines and circles,
or circles and circles, depend continuously on the (points determining the) lines and circles.

\subsection{Continuous Tarski geometry}
Let ``continuous Tarski geometry'' refer to classical
Tarski geometry with two-point line-circle continuity, with the following modifications:  
\begin{itemize}
\item  (A4-i) instead of (A4)  (extending non-null segments)
\item  (A7-i) (strict inner Pasch) instead of (A7).  That is, use \,$\B$ instead of $\T$ in 
two of the three occurrences of $\T$ in inner Pasch, and require $\neg Col(a,b,c)$.
\item  Take (A14) and (A15) as axioms (symmetry and transitivity of betweenness)
\item  Use the triangle circumscription principle (A10${}_3$) for the parallel axiom
\end{itemize}
The reason for the name ``continuous Tarski geometry'' will be apparent 
eventually, when we show what seems intuitively obvious: that Skolem functions for these axioms
can be implemented by ruler and compass constructions.    

\begin{Theorem} \label{theorem:continuousTarskiGeometry}
Continuous Tarski geometry has the same theorems as Tarski geometry.
\end{Theorem}

\noindent{\em Proof}.  To extend a null segment $bb$ by $cd$, first select any 
point $a$ different from $b$, then extend $ab$ by $cd$.  Hence the restriction to (A4-i)
costs nothing.  By Lemma~\ref{lemma:collinearA7}, the restriction to the non-collinear 
and non-degenerate case of (A7) is made up for by the inclusion of (A14) and (A15) as axioms. 
That completes the proof of the theorem.

\subsection{Skolemizing Tarski's geometry}
Since Tarski's axioms are already in existential form, one can add Skolem functions to make 
them quantifier-free.   Perhaps the reason why Tarski did not do so is his desire that there should 
be just one model of his theory over the real plane $\R^2$.  If one introduces Skolem functions 
for the intersection points of two circles, then those Skolem functions can be interpreted 
quite arbitrarily, unless one also adds further axioms to guarantee their continuity, and even 
then, one has a problem because those Skolem functions will be meaningless (have arbitrary values)
when the circles do not intersect.  Tarski did not have circle-circle continuity, but the same 
problem arises with Skolem functions for inner Pasch, when the hypotheses are not satisfied.

The problem can be seen in a simpler context, when we try to axiomatize field theory 
with a function symbol $i(x)$,  the official version of $x^{-1}$.  The point is that 
0 has no multiplicative inverse, yet Skolem functions are total, so $i(0)$ has to 
denote {\em something}.   We phrase the axiom as $x \neq 0 \implies x\cdot i(x) = 1$,
so we can't prove $0 \cdot i(0) = 1$, which is good, since we {\em can} prove $0\cdot i(0) = 0$.
In spite of this difficulty, the theory with Skolem functions is a conservative extension of 
the theory without Skolem functions,  as one sees (for theories with classical logic) 
from the fact that every model of the theory without can be expanded by suitably interpreting
the Skolem function symbols.  We return below to the question of how this works for 
intuitionistic theories in Lemma~\ref{lemma:conservativeSkolem} below.

Papers on axiomatic geometry often use the phrase ``constructive theory'' to mean one with 
enough function symbols to be formulated with quantifier-free axioms.  While this is not sufficient 
to imply that a theory is ``constructive'' in the sense of being in accordance with Bishop's
constructive mathematics (or another branch of constructive mathematics), it is a desirable feature,
in the sense that a constructive theory should provide terms to describe the objects it can prove
to exist.  In finding a constructive version of Tarski's theories, therefore, we will wish to 
produce a version with function symbols corresponding to ruler and compass constructions.  
In order to compare the constructive theory with Tarski's classical theory,  we will first 
consider a Skolemized version of Tarski's theory, with classical logic.

\subsection{Skolem functions for classical Tarski}
One introduces Skolem functions and reformulates the axioms to be quantifier-free.
But we want these Skolem functions to be meaningful as ruler and compass constructions.
Hence, we do not Skolemize Tarski's theory as he gave it, but rather the modified
version we called ``continuous Tarski geometry.''  The axioms are listed for reference
in \S\ref{section:axioms}; here we just give a list of the Skolem functions:
\begin{itemize}
\item
$ext(a,b,c,d)$ is a point $x$ such that for $a\neq b$, we have 
 $\T(a,b,x) \land bx = cd$.
\item
$ip(a,p,c,b,q)$  is the point asserted to exist by inner Pasch (see Fig~\ref{figure:InnerOuterPaschFigure}),
provided $a$, $b$, and $c$ are not collinear, and $\B(a,p,c)$.
\item
Three constants $\alpha$, $\beta$, and $\gamma$ for three non-collinear points.
(In this paper we consider only plane geometry, for simplicity.)
\item $center(a,b,c)$ is a point equidistant from   $a$, $b$, and $c$, provided 
$a$, $b$, and $c$ are not collinear. 
\item
$\ilcone(a,b,c,d)$ and $\ilctwo(a,b,c,d)$ for the 
two intersection points of  $\Line(a,b)$ and $\Circle(c,d)$, the circle with center $c$   
passing through $d$.  
\end{itemize}
The function $center$ is needed to remove the existential quantifier in Szmielew's
parallel axiom (A$10_2$), which says that if $a$, $b$, and $c$ are not collinear, there exists a 
circle through $a$, $b$, and $c$.  For the version (A10) of the parallel axiom used in \cite{schwabhauser},
we would need two different Skolem functions.  The points asserted to exist by that version are 
not unique and do not correspond to any natural ruler and compass construction, which is a  
reason to prefer triangle circumscription as the parallel axiom.

The question arises,  what do we do about ``undefined terms'',  e.g., $\ilcone(a,b,c,d)$ when 
the line and circle in question do not actually meet?  One approach is to modify the logic,
using the ``logic of partial terms'',  introducing a new atomic statement $t \defined$
(read ``$t$ is defined'')  for each term $t$.  In Tarski's geometry as described here, 
that is not necessary, since we can explicitly give the conditions for each term to be defined.
In that way, $t \defined$ can be regarded as an abbreviation at the meta-level, rather than
an official formula.  We write the formula as $(t \defined)^\circ$ to avoid confusion and 
for consistency of notation with another section below. 

\begin{Definition} \label{definition:definedness}
 When the arguments to the Skolem 
functions are variables or constants, we have
\begin{eqnarray*}
(ext(a,b,c,d) \defined)^\circ  &:=&   a\neq b  \\
ip(a,p,c,b,q) \defined)^\circ  &:=& \B(a,p,c) \land \T(b,q,c)   \land \neg Col(a,b,c) \\
(center(a,b,c) \defined)^\circ  &:=& \neg Col(a,b,c)
\end{eqnarray*}
If the arguments $a$, $b$, $c$, $d$ are not variables or constants, then we need
to add (recursively) the formulas expressing their definedness on the right.
\end{Definition}
In addition to the obvious ``Skolem axioms'' involving these function symbols,
we need additional axioms to ensure that the two intersection points of a line 
and circle are distinguished from 
each other (except when the intersection is of a circle and a tangent line), and that
the intersection points depend continuously on the (points determining the) lines and circles.

We discuss the two points of intersection of $\Line(a,b)$ and $\Circle(c,d)$,
which are denoted by $\ilcone(a,b,c,d)$ and $\ilctwo(a,b,c,d)$.  We want an axiom
asserting that these two points occur on $\Line(a,b)$ in the same order as $a$ and $b$ do;
that axiom serves to distinguish the two points and ensure that they depend continuously on 
$a, b, c$, and $d$.  To that end we need to define $SameOrder(a,b,c,d)$, assuming $a \neq b$
but allowing $c=d$. This can be done as follows:
\begin{eqnarray*}
SameOrder(a,b,c,d) &:=& (\T(c,a,b) \implies \neg \B(d,c,a))  \\
&& \land \  (\T(a,c,b) \implies \neg \B(d,c,b)) \\
&& \land \  (\T(a,b,c) \implies \T(a,c,d))
\end{eqnarray*}
The axiom in question is then
$$ SameOrder(a,b,\ilcone(a,b,c,d),\ilctwo(a,b,c,d)).$$

\subsection{Continuity of the Skolem functions }
We will investigate what additional axioms are necessary to guarantee that the Skolem functions
are uniquely defined and continuous.  Unless we are using the logic of partial terms,
technically Skolem functions are total, in which case we cannot avoid some arbitrariness 
in their values, but when their ``definedness conditions'' given above are satisfied,
we expect them to be uniquely defined and continuous. 
This will be important for metatheorems
about the continuous dependence on parameters of things proved constructively to exist; 
but we think it is also of interest even 
to the classical geometer.

Evidently for this purpose we should use 
the version of the axioms that has been sanitized of degenerate cases.  Thus,  $ext$ only 
Skolemizes axiom (A4-i),  for extending non-degenerate segments, and $ip$ only Skolemizes
axiom (A7-i) rather than A7.  These Skolem functions will then be uniquely defined (and provably so).

Continuity of a term $t(x)$, where $x$ can be several variables $x_1,\ldots,x_n$, can be defined in geometry: 
it means that for every circle $C$ with $t(x)$ as center,  where $C$ is given by $t(x)$ 
and a point $p$ on $C$, there exist circles $K_i$ about $x_i$ such that if $z_i$ is inside
$K_i$, for $i=1,\ldots,n$, then $t(z)$ is inside $C$.

\begin{Lemma} \label{lemma:ic-continuity}
The terms $\ilcone(a,b,c,d)$ and $\ilctwo(a,b,c,d)$  are provably continuous in 
$a$, $b$, $c$, and $d$ (when their definedness conditions hold).
\end{Lemma}

\noindent{\em Proof}.  Once we have defined multiplication and addition, this
proof can be carried out within geometry, using ordinary algebraic calculations.
 It is very much easier to believe that these (omitted)
proofs can be carried out, than it is to actually get a theorem-prover or proof-checker
to do so.  See \cite{beeson-edinburgh} for a full discussion of the issues involved.
\medskip

\begin{Theorem}[Continuity of inner Pasch]
Tarski's geometry, using axioms (A4-i) and (A7-i), proves the continuity of $ip(a,p,c,b,q)$
as a function of its five parameters, when the hypotheses of inner Pasch are satisfied.
\end{Theorem}

 {\em Remark}. If we use axiom (A7), without the modifications in (A7-i), then
$ip(a,p,c,b,q)$ is not continuous as $(p,q)$
approaches $(a,b)$, as discussed above.
\medskip

\noindent{\em Proof}.  This also can be carried out by introducing coordinates and making ordinary
algebraic computations within Tarski geometry.

\subsection{Continuity and the triangle circumscription principle}
Above we have given the triangle circumscription principle with the hypothesis that 
$a$, $b$, and $c$ are non-collinear (and hence distinct) points.  What happens when that
requirement is relaxed?  If $a$ and $b$ are allowed to approach each other without restriction
on the direction of approach, then $center(a,b,c)$ does not depend continuously on its parameters.
But if $a$ and $b$ are restricted to lie on a fixed line $L$ (as is the case when using 
triangle circumscription to define multiplication as Hilbert did), then as $a$ approaches $b$ 
(both remaining away from $c$),
the circle through
$a$, $b$, and $c$ nicely approaches the circle through $a$ and $c$ that is tangent to $L$
at $a$.   The {\em strong triangle circumscription principle} says that there is a term 
$C(a,b,c,p,q)$ such that when $a$ and $b$ lie on $L = \Line(p,q)$ and $c$ does not lie on $L$,
then $e=C(a,b,p,q)$ is equidistant from $a$, $b$, and $c$, and moreover, if $a=b$ then $ea$
is perpendicular to $L$ at $a$ (i.e., the circle is tangent to $L$ at $a$).  In \cite{beeson-bsl},
it is shown how to construct the term $C$, using segment extensions and the uniform perpendicular;
so this construction can be carried out in  Tarski geometry with Skolem functions.

 \section{A constructive version of Tarski's theory} \label{section:constructivity}
Finally we are ready to move from classical to intuitionistic logic.
Our plan is to give two intuitionistic versions of Tarski's theory, 
one with function symbols as in the Skolemized version above, and one 
with existential axioms as in Tarski's original theory.
The underlying logic will be intuitionistic predicate logic.
We first give the specifically intuitionistic parts of our theory, which are very few in number.
We do not adopt decidable equality ($a=b \lor a \neq b$), 
nor even the substitute concept of ``apartness'' introduced by Brouwer and Heyting (and discussed below),  primarily because we 
aim to develop a system in which definable terms (constructions) denote continuous functions, but also because we
wish to keep our system closely related to Euclid's geometry, which contains nothing like apartness.

\subsection{Introduction to constructive geometry} 
Here we discuss some issues particular to geometry with intuitionistic logic.
The main point is that we must avoid case distinctions in existence proofs.  
 What one has 
to avoid in constructive geometry is not proofs of equality or inequality by contradiction,
but rather constructions (existence proofs) that make a case distinction. For example, 
classically we 
have two different constructions of a perpendicular through point $p$ to line $L$, one 
for when $p$ is not on $L$, and another for when $p$ is on $L$.  Pushing a double negation 
through an implication, we only get not-not a perpendicular exists, which is not enough.
To constructivize the theorem, we have to give a uniform construction of the perpendicular,
which works without a case distinction. (Two different such constructions are given in 
this paper, one using line-circle continuity, but not the parallel axiom, and one using the parallel axiom, but not line-circle.)

In particular, in order to show that the models of geometry are planes over Euclidean 
fields, we need to define addition and multiplication by just such uniform constructions,
without case distinctions about the sign of the arguments.  The classical definitions
due to Descartes and Hilbert do depend on such case distinctions; in \cite{beeson-bsl} we 
have given uniform definitions; here we check that their properties can be proved
in intuitionistic Tarski geometry.   To actually carry out the complete development
directly would be a project
of about the length and scope of Szmielew's comparable development of classical geometry
from Tarski's axioms, in Part I of \cite{schwabhauser}.  Therefore it is important
that the double-negation interpretation can be made to carry the load.

We mention here two principles which are not accepted by 
all constructivists, at least in the context of real analysis.
 Here $x < y$ refers to points on a fixed line $L$, and can be 
defined in terms of betweenness.
\smallskip

 \axioms
$ \neg \neg x > 0 \implies x > 0$ & (Markov's principle) \\
$ x \neq 0 \implies x < 0 \lor x > 0$ & (two-sides)
\endaxioms
\smallskip

We accept the former, but not the latter.
Markov's principle follows from the stability of betweenness and is a fundamental principle of 
constructive geometry.  It allows us to avoid distinguishing more than one sense of 
inequality between points.  Geometry without it would be much more complicated.

The  principle ``two-sides'' (which we do not accept)  is closely related to ``a point not on a given line is on one side or the other 
of the line''.  (Here the ``line'' could be the $y$-axis, i.e.,a line perpendicular to $L$ at the point 0.)  This principle is not needed in the 
formalization of Euclid, or the development of the geometrical theory of arithmetic, and as 
we will show, it is not a theorem of intuitionistic Tarski geometry.

One might consider adopting two-sides as an axiom, on grounds similar to those sometimes
used to justify Markov's principle or apartness, namely that if we ``compute $x$ to sufficient accuracy
we will see what sign it has.''  That justification applies only to the model of computable reals, 
not to various more general intuitionistic models of sequences generated by free choices of 
approximations to points.  Brouwer argued against this principle in one of his later papers
\cite{brouwer1949b}
on those grounds; and our development of constructive geometry shows that it is not needed
for the usual theorems, including the geometric definitions of addition and multiplication.
In our opinion, not only is it unnecessary, it is also constructively undesirable, as 
the choice of which disjunct holds cannot depend continuously on $x$,  so anyone claiming its
validity must make some assumptions about how points are ``given'', e.g. by a computable 
sequence of rational
approximations;  we do not want to make such assumptions.

On the other hand, the following principle {\em has} been accepted by all constructivists in 
the past who considered geometry:
\smallskip

\axioms
$ a < b \implies x < b \lor a < x$ & (apartness) 
\endaxioms
\smallskip

It turns out that apartness is completely unnecessary for the formalization of Euclid, 
and is not a theorem of intuitionistic Tarski geometry.
  The desire to use apartness probably arose from an 
unwillingness to the trichotomy law of order, and to find some replacement for it.
In our work,  the law of trichotomy of order is replaced by the stability of equality
and betweenness.   If we want to formalize one of Euclid's proofs where two points 
are proved equal by contradiction (consider III.4 for a specific example),  the proof
in Euclid shows $\neg a \neq b$;  in other words $\neg \neg a=b$.  What we need
to formalize such proofs is the principle 
$$ \neg a \neq b \implies a=b$$
or, otherwise expressed, $ \neg \neg a=b \implies a=b$.  This principle,
already mentioned in the introduction, is called the ``stability of 
equality.''   The trichotomy law can also be double negated, each case  but one
shown contradictory, and the final double negation removed by the stability of 
betweenness, $\neg \neg \B(a,b,c) \implies \B(a,b,c)$.   That is the 
fundamental reason why apartness is not needed in constructive geometry.

\subsection{Stability} \label{section:stability}
The word ``stable'' is applied to a predicate $Q$ if $\neg \neg Q \implies Q$. 
 Our
intuitionistic versions of Tarski geometry will all have axioms of stability for the basic predicates.
That is, we include the axioms
\begin{eqnarray*}
\neg a \neq b \implies a =b \\
\neg \neg \B(a,b,c) \implies \B(a,b,c) \\
\neg \neg ab=cd \implies ab=cd
\end{eqnarray*}
In this section we justify accepting these axioms.  Our 
 intuition is that there is nothing asserting existence in the meaning 
of equality, congruence, or betweenness; hence assertions of equality, congruence, or betweenness
 can be constructively proved by contradiction.  There are many examples in Euclid%
 \footnote{Just to mention one, Euclid III.4}
  where 
 Euclid argues that two points, differently constructed, must coincide; such examples use 
 the stability of equality.  Similarly, if point $x$ lies on $\Line(a,b)$, we may 
 wish to argue by cases as to its position on the line relative to $a$ and $b$.  We
  double-negate the disjunction of the five possible
 positions, argue each case independently, and arrive at the double negation of the 
 desired conclusion.  As long as what we are proving is a 
 betweenness, congruence, or equality,  stability allows us to remove the double 
 negation and reach the desired conclusion.
 
 We explain this point with more detail,
 for those inexperienced with intuitionistic reasoning:  Suppose $P \implies Q$,
 and $R \implies Q$.   Then $(P \lor R) \implies Q$ (both classically and intuitionistically).
 Taking $R$ to be $\neg P$, if $P \implies Q$ and $\neg P \implies Q$, then 
 $P \lor \neg P \implies Q$.  So, classically, $Q$ holds.  But intuitionistically,
 we may not be able to prove that the cases $P$ and $\neg P$ are exhaustive;  for 
 example we cannot assert in general that point $p$ is on line $L$ or it is not.
 But intuitionistically,
 we still have $\neg\neg\,(P \lor \neg P)  \implies \neg \neg\, Q$, since if $\neg Q$ then 
 $\neg P$ and $\neg \neg P$, which is contradictory.
 Now if $Q$ is stable we can still conclude $Q$,  since $\neg \neg (P \lor \neg P)$
 is intuitionistically valid.  
 
 What we are not allowed to do, constructively, is argue by cases for an existential 
 conclusion, using a different construction for each case.  (In the previous 
 paragraph, if $Q$ begins with $\exists$, then $Q$ will not be stable.) 
 This observation makes it 
 apparent why the constructivization of geometry hinges on the successful discovery of 
 {\em uniform} constructions,  continuous in parameters.
 
 As we mentioned above, angles can be defined in Tarski's theory, and one can show 
 that  the equality and 
ordering
of angles is stable.  That is, $$\neg \neg \alpha < \beta \implies \alpha < \beta$$
for angles $\alpha$ and $\beta$.  Thus, when Euclid wants to prove $\alpha = \beta$, 
and says,  if not, then one of them is greater; let $\alpha > \beta$, and so on, 
the reasoning is constructive, because we have 
$$\neg \neg(\alpha < \beta \ \lor \ \alpha = \beta \ \lor \ \beta < \alpha)$$
and if $\alpha < \beta$ and $\beta < \alpha$ lead to contradictions, then 
$\neg \neg \alpha = \beta$,  whence by stability, $\alpha = \beta$.  Similarly 
if what is to be proved is an inequality of angles.

Julien Narboux pointed out that the stability of equality can be derived from 
the stability of congruence:

\begin{Lemma} \label{lemma:stability2981}
With the aid of axioms A1 and A3,  stability of congruence implies stability of 
equality.
\end{Lemma}

\noindent{\em Proof}.  Suppose $\neg a \neq b$.  We want to prove $a=b$.
By A3, it suffices to prove $ab=aa$.  By the stability of congruence,
we may prove this by contradiction. Suppose, for proof by contradiction,
that $ab \neq aa$.  We claim $a \neq b$.  To prove it, suppose
 $a=b$. Then from $ab \neq aa$ we obtain $ab \neq ab$,  contradicting A1. 
Therefore $a \neq b$. But that contradicts the hypothesis $\neg a \neq b$ from 
the first line of the proof.  That completes the proof of the lemma.
\medskip

We could therefore drop stability of equality as an axiom, but we retain it 
anyway, because of its fundamental character, and to emphasize that it is 
perhaps even more fundamental than the facts expressed in A1 and A3.
\medskip

{\em Stability of incidence}. 
Tarski's theory has variables for points only, so when we discuss lines,
implicitly each line $L$ is given by two points, $L = \Line(a,b)$.  When 
we say point $x$ lies on line $L$, that abbreviates 
$$ \neg \neg\, (\T(x,a,b) \lor \T(a,x,b) \lor \T(a,b,x)). $$
Since logically, $\neg\neg\neg\, P$  is equivalent to $\neg\, P$,  
the relation $x$ lies on $L$  is stable. (Four negations is the same as two negations.)
We refer to this as the ``stability of incidence.''  When the definition of 
incidence is expanded, this is seen to be a logical triviality, not even worth
of the name ``lemma.''  But when working in Tarski's theories with less than
complete formality, we do mention lines and incidence and justify some proof steps
by the ``stability of incidence.''  In other words, we are allowed to prove that 
a point $x$ lies on a line $L$ by contradiction.

\subsection{Strict and non-strict betweenness}
Should we use strict or non-strict betweenness in constructive geometry? 
The answer is, it doesn't matter much, because of the stability of $\B$. What we 
do officially is use strict betweenness $\B$,  and regard $\T$ as defined by 
$$ \T(a,b,c) :=  \neg(a \neq b \land b \neq c \land \neg \B(a,b,c)).$$
We could also have taken $\T$ as primitive and defined $\B$ by 
$$ \B(a,b,c) := \T(a,b,c) \land a \neq b \land b \neq c.$$

\subsection{Intuitionistic Tarski geometry with existential axioms}
The language of this theory takes strict betweenness $\B(a,b,c)$ as primitive,
and $\T(a,b,c)$ will then be a defined concept, given by Definition~\ref{defn:T}.
Some of the axioms will be ``unmodified'' from Tarski's theory, by which we mean that 
the only change is to define $\T$ in terms of $\B$.  The other modifications 
are described in detail in \S\ref{section:axioms2}, just before the listing of the modified
axioms.  Here we  summarize the modifications:

\begin{itemize}
\item  Modify Axiom A4 (segment extension) so only non degenerate segments are extendable.
\item Axiom (A6) becomes $\neg \B(a,b,a)$.  
\item Replace inner Pasch (A7) by (A7-i), which  requires $\neg Col(a,b,c)$ and replaces
two of the three occurrences of $\T$ by $\B$.
\item Add (A14-i) and (A15-i), the symmetry and inner transitivity axioms for betweenness:
$$ \B(a,b,c) \implies \B(c,b,a)$$
$$ \B(a,b,d) \land \B(b,c,d) \implies \B(a,b,c)$$
\item Use  the triangle circumscription principle as the parallel postulate.
\item Add two-point line-circle continuity (instead of segment-circle).
\item Use intuitionistic logic only.  
\item add the stability of equality, betweenness, and congruence.
\end{itemize}

The resulting theory is called ``intuitionistic Tarski geometry'', or 
``intuitionistic Tarski geometry with existential axioms.''
Another way of describing it is:  restrict continuous Tarski geometry to intuitionistic logic,
and add the stability axioms for equality, betweenness, and equidistance, use 
the triangle circumscription principle for the parallel axiom, and add two-point line-circle 
continuity.  We use the phrase ``intuitionistic Tarski geometry without any continuity''
to refer to the theory obtained by dropping the line-circle continuity axiom.

\begin{Theorem} \label{theorem:intuitionisticTarskiImpliesClassical}
Intuitionistic Tarski geometry plus classical logic is equivalent to Tarski geometry
(with or without line-circle continuity).
\end{Theorem}

\noindent{\em Proof}.  This follows from Theorem~\ref{theorem:continuousTarskiGeometry},
Since intuitionistic Tarski geometry is continuous Tarski geometry with intuitionistic
logic and stability axioms, it is classically equivalent to continuous Tarski geometry.
But by Theorem~\ref{theorem:continuousTarskiGeometry},  that theory is classically
equivalent to Tarski geometry.  That completes the proof.

\subsection{Intuitionistic quantifier-free Tarski geometry}
We can use the same Skolem  functions as for the classical theory, since we 
already made the necessary restrictions to the Skolem functions for segment extension
and Pasch's axioms.  For the same reason, the
conditions for definedness of Skolem terms are not changed.

\begin{Lemma} \label{lemma:stabilityOfDefined}
For every term $t$ of intuitionistic quantifier-free Tarski geometry, the
sentence $ \neg \neg (t \defined)^\circ \implies (t \defined)^\circ$ is provable.
\end{Lemma}

\noindent{\em Proof}.  By induction on the complexity of $t$, using the stability of 
$\B$, $E$, and equality for the base case.

Since the conditions for the definedness of Skolem terms are definable, there is no
logical problem about using (total) Skolem functions in this intuitionistic theory,  without 
modifying the logic, which is the ordinary intuitionistic first-order predicate calculus.
However, there might be a philosophical problem, as one might ask, what is the intended 
interpretation of those total Skolem symbols?  One cannot specify a total (everywhere defined)
construction to interpret, for example, the Skolem symbol for inner Pasch.  Therefore it 
is more philosophically correct to use the ``logic of partial terms'', which is explained
in another section below.  However, it is possible to consider the Skolem symbols as mere
syntactic tools, which, even if not meaningful, at least cause no unwanted deductions, 
according to the following lemma:

\begin{Lemma} \label{lemma:conservativeSkolem} [Conservativity of Skolem functions]
Suppose intuitionistic Tarski with Skolem functions proves a theorem $\phi$ that does 
not contain Skolem functions.  Then intuitionistic Tarski (with existential quantifiers
and no Skolem functions) also proves $\phi$.\  In fact, the same is true of any 
intuitionistic theory whose axioms before Skolemization have the form 
$P(x) \implies \exists y\, Q(x,y)$, with $P$ quantifier-free.
\end{Lemma}

\noindent{\em Proof}.  Consider a Skolem symbol with axiom $P(x) \implies Q(x,f(x))$,
Skolemizing the axiom $P(x) \implies \exists y\, Q(x,y)$.  
The corresponding lemma for classical theories needs no restriction
on the form of $P$;  one simply shows that every model of the theory without Skolem 
functions can be expanded to a model of the theory with Skolem functions.  The interpretation
of the values of a Skolem symbol, say $f(b)$ are just taken arbitrarily when $P(b)$ is not 
satisfied.  Then one appeals to the completeness theorem.  One can use the Kripke completeness
theorem to make a similar argument for theories with intuitionistic logic;  but in general 
one cannot define $f(b)$ at a node $M$ of a Kripke model where $P(b)$ fails, because $P(b)$
might hold later on, and worse, there might be nodes $M_1$ and $M_2$ above $M$ at which different
values of $y$ are required, so there might be no way to define $f(b)$ at $M$.  That cannot 
happen, however, if $P$ is quantifier-free, since then, if $P(b)$ does not hold at $M$,
it also doesn't hold at any node above $M$.  Hence if $P$ is quantifier free, we can complete
the proof, using Kripke completeness instead of G\"odel completeness.  (For an 
introduction to Kripke models and a proof of the completeness theorem, see 
\cite{troelstra}, Part V, pp.~324{\em ff}.)

\section{Perpendiculars, midpoints and circles} \label{section:linecircle}
We have included two-point line-circle continuity in our axiom systems for 
ruler and compass geometry, since this corresponds to the the physical use of ruler
and compass.  Tarski, on the other hand, had  segment-circle continuity.   In this 
section, we will show how to construct perpendiculars and midpoints, using two-point
line-circle continuity.  By definition, $ab$ is perpendicular to $bc$ if $a \neq b$ and 
$abc$ is 
a right angle, which means that if $da = ac$ and $\B(d,a,c)$ then $ad=ac$.  We 
sometimes write this as $ab \perp bc$.  More generally, two lines $K$ and $L$ 
are perpendicular at $b$ if $b$ lies on both lines, and there are points $a$ and $c$ 
on $K$ and $L$ respectively such that $ab \perp bc$.  We sometimes write this as $K \perp L$.
It can be shown that $K \perp L$ if and only if $L \perp K$.

We will deal with ``dropped perpendiculars''
(perpendiculars from a point $p$ to a line $L$, assuming $p$ is not on $L$), and also with 
``erected perpendiculars'' (perpendiculars to a line $L$ at a point $p$ on $L$).  
Our results are of equal interest for classical and constructive geometry.
As in \cite{schwabhauser}, we abstain  from the use of the dimension axioms 
or the parallel axiom (except in one subsection, where we mention our use of 
the parallel axiom explicitly); although these restrictions are not stated in the lemmas they
are adhered to in the proofs. 

A fundamental fact about perpendiculars is the uniqueness of the dropped perpendicular:
\begin{Lemma} \label{lemma:perpunique}
Suppose $p$ does not lie on line $L$, but $a$ and $b$ do lie on $L$, and 
both $pa$ and $pb$ are perpendicular to $L$.  Then $a=b$.
\end{Lemma}

\noindent{\em Proof}.  This is Satz~8.7 of \cite{schwabhauser}; the proof 
offered there is completely constructive, but it does appeal to some earlier 
theorems.  One can either check these proofs directly, or appeal to the 
double-negation interpretation (Theorem~\ref{theorem:doublenegation} below).
\medskip

We need to verify that the midpoint of a segment can be constructed by a term 
in intuitionistic Tarski geometry with Skolem functions.  If we were willing 
to use the intersection points of circles, the problem might seem simple: we could use
the two circles drawn in Euclid's Proposition I.1, and connect the two intersection 
points.  (This is not Euclid's construction of midpoints, but still it is commonly used.)
This matter is not as simple as it first appears, as we shall now explain.
 
We try to find the midpoint of segment $pq$. 
 Let $K$ be the circle with center $p$ passing through $q$, and let
$C$ be the circle with center $q$ passing through $p$, and let $d$ and $e$ be the two 
intersection points of these circles.  
Now the trick would be to prove that $de$ meets $pq$ in a point $f$;  if that could 
be done, then it is easy to prove $f$ is the desired midpoint, 
by the congruence of triangles $pef$ and $qef$.  But it seems at first that 
the the full Pasch axiom 
is required to prove the existence of $f$.  True,  full Pasch follows from inner Pasch
and the other axioms, 
at least classically, but we would have to verify that constructively using only the 
axioms of intuitionistic Tarski, which does not seem trivial.  In particular, 
we will need the existence of midpoints of segments to do that!  (Moreover, the dimension
axiom would need to be used; without it, circle-circle continuity is sphere-sphere 
continuity and not every two intersection points $d$ and $e$ of two spheres will lie
in the same plane as $pq$. Full Pasch fails in three dimensions, while inner and 
outer Pasch hold.)

In fact, the existence of midpoints has been the subject of much research, and 
it has been shown that one does not need circles and continuity at all!  Gupta \cite{gupta1965}
(in Chapter 3) showed that inner Pasch suffices to construct midpoints, i.e.,classical
Tarski proves the existence of midpoints.  Piesyk (who was a student of Szmielew) proved it \cite{piesyk1965}, using outer Pasch instead of inner Pasch.   Later
Rigby  \cite{rigby1975} reduced the assumptions further.  At the end of this section,
we will give Gupta's construction, but not his proof.  Since the proof just proves 
that the constructed point $m$ satisfies $ma = mb$, by the stability of equality
(and the double-negation interpretation, technically, which we shall come to 
in \S\ref{section:doublenegation}) we know that Gupta's classical 
proof \cite{gupta1965, schwabhauser} can be made constructive.   The simpler construction 
using line-circle continuity is adequate for most of our purposes. 

\subsection{The base of an isosceles triangle has a midpoint}
Euclid's own midpoint construction is to construct an isosceles triangle on $pq$ and 
then bisect the vertex angle. 
One of Gupta's simple lemmas enables us to justify the second part of this 
Euclidean midpoint construction, and we present that lemma next.

\begin{Lemma} \label{lemma:midpoint-helper} [Gupta]
Intuitionistic Tarski geometry with Skolem functions, and without continuity, 
proves that for some term $m(x,y,z)$,  if $y \neq z$ and  $x$ is equidistant from 
$y$ and $z$, with $x$, $y$, and $z$ not collinear, then $m(x,y,z)$ is the midpoint of $yz$.
\end{Lemma}

\begin{figure}[h]
\caption{To construct the midpoint $w$ of $yz$, given $x$ with $xz=xy$, using 
two applications of inner Pasch.}
\label{figure:GuptaMidpointFigure}
\GuptaMidpointFigure
\end{figure}
\medskip

\noindent
{\em Proof}.  See page 56 of \cite{gupta1965}.  But the proof is so simple and beautiful 
that we give it here.
Let $\alpha$ and $\beta$ be two of the three distinct points guaranteed by Axiom A8.  Let 
\begin{eqnarray*}
t &=& ext(x,y,\alpha,\beta) \\
u &=& ext(x,z,\alpha,\beta) \\  
v &=& ip(u,z,x,t,y) \\
w &=& ip(x,y,t,z,v)
\end{eqnarray*}
Then $w$ is the desired midpoint.  The reader can easily check this; see
 Fig.~\ref{figure:GuptaMidpointFigure} for illustration. 
 Thus we can define 
\begin{eqnarray*}
m(x,y,z) &=& ip(x,y,t,z,v) \\
&=& ip(x,y,t,z,ip(u,z,x,t,y)) \\
&=& ip(x,y,ext(x,y,\alpha,\beta) ,z,ip(ext(x,z,\alpha,\beta),z,x,ext(x,y,\alpha,\beta) ,y)) 
\end{eqnarray*}
That completes the proof.
\medskip

Since circle-circle continuity enables us to construct an equilateral triangle on any 
segment (via Euclid I.1),
we have justified the existence of midpoints and perpendiculars if circle-circle
continuity is used.  But that is insufficient for our purposes, since intuitionistic
Tarski geometry does not contain circle-circle continuity as an axiom.  (We 
show below that it is a theorem, but to prove it, we need midpoints and perpendiculars,
so they must be obtained some other way.)

\subsection{A lemma of interest only constructively}

In erecting a perpendicular to line $L$ at point $a$, we need to make use of 
a point not on $L$ (which occurs as a parameter in the construction).
 In fact, we need a point $c$ not on $L$ such that $ca$ is not perpendicular to $L$.
 Classically, we can make the case distinction whether $ca$ is perpendicular to $L$,
 and if it is, there is ``nothing to be proved''. But constructively, this case
 distinction is not allowed, so we must first construct such a point $c$.  Our 
 first lemma does that.
 
 \begin{Lemma} \label{lemma:droppedperpstarter}
 Let point $a$ lie on line $L$, and point $s$ not lie on line $L$.
 Then there is a point $c$ not on $L$ such that $ca$ is not perpendicular to $L$, 
 and a point $b$ on $L$ such that $\B(c,s,b)$.
 \end{Lemma}
 
 \noindent{\em Remark}.  Classically, this lemma is trivial, but constructively,
 there is something to prove.
 \medskip

 \noindent{\em Proof.}  Let $b$ be a point on $L$ such that $ab = as$.
 (Such a point can be constructed using only the segment extension axiom.)
Then by Lemma~\ref{lemma:midpoint-helper}, $sb$ has a midpoint $c$,
and $ac \perp bc$.
We claim that $ac$ cannot be perpendicular to $L$; for it were,
then triangle $cab$ would have two right angles, one at $c$ 
and one at $a$.   That contradicts Lemma~\ref{lemma:perpunique}.
That completes the proof. 

\subsection{Erected  perpendiculars from triangle circumscription}
Following Szmielew, we have taken as our form of the parallel axiom, the 
axiom that given any three non-collinear points, there is another point equidistant
from all three.   (That point is then the center of a circumscribed circle containing
the three given points.)  An immediate corollary is the existence of midpoints and 
perpendiculars.

\begin{Lemma} \label{lemma:midpointcircum} 
Every segment has a midpoint and a perpendicular bisector.  
If $p$ is a point not on line $L$,  then there is a point $x$ on $L$ with $px \perp L$.
\end{Lemma}

\noindent{\em Proof}.  Let $ab$ be given with $a \neq b$.  By Lemma~\ref{lemma:droppedperpstarter},
there is a point $c$ such that $a$, $b$, and $c$ are not collinear.   By
triangle circumscription, there exists a point $e$ such that $ea=eb=ec$.  Then $eab$
is an isosceles triangle.  By Lemma~\ref{lemma:midpoint-helper}, $ab$ has a 
midpoint $m$.   Since $ea = eb$, we have $em \perp ab$, by definition of perpendicular.
That completes the proof.
\medskip

While this is formally pleasing, there is something unsatisfactory here,  because we 
intend that the axioms of our theories should correspond to ruler and compass 
constructions, and in order to construct the point required by the triangle circumscription 
axiom, we need to construct the perpendicular bisectors of $ab$ and $ac$, and find 
the point of intersection (whose existence is the main point of the axiom).
So from the point of view of ruler and compass constructions, our argument has been 
circular.   
Therefore the matter of perpendiculars cannot be left here.  We must consider how 
to construct them with ruler and compass.  

 \subsection{Dropped perpendiculars from line-circle continuity}
In this section we discuss the following method (from Euclid~I.12) of dropping a perpendicular from 
point $p$ to line $L$:  draw a large enough circle $C$ with center $p$ that 
some point of $L$ is strictly inside $C$.  Then apply two-point line-circle continuity to get two 
points $u$ and $v$ where $C$ meets $L$.  Then $puv$ is an isosceles triangle, so 
it has a midpoint and perpendicular bisector.  This argument is straightforwardly 
formalized in intuitionistic Tarski geometry:

\begin{Theorem}[dropped perpendiculars from line-circle] \label{theorem:droppedperp}
One can drop a perpendicular 
to line $L$ from a point not on $L$ by a ruler and compass construction 
(essentially the construction of Euclid~I.12), and prove the construction
correct in  intuitionistic Tarski geometry.   
\end{Theorem}

\noindent{\em Proof}.  Let $p$ be a point, and let $L$ be the line through two 
distinct points $a$ and $b$.  Suppose $p$ is not on $L$.
 Let $r = ext(p,a,a,b)$, so that $pr$ is longer than 
$pa$.  Let $C$ be the circle with center $p$ passing through $r$.  Then $a$ is 
strictly inside $C$.  By line-circle continuity,  there are two points 
$x$ and $y$  on $L$ where $L$ meets $C$, i.e.,$px=py$.    Since by 
hypothesis $p$ is not on $L$,  the  segment $xy$ is the 
base of an isosceles triangle, so by Lemma~\ref{lemma:midpoint-helper}
it  has a midpoint $m$.  Then $pm \perp L$, by definition of perpendicular.
That completes the proof.

We note that the analogous lemmas with one-point line-circle or segment-circle
in place of two-point line-circle will not be proved so easily.   In the case
of segment-circle,  we would need to construct points outside $C$ in order to 
apply the segment-circle axiom; but to show those points are indeed outside $C$,
we would need the triangle inequality, which requires perpendiculars for its proof.
In the case of one-point line-circle,  we would need to construct the other intersection 
point, and the only apparent way to do that is to first have the dropped perpendicular
we are trying to construct.

Of course, these difficulties are resolved if we are willing to use the   
1965 discoveries of Gupta, who showed how to construct perpendiculars without any
use of circles.  But that is beside the point here, since we are considering whether 
our choice of axioms corresponds well to Euclid or not.  Also,  the use of 
the triangle circumscription axiom is no help, since although we could prove 
the existence (if not the construction) of erected perpendiculars, we cannot do the 
same for dropped perpendiculars.

\subsection{Reflection in a line is an isometry}
Since we can drop perpendiculars,
we can define reflection in a line (for points not on the line). We need to know that reflection preserves
betweenness and equidistance.   That reflection preserves equidistance (i.e., is an isometry),  
 is  Satz~10.10 in \cite{schwabhauser}, but the proof uses nothing but things proved before 
the construction of perpendiculars and midpoints late in Chapter~8.   Strangely, it 
is not explicitly stated in \cite{schwabhauser} that reflection in a line 
preserves betweenness, so we begin by proving that.

\begin{Lemma} \label{lemma:reflectionB} Reflection in a line preserves betweenness.
\end{Lemma}

\noindent{\em Proof}.  Suppose $\B(a,b,c)$, and let $p,q,r$ be the reflections of $a,b,c$ 
in line $L$.     Since reflection
is an isometry, we have $pq = ab$ and $qr = bc$ and $pr = ac$.  We wish to prove 
$\B(p,q,r)$.  By the stability of betweenness, we may use proof by contradiction,
so assume $\neg \B(p,q,r)$. Since $pq = ab < ac = pr$, we do not have $\B(p,r,q)$.  
Similarly we do not have $\B(q,p,r)$.  Since $a,b,c$ are distinct points, the 
reflections $p,q,r$ are also distinct.  Therefore $q$ is not on $Line(p,r)$
(see the discussion of stability of incidence after Lemma~\ref{lemma:stability2981}).
 
Therefore $pqr$ is a triangle in which the sum of two sides $pq$ and $qr$
equals the third side $pr$.  By Lemma~\ref{theorem:droppedperp},
we can drop a perpendicular from $q$ to $pr$.  Let $t$ be the foot of that perpendicular.
Each of the two right triangles $qtp$ and $qtr$  has its base less than its 
hypotenuse, so $pr < pq + qr$, contradiction.  (We did not assume $\B(p,t,r)$ in this 
argument.)  That completes the proof of the lemma.

\subsection{Erected perpendiculars from dropped perpendiculars}

Next we prove the existence of erected perpendiculars, assuming only that 
we can drop perpendiculars, and without using any form of the parallel axiom.
   Gupta's proof,
as presented in Satz~8.21 of \cite{schwabhauser}, accomplished this. It
is much simpler than Gupta's beautiful circle-free construction of 
dropped perpendiculars, but still fairly complicated.   Gupta uses 
his Krippenlemma.  Below we give a new proof, avoiding the use of the
Krippenlemma, and using only simple ideas that 
might have been known to Tarski in 1959.  However, we do use both 
inner and outer Pasch; but in 1959, Tarski knew that outer Pasch implies inner.
The proof of outer Pasch from inner given by Gupta relies on the existence 
of dropped perpendiculars, but is far simpler than Gupta's circle-free,
parallel-axiom-free construction of dropped perpendiculars.  Therefore,
the work in this section gives us simple constructions in intuitionistic Tarski
of erected perpendiculars, assuming either the parallel axiom or line-circle.
In other words, the Skolem terms for the constructions below will involve
either $center$ or $ilc_1$ and $ilc_2$ (the symbols for intersections of lines and circles),
 depending on how we construct dropped perpendiculars;
if we use Gupta's construction, then we only need the Skolem symbol $ip$ for inner Pasch.  

The proof we are about to give requires two lemmas,
which we prove first.

\begin{figure}[ht]
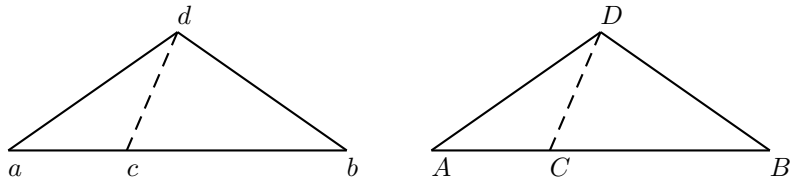

\caption{The interior 5-segment lemma.  $cd = CD$.}
\label{figure:interior5}
\InteriorFiveFigure
\end{figure}

\begin{Lemma}[Interior 5-segment lemma] \label{lemma:interior5}
The 5-segment axiom is valid if 
we replace $\B(a,b,c)$ and $\B(A,B,C)$ by $\B(a,c,d)$ and $\B(A,C,B)$, respectively.
That is, in Fig.~\ref{figure:interior5}, if the corresponding solid segments are 
congruent, so are the dashed segments.
\end{Lemma}
\noindent{\em Proof}.  This is Satz~4.2 in \cite{schwabhauser}. The proof is an 
immediate consequence of axioms (A4) and (A6); we do not repeat it.

\begin{Theorem}[Erected perpendiculars] \label{theorem:erectperp}
In intuitionistic Tarski geometry with line-circle continuity but 
without circle-circle continuity and without the parallel axiom, 
one can prove the following:  Let $L$ be a 
line; let $a$ be a point on $L$ and let $s$ be a point not on $L$.  Then there exists
a point $p$ on the opposite side of $L$ from $s$ such that $pa \perp L$.

The point $r$ on $L$ such that $\B(p,r,s)$, as well as the point $p$,
are given by terms of intuitionistic Tarski geometry with Skolem functions.
\end{Theorem}

\noindent{\em Remarks.}  No dimension axiom is used, and no parallel axiom.
 The point $s$ and line $L$
determine a plane in which $p$ lies.  We need $s$ constructively also because we 
do not (yet) know how, without circle-circle continuity, to construct a point
not on $L$ by a uniform construction.  We need $pa$ longer than the perpendicular from $s$ to 
$L$ to avoid a case distinction in the construction of midpoints  in the next lemma.  
\medskip

\noindent{\em Proof}.  By Lemma~\ref{lemma:droppedperpstarter}, we can find 
 a point $c$ not on $L$, and a point $b \neq a$ on $L$,
  such that $ca$ is not perpendicular to $L$, and $\B(c,s,b)$.

\begin{figure}[ht]
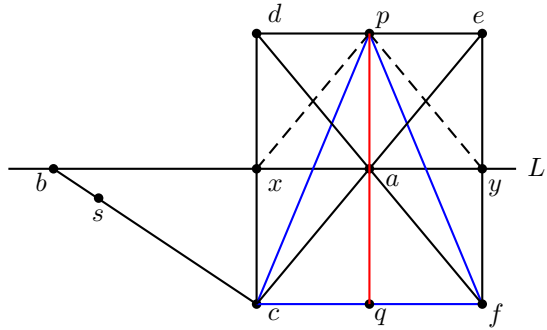

\caption{Erecting a perpendicular to $L$ at $a$, given $b$, $s$, and $c$.}
\label{figure:erectperp}
\ErectPerpFigure
\end{figure}

By Lemma~\ref{theorem:droppedperp}, let $x$ be a point on $L$
such that $cx \perp L$.  Then $x \neq a$, since $ca$ is not perpendicular to $L$.
We have $c \neq x$ since $c$ is not on $L$ but $x$ is on $L$. Hence we can construct
the reflection $d$ of $c$ in $x$, namely
$d = ext(c,x,c,x)$.
Similarly, $c \neq a$, so we can define $e$ to be the reflection of $c$ in $a$.
Since angle $dxa$ is a right angle, $da = ca$.
Since $ea=ca$, we have $da = ea$.  
Then triangle $ade$ is isosceles, and hence segment $de$ has a midpoint, 
by Lemma~\ref{lemma:midpoint-helper}.  Let $p$ be that midpoint. 
Let $y$ be the reflection of $x$ in $a$, and $f$ the reflection 
of $d$ in $a$.   Fig.~\ref{figure:erectperp} illustrates the situation.

Since reflection preserves congruence, we have 
 $xd=xc = ey $, and $yf = xd = xc = ey$, and $ca = da = fa$.  Since reflection
 preserves betweenness, we have $\B(e,y,f)$. (Thus $f$ is the reflection of $e$ 
 in $y$ as well as the reflection of $d$ in $a$.)  Let $q$ be the reflection
 of $p$ in $a$.  Since $ca = fa$,
 angle $aqf$ is a right angle. Therefore  $pc = pf$.  
Now consider the four-point configurations $(c,x,d,p)$ and 
$(f,y,e,p)$.  Of the five corresponding pairs of segments, four are congruent,
so by the interior 5-segment lemma (Lemma~\ref{lemma:interior5}, $px = py$.  Since $xa = ya$,  angle $pax$ is 
a right angle and $pa \perp L$, as desired.

It remains to prove that $p$ is on the other side of $L$ from $s$.
Fig.~\ref{figure:erectperpfinish} illustrates the situation.

\begin{figure}[ht]
\caption{$p$ is on the other side of $L$ from $s$, by constructing first $t$ and then $r$.}
\label{figure:erectperpfinish}
\ErectPerpFinishFigure
\end{figure}

Segment $xa$ connects two sides of triangle $cde$, and $cp$ connects the vertex $c$
to the third side, so by the crossbar theorem (Satz~3.17 of \cite{schwabhauser},
derived with two applications of inner Pasch), there is a point $t$  on $xa$
(and hence on $L$) with $\T(c,t,p)$.   That is, $p$ is on the other side of $L$
from $c$.  Using $\B(c,s,b)$, we can apply  inner Pasch to the five points $ptcbs$,
yielding a point 
$r$ such that $\B(t,r,b)$ and $\B(p,r,s)$.  Since $t$ and $b$ lie on $L$,
that shows that $p$ and $s$ are on opposite sides of $L$.   
That completes the proof.

\subsection{Construction of parallels}

We take the opportunity to point out that even in neutral geometry (that is, 
without any form of the parallel postulate) we can always construct a parallel $K$ to 
a given line $L$ through a given point $p$ not on $L$,  by dropping a perpendicular 
$M$ to $L$ through $p$ and then erecting $K$ perpendicular to $M$ at $p$.  In the 
sequel, when a reference is made to ``constructing a parallel to $L$ through $p$'',
this is what is meant. 

Using the uniform perpendicular construction to construct $K$ permits a similar
``uniform parallel'' construction, which, when given $L$ and $p$, constructs  
  a line $K$ through $p$ such that if $p$ is not on $L$,
then $K$ is parallel to $L$.   (But $K$ is constructed anyway, whether or not $p$ is on $L$;
if $p$ is on $L$, then $K$ coincides with $L$.)  The uniform parallel 
construction is important in establishing the properties of coordinates, and then
addition and multiplication, in \cite{beeson-bsl},  so the fact that the uniform parallel
can be defined in Tarski geometry permits the construction of coordinates and arithmetic.

\subsection{Midpoints from perpendiculars}
 Next we intend to construct the midpoint
of a non-degenerate segment $ab$.   This is Satz~8.22 in \cite{schwabhauser}.  That proof 
consists of a simple construction that goes back to Hilbert \cite{hilbert1899} (Theorem~26),
and a complicated proof that it really constructs
the midpoint. The proof of the correctness of Hilbert's construction from Tarski's axioms given in 
\cite{schwabhauser} is complicated, appealing to Gupta's ``Krippenlemma'', whose proof
is not easy.   Here we give another proof, not relying on the Krippenlemma.  We do 
need to use the properties of reflection in a line.  

\begin{Lemma} \label{lemma:midpointprep}
Let $a \neq b$, and suppose $ap \perp ab$ and $br \perp ab$, and $\B(a,m,b)$ 
and $\B(p,m,r)$ and $br = ap$, and $p$ is not on $Line(a,b)$.
   Then $m$ is the midpoint of $ab$, i.e.,$am=mb$.
\end{Lemma}

\noindent{\em Proof}. 
  By the stability of equality, 
we may use classical logic to prove $am = mb$; hence we may refer to the   
proof in \cite{schwabhauser}, page 65 (Abb. 31) (which appeals to the Krippenlemma,
whose proof is complicated) 
without worrying whether it is constructive.  
  But in the interest of giving a self-contained
and simple proof,  we will show how to finish the proof without appealing
to the complicated proofs of Gupta.

\begin{figure}[ht]
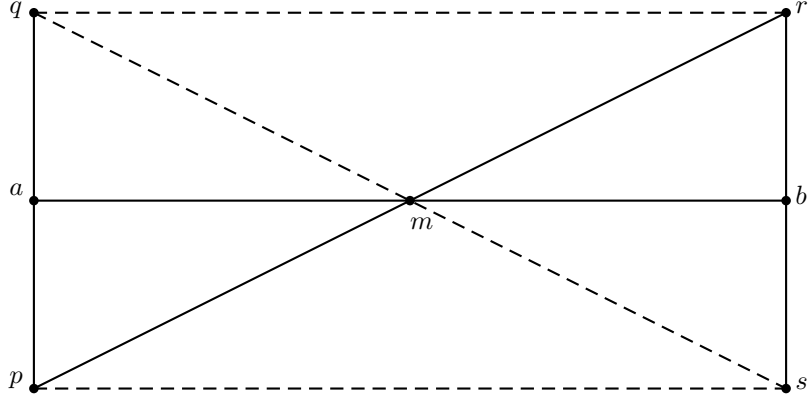

\caption{Given $ap\perp ab$ and $br \perp ab$ and $br=ap$, and $ab$ meets $pr$ at $m$, 
then $m$ is the midpoint of $ab$.}
\label{figure:MidpointButterflyFigure}
\MidpointButterflyFigure
\end{figure}

Let $q$ and $s$ be the reflections of $p$ and $r$, respectively, in $\Line(a,b)$.
Fig.~\ref{figure:MidpointButterflyFigure} illustrates the construction.

\noindent
Since $ap \perp ab$ and $br \perp ab$ we have $\B(q,a,p)$ and $\B(r,b,s)$.
Since reflection in a line preserves betweenness, and $m$ is its own reflection 
since $\B(a,m,b)$, we have $\B(q,m,s)$.  Hence $m$ is the intersection point of the 
diagonals of the quadrilateral $qpsr$.   Since reflection in a line is an isometry,
we have $qr = ps$.  Since $qa = ap = rb = bs$, we have $qp = rs$.  Hence the opposite
sides of quadrilateral $qpsr$ are equal.  By Satz~7.29,  the diagonals bisect each 
other.  Hence $mp = mr$ and $mq = ms$.  Now by the inner five-segment lemma 
(Lemma~\ref{lemma:interior5}), applied to the configurations $qapm$ and $sbrm$, we have $ma = mb$.
Then $m$ is the midpoint of $ab$.  That completes the proof.

\begin{Lemma} [midpoint existence] \label{lemma:midpointsG}
In intuitionistic Tarski geometry with only line-circle continuity,
midpoints exist.  More precisely, given segment $ab$ (with $a \neq b$)
and point $s$ not collinear with $ab$,
one can construct the midpoint $m$ of $ab$.  
\end{Lemma}

\noindent{\em Remark}. The proof shows how to derive the existence of midpoints
from the existence of (erected) perpendiculars;  we have shown above that
line-circle continuity enables us to erect perpendiculars.  But we use both
inner and outer Pasch.
\medskip 

\begin{figure}[ht]
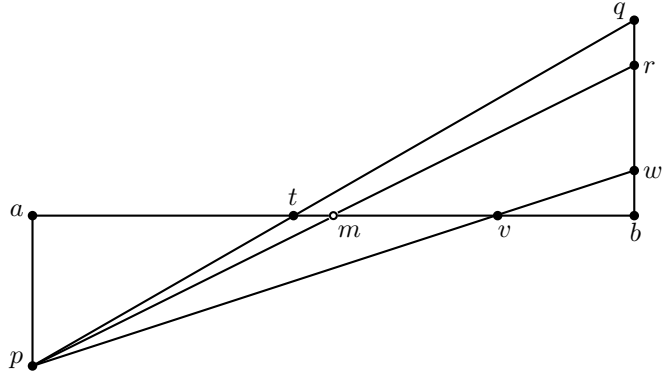

\caption{Midpoint from erected perpendiculars; construct $t$ and $m$ by outer and inner Pasch.}
\label{figure:midpointPasch}
\MidpointPaschFigure
\end{figure}

\noindent{\em Proof}. The construction is illustrated in Fig.~\ref{figure:midpointPasch}.
 Let segment $ab$ be given along with point $s$ not on $L=\Line(a,b)$.
Erect a perpendicular $ap$ to $ab$ at $a$ (on the opposite side of $L$ from $s$).
Then erect a perpendicular $wb$ to $L$ at $b$, with $w$ on the opposite side of $L$ from $p$,
  by Lemma~\ref{theorem:erectperp}.   Let $q = ext(b,w,a,p)$.  Then $\B(b,w,q)$ 
  and $ap < qb$.  Since $w$ is on the opposite side of $L$ from $p$, there is  
  a point $v$ on $L$ between $w$ and $p$.  Applying outer Pasch to $bwqpv$,  there 
exists a point $t$ with $\B(a,t,b)$ and $\B(p,t,q)$.
Then, construct point $r$
on segment $qb$ so that $rb = ap$, which is possible since $ap < qb$.

Applying inner Pasch to the five-point configuration
$ptqrb$,  we find a point $m$ such that $\B(a,m,b)$ and $\B(p,m,r)$.  
 By Lemma~\ref{lemma:midpointprep}, $m$ is the midpoint of $ab$.
That completes the proof.  
\medskip

\noindent{\em Constructions and terms}.  The above explicit construction 
corresponds to a term of intuitionistic Tarski geometry, built up by composing 
the terms for the construction steps.  Such a term can be written explicitly,
and we will explain that point now.  First, let
$Perp(a,b,s)$ and $wit(a,b,s)$ be  terms giving the construction 
in Lemma~\ref{theorem:erectperp}; that is, the first is a point on the perpendicular and 
the second is a witness that it is on the other side of $\Line(a,b)$.
Let $op$ and $ip$ be Skolem functions corresponding to outer Pasch and inner Pasch.
Then the following script corresponds to the proof above:

\begin{verbatim}
midpoint(a,b,s){
   p = Perp(a,b,s)
   w = Perp(b,a,p)
   v = wit(b,a,p)
   q = ext(b,w,a,p)
   t = op(b,w,q,p,v)
   r = ext(ext(w,b,w,b),b,a,p)
   m = ip(b,r,q,p,t)
   return m
  }
\end{verbatim}
Such a script can be converted to a (long) term by composing the right sides of the 
equations, eliminating the variables on the left.

From Gupta's proof of outer Pasch by means of inner Pasch, we can extract a term 
built up from $ip$ and $ext$ to substitute for $op$; and if desired, the terms $Perp$ 
and $wit$, which here involve $center$ or $ilc_1$ and $ilc_2$,  can be replaced by 
longer terms from Gupta's perpendicular construction (discussed below), built up from $ext$ and $ip$ only.

Regarding the role of $s$:  Of course since the midpoint is unique, the value does 
not depend on the parameter $s$.  Nevertheless we will not know until 
Lemma~\ref{lemma:pointoffline} below how to get rid of $s$ 
in the term, as  we need a point not on $\Line(a,b)$ to construct a perpendicular to $ab$.
That lemma gives us, in principle, the means to replace $s$ by a term in $a$ and $b$, 
but that term will involve $center$ and rely on the parallel postulate for its correctness.
\medskip
 
\begin{Corollary} In intuitionistic Tarski geometry (with two-point line-circle but 
without any parallel axiom)  
we can construct both dropped and erected perpendiculars, and midpoints.
\end{Corollary}

\noindent{\em Proof}.  We have justified Euclid~I.12 for dropped perpendiculars, and 
shown how to construct erected perpendiculars from dropped perpendiculars, and midpoints
from erected perpendiculars.   That completes the proof.

\medskip  
That result is pleasing, but it leaves open the question of whether the axiom system 
required can be weakened.  It is not obvious how to weaken it even to one-point line-circle
or segment-circle continuity; but Gupta showed in 1965 that no continuity whatsoever is 
required.  We discuss his construction in the next section.

\subsection{Gupta's perpendicular construction} \label{section:gupta}
Here we give Gupta's beautiful construction of a dropped perpendicular.  It 
can be found in his 1965 thesis \cite{gupta1965} and again in \cite{schwabhauser}, p.~61.

\begin{figure}[ht]
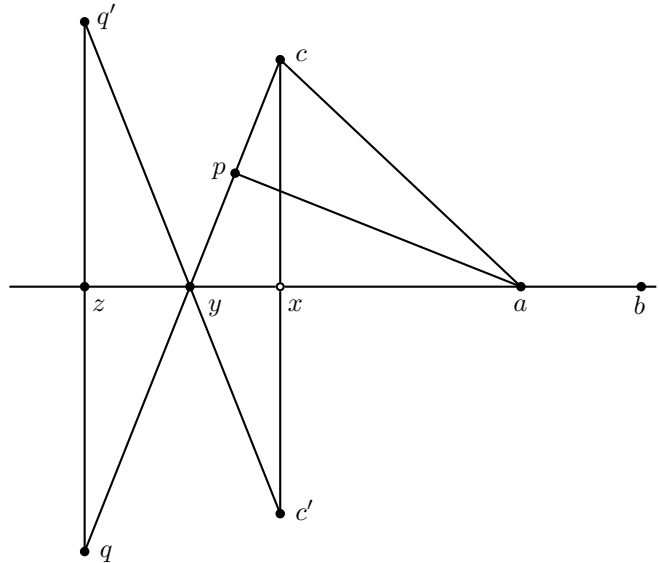

\caption{Gupta's construction of a dropped perpendicular from $c$ to $ab$.}
\label{figure:GuptaPerpendicularFigure}
\GuptaPerpendicularFigure
\end{figure}

 The initial data are two distinct points 
$a$ and $b$, and a third point $c$ not collinear with $ab$.
The construction goes as follows: 
Extend $ba$ by $ac$ to produce point $y$.
Then $acy$ is an isosceles triangle, so by Lemma~\ref{lemma:midpoint-helper},
we can construct its midpoint $p$, and $apc$ is a right angle. Now extend
segment $cy$ by $ac$ to point $q$, and extend $ay$ by $py$ to produce point $z$.
Then reflect $q$ in $z$ to produce $q^\prime$.  Then extend $q^\prime y$ by $yc$
to produce point $c^\prime$.  By construction $cyc^\prime$ is an isosceles triangle,
so its base $cc^\prime$ has a midpoint $x$, which is the final result of the 
construction.  It is not immediately apparent either that $x$ is collinear
with $ab$ or that $cx \perp ax$,  but both can be proved.   The proof 
occupies a couple of pages, but it uses only (A1)-(A7);  in other words,
no continuity,  no dimension axioms, no parallel axiom.  Is the proof 
constructive?  That point is discussed in  \S\ref{section:doublenegation}
below, where it is shown in Corollary~\ref{corollary:gupta} that at least,
Gupta's proof can be mechanically converted to a constructive proof.

\subsection{Erected perpendiculars from line-circle and the parallel axiom}

\begin{figure}[ht]
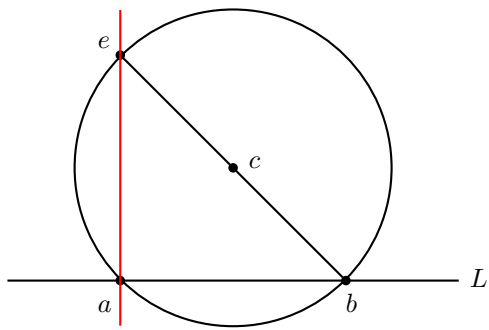

\caption{Erecting a perpendicular to $L$ at $a$, given $c$. First $b$ and then $e$ 
are constructed by line-circle continuity.}
\label{figure:perpfromcircle}
\PerpFromCircleFigure
\end{figure}

There is a very simple construction of an erected perpendicular based
on line-circle continuity and the parallel axiom.  It has only two ruler and compass
steps, and is therefore surely the shortest possible construction of an erected perpendicular.%
\footnote{It has three steps if you count drawing the final perpendicular line, but that step
constructs no new points.}
This construction is illustrated in Fig.~\ref{figure:perpfromcircle}. 

The construction starts with a line $L$, a point $a$ on $L$, and a point $c$
not on $L$ such that $ca$ is not perpendicular to $L$.  Then we draw the circle
with center $c$ through $a$, and let $b$ be the other point of intersection with 
$L$.  Then the line $bc$ meets the circle at a point $e$, and  $ea$ is the desired
perpendicular to $L$ (although we have not proved it here). 

The correctness of this construction, i.e.,that $ea$ is indeed perpendicular to $L$,
is equivalent to Euclid III.31, which says an angle inscribed in a semicircle
is a right angle.  Euclid's proof uses I.29, which in turn depends on I.11, which 
requires erecting a perpendicular.  Since III.31 is certainly not valid in hyperbolic
geometry, we will need to use the parallel postulate in any proof of III.31, and 
it seems extremely unlikely that we can prove III.31 without first having proved the 
existence of erected perpendiculars.   Thus,  this construction cannot replace 
the ones discussed above in the systematic development of geometry.  Still, it is 
of interest because it is the shortest possible ruler and compass construction.

Once we have perpendiculars, several
basic theorems follow easily.  The following can be proved
without using the parallel axiom or any continuity, from 
the assumption that dropped perpendiculars exist.  Hence,
without appealing to Gupta's construction, they can be proved
in intuitionistic Tarski geometry (using triangle circumscription).

\subsection{Angles and congruence of angles} \label{section:angles}
It is a curiosity that perpendicularity can be extensively studied without needing 
to discuss angles in general.  But angles are fundamental to Euclid and Hilbert.
Tarski's method of treating angles as triples of points means that we need a notion 
``$abc$ and $ABC$ are the same angle''.  To define that notion, we first define 
``$x$ lies on $Ray(b,a)$'' by $\neg (\neg \T(b,x,a) \land \neg \T(b,a,x))$.  
Then ``$abc$ and $ABC$ are the same angle'' means that the same points lie on $Ray(b,a)$
as on $Ray(B,A)$ and the same points lie on $Ray(b,c)$ as on $Ray(B,C)$.  Then there 
are several ways of defining congruence of angles, all equivalent. 
Specifically, two angles $abc$
and $ABC$ are congruent if by adjusting $a$ and $c$ on the same rays we can make $ab=AB$
and $bc = BC$ and $ac = AC$; or equivalently, if the points on all four rays can be so adjusted;
or equivalently, if any adjustment of $a$ and $b$ can be matched by an adjustment of $A$ and $B$.

 To avoid needing
a quantifier, Tarski proceeds as follows:  first extend $ba$ by $BA$ and $bc$ by $BC$.
Then extend $BA$ by $ba$ and $BC$ by $bc$.  Now we have two angles in the ``same angle''
relation to the original angles, whose corresponding sides are congruent segments.
  Hence it suffices to define congruence for angles 
$abc$ and $ABC$ with $ab = AB$ and $bc = BC$.  That definition is just $ac = AC$.
See \cite{schwabhauser}, Chapter 11, for a formal development of the 
basic properties of angles (which we do not use here).

The definition of angle ordering is given in \cite{schwabhauser}, Definition 11.27,
and the well-definedness depends on an argument by cases using inner and outer Pasch.
It would appear on the face of the matter that one needs a more general version of 
Pasch than inner and outer Pasch.  However, that is not actually the case: as soon 
as we have both inner and outer Pasch, we can derive any and every conceivable version
of Pasch.  For example, to constructivize the theory of angle congruence we need
a version that combines inner and outer Pasch, in the following sense.

\begin{Theorem}[Continuous Pasch] \label{theorem:continuousPasch}
Using inner and outer Pasch, we can derive the following:
Suppose $c \neq a$, and suppose $p$ is on the ray from $a$ through $c$ (but we assume nothing about 
the order of $p$ and $c$ on that ray).  Suppose that $b$ is not on $Line(a,d)$
and $\B(c,q,b)$.  
Then $bp$ meets the ray from $a$ through $q$ in a point $x$ such that $\B(b,x,p)$,
and if $\B(a,p,c)$ then $\B(a,x,q)$, and if $\B(a,c,p)$ then 
$\B(a,q,x)$.  
\end{Theorem}

\begin{figure}[ht]
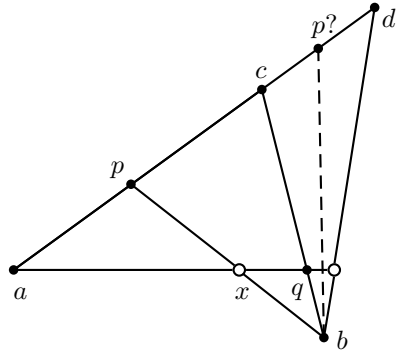

\caption{Continuous Pasch. Point $x$ exists regardless of which side of $c$ point $p$ is on.}
\label{figure:continuousPasch}
\ContinuousPaschFigure
\end{figure}

\noindent{\em Proof}.  Classically, we can simply argue by cases: if $\B(a,p,c)$,
$x$ exists by inner Pasch, and if $\B(a,c,p)$, then $x$ exists by outer Pasch.
Constructively, this case distinction is illegal.  Instead we proceed as follows.
First construct point $d$ such that $\B(a,p,d)$, for example by extending the 
non-null segment $ac$ by $ap$ to get $d$.  Then by outer Pasch, there exists
a point $r$ such that $\B(b,r,d)$ and $\B(a,q,r)$.  Then by inner Pasch 
there exists point $x$ such that $\B(a,x,r)$ and $\B(b,x,p)$.  That is the 
desired point $x$.  The additional properties of $x$ follow from inner and 
outer Pasch as in the classical argument.  That completes the proof of the 
theorem.

The SAS (side-angle-side) criterion for angle congruence is a postulate in Hilbert,
and Euclid's failed ``proof'' in I.4 shows that he, too, should have taken it as a postulate.
The five-segment axiom in Tarski is in essence a version of SAS, as discussed above where
the axiom is introduced.   But it should be noted that the five-segment axiom does not 
immediately cover the case of two right triangles with corresponding legs congruent.
That requires in addition the theorem that all right angles are congruent. 
This theorem is also needed to prove the uniqueness of the perpendicular to a line at a given point.
 
Euclid took the congruence of all right angles 
 as his Postulate 4.  Hilbert (\cite{hilbert1899}, p. 20) remarks 
that this was ``unjustified'', and says that the proof of it goes back to Proclus.  

 \begin{Lemma} \label{lemma:allrightanglescongruent} All right angles are congruent.   In other words, 
if $abc$ and $ABC$ are right angles with $ab = AB$ and $bc = BC$ then $ac = AC$.
\end{Lemma}

\noindent{\em Proof}.  This is Satz 10.12 in \cite{schwabhauser}.  However, the proof
appeals only to the definition of angle congruence and simple theorems, such as the 
fact that reflections in points and in lines are isometries.

\begin{Lemma} \label{lemma:triangleinequality}
(i) An exterior angle of a triangle is greater than either of the 
opposite interior angles.
\smallskip

(ii) The leg of a right triangle is less than the hypotenuse. 
\smallskip

(iii) If $a$, $b$, and $c$ are not collinear, the triangle inequality holds: $ac < ab + bc$.
\smallskip

(iv)  Whether or not $a$, $b$, and $c$ are collinear, we have $ac \le ab + bc$.
\end{Lemma}

\noindent{\em Proof}.   (i) is  Satz~11.41 of  \cite{schwabhauser};
(ii)  is a special case
of Satz~11.53.   Neither of these proofs uses anything but elementary 
betweenness and congruence, and the existence of perpendiculars.

Turning to the triangle inequality, let the non-degenerate triangle $abc$ be given. Drop a 
perpendicular $cx$ from $c$ to $\Line(a,b)$. 
  Since segment inequality is stable,
we are allowed to argue by cases on the position of $x$ relative to $a$ and $b$.
 If $\T(a,x,b)$ then 
$ab = ax + xb < ac + bc$ by (ii).  If $\B(a,b,x)$, then 
$ab < ax < ac$;  and similarly if $\B(x,a,b)$, then $ab < bx < bc$.  
 Hence constructively not not $ab < ac + bc$; 
   so $ab < ac +bc$.
That proves (iii).

Turning to (iv), we use the uniform perpendicular construction to get a perpendicular
to $Line(a,b)$ through $c$.  Let $x$ be the point where this perpendicular meets $\Line(a,b)$,
and proceed as for (iii), but with $\le$ in place of $<$.  That completes the proof of the lemma.

The following lemma is presented, not for its intrinsic interest, but 
because it is needed in what follows.  We needed the exterior angle theorem 
to prove it.

\begin{Lemma}\label{lemma:saccheri-helper}   Let $a,b,c$, and $d$ be   
distinct points, with $ab \perp ad$ and   $cd \perp ad$.  If $bc$ meets $Line(a,d)$,
then the intersection point $m$ is between $a$ and $d$.
\end{Lemma}

\noindent{\em Proof}.  
By the stability of betweenness, we may prove $\B(a,m,d)$ by contradiction,
so suppose $\neg \B(a,m,b)$.  Without loss of generality, we may
assume $\B(m,a,d)$.  Then angle $amb$ is an exterior angle of triangle $mcd$.
This is less than a right angle, since angle $mab$ is a right angle, and the 
angles of a triangle are less than two right angles. 
By the exterior angle theorem (Euclid I.16, or Lemma~\ref{lemma:triangleinequality} above, or Satz 11.41 in \cite{schwabhauser}), 
angle $amb$ is greater than
the interior angle $mdc$, which is a right angle; contradiction.
That completes the proof.

 \section{Other forms of the parallel axiom}
Within neutral geometry (that is, geometry without any form of the parallel postulate),
we can consider the logical relations between various forms of the parallel axiom.
In  \cite{beeson-bsl}, we considered  the Playfair axiom, Euclid 5, and the strong parallel axiom,
which are all classically equivalent to Euclid 5.  Constructively, Playfair is weaker, 
as shown in \cite{beeson-bsl}; a formal independence result confirms the intuition 
that it should be weaker because it makes no existential assertion.  The other versions
of the parallel postulate, which do make existential assertions, each turn out to be 
fairly easy to prove equivalent to either Euclid 5 or the strong parallel postulate. 
In \cite{beeson-bsl}, we prove that Euclid 5 and the strong parallel postulate are 
actually constructively equivalent, although the proof requires first developing the 
geometrical definitions of arithmetic using only Euclid 5. 

In \cite{beeson-bsl}, we gave the proof that the triangle circumscription principle is 
equivalent to the strong parallel axiom; below we prove that    
Tarski's version of the parallel
postulate taken as axiom (A10) in \cite{schwabhauser}, is equivalent to Euclid 5.
Lemmas in this section are proved in neutral geometry, i.e., without any form of the parallel postulate.
It follows that all the known versions of the parallel postulate (that are equivalent 
in classical Tarski geometry with line-circle continuity) that make an existential 
assertion are also equivalent in constructive Tarski geometry, although some of the proofs 
are much longer.

\begin{Lemma} \label{lemma:alternateinteriorangles}
Playfair's axiom implies the alternate interior angle theorem, that any line 
traversing a pair of parallel lines makes alternate interior angles equal.
\end{Lemma}

\noindent{\em Proof}.  Since ordering of angles is stable, we can argue by 
contradiction.  Hence the usual classical proof of the theorem applies.

\subsection{Euclid 5 formulated in Tarski's language}
\begin{figure}[ht]
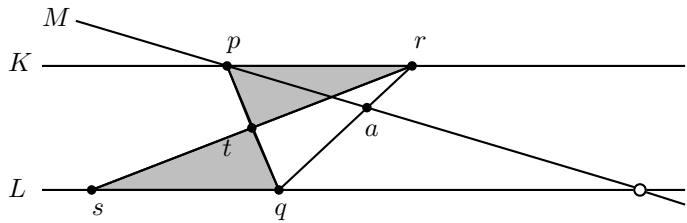

\caption{Euclid 5.  Transversal $pq$ of lines $M$ and $L$ makes corresponding interior angles less than 
two right angles, as witnessed by $a$. The shaded triangles are assumed congruent. Then $M$ 
meets $L$ as indicated by the open circle.}
\label{figure:EuclidParallelRawFigure}
\hskip 2.5cm
\EuclidParallelRawFigure
\end{figure}

Here we give a formulation of Euclid's parallel postulate, expressed in 
Tarski's points-only language.  Euclid's version mentions angles, and the 
concept of ``corresponding interior angles'' made by a transversal.  The following 
is a points-only version of Euclid 5.   See Fig.~\ref{figure:EuclidParallelRawFigure}.

\medskip

\axioms
$   \B(q,a,r) \land \B(p,t,q) \land pr=qs \land pt=qt \land rt=st   \land \neg\, Col(s,q,p)   $&(Euclid 5)\\
$  \qquad  \implies \exists x\, \B(p,a,x) \land \B(s,q,x)$&
\endaxioms
\smallskip

\subsection{Tarski's parallel axiom}
\begin{figure}[h]
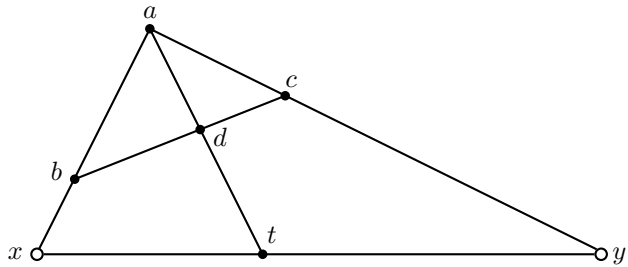
   
\caption{Tarski's parallel axiom  }
\label{figure:TarskiParallelFigure}  
\TarskiParallelFigure
\end{figure}

As mentioned above, Tarski in \cite{tarski1959} and later \cite{schwabhauser} took a different
form of the parallel postulate, illustrated in Fig.~\ref{figure:TarskiParallelFigure}.%
\footnote{Technically, according to \cite{tarski-givant}, the versions of the parallel 
axiom taken in the two 
cited references 
differed in the order of arguments to the last betweenness statement, but that is of no consequence.}
The following axiom is similar to Tarski's axiom, and we give it his name, but his axiom used
non-strict betweenness and did not include the hypothesis that $a$, $b$, and $c$ are not 
collinear.  It is intended to say that if $t$ is in the interior of angle $bac$, then 
there is a line through $t$ that meets both sides of the angle.    To express this using 
variables for points only, Tarski used the point $d$ to witness that $t$ is in the interior
of the angle.   See Fig.~\ref{figure:TarskiParallelFigure}.

The degenerate cases are trivial: if $a$, $b$, and $c$ are collinear, then we 
can (classically, or with more work also constructively) find $x$ and $y$ without any 
parallel axiom, and if (say) $d=b$ then we can take $x=t$ and $y=c$, etc.  Hence the following 
axiom is classically equivalent in neutral geometry to the one used by Tarski:
\smallskip

\axioms
$\B(a,d,t) \land \B(b,d,c) \land a \neq d  $&(Tarski parallel axiom) \\
$ \land\ (\neg\,\B(a,b,c) \land \neg\,\B(b,c,a) \land \neg\,\B(c,a,b)\land a \neq c) $ & \\
$ \qquad \implies \exists x \,\exists y\,(\B(a,b,x) \land \B(a,c,y) \land \B(x,t,y) )$&
\endaxioms

Something like this axiom was first considered by Legendre (see \cite{greenberg}, p. 223), but 
he required angle $bac$ to be acute, so Legendre's axiom is not exactly the same as Tarski's 
parallel axiom.  The axiom says a bit more than just that $xy$ meets both sides of the angle,
because of the betweenness conditions in the conclusion; but it would be equivalent to 
demand just that $x$ and $y$ lie on the rays forming angle $bac$, as can be shown using 
Pasch.

\subsection{Euclid 5 implies Tarski's parallel axiom}

\begin{Theorem} \label{theorem:8.2}
Euclid 5 implies Tarski's parallel axiom in neutral
intuitionistic Tarski geometry.
\end{Theorem}
\medskip

\begin{figure}[h]
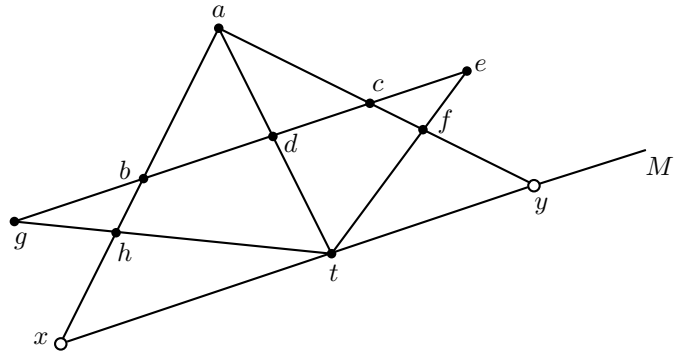
   
\caption{Constructive proof of Tarski's parallel axiom from Euclid 5. $M$ is constructed
parallel to $\Line(b,c)$ and $cd=ce$ and $bd=bg$. Then $x$ and $y$ exist by Euclid 5. }
\label{figure:TarskiParallelProofFigure}  
\TarskiParallelProofFigure
\end{figure}

\medskip

\noindent{\em Proof}. Assume the hypothesis of Tarski's parallel axiom.  Construct line $M$
parallel to $\Line(b,c)$ through $t$.  Construct point $e$ by extending segment $dc$ by $dc$;
then $ec = dc$ and $\B(d,c,e)$, as
 illustrated in Fig.~\ref{figure:TarskiParallelProofFigure}.   Let $L$ be $\Line(a,c)$.
 Then $\Line(d,t)$ meets $L$ at $a$, and does not coincide with $L$ since, if it did 
 coincide with $L$, then points $d$ and $c$ would be on $L$, and hence point $b$,
 which is on $\Line(b,c)$, would lie on $L$ by Axiom I3; but that would contradict the 
 hypothesis that $a$, $b$, and $c$ are not collinear.   Hence $\Line(d,t)$ meets $L$ 
 only at $a$, by Axiom I3. Hence segment $dt$ does not meet $L$. By outer Pasch 
 (applied to $adtec$),
 there is a point $f$ on $L$ with $\B(e,f,t)$.   Now we apply Euclid 5; 
 the two parallel lines are $\Line(b,c)$ and $M$, and the conclusion is that $L$ meets 
 $M$ in some point, which we call $y$.   Specifically we match the variables $(L,K,M,p,r,a,q)$
 in Euclid 5 to the following terms in the present situation: $(M,\Line(b,c),\Line(a,c),c,e,f,t)$.
 Then all the hypotheses have been proved, except that we have $\B(e,f,t)$ while what is 
 required is $\B(t,f,e)$; but those are equivalent by the symmetry of betweenness.  Hence 
 Euclid 5 is indeed applicable and we have proved the existence of point $y$ on $M$ and $L$.
 
 Now, we do the same thing on the other side of angle $bac$,  extending segment $db$ to 
 point $g$ with $db=bg$ and $\B(g,d,b)$, and using the plane separation property to
 show that $gt$ meets $\Line(a,b)$ in a point $h$ with $\B(g,h,t)$.   Then Euclid 5 
 applies to give us a point $x$ on $M$ and $\Line(a,b)$. 

 It only remains to prove $\B(x,t,y)$.  By outer Pasch applied to $xbacd$, there 
 exists a point $u$ with $\B(x,u,c)$ and $\B(a,d,u)$ (the point $u$ is not shown 
 in the figure.)   Then by outer Pasch applied 
 to $acyxu$, we obtain a point $v$ with $\B(a,u,v)$ and $\B(x,v,y)$.  But then $v=t$,
 since both lie on the non-coincident lines $ad$ and $xy$.  Hence $\B(x,t,y)$.
   That completes the proof of the theorem.
 \smallskip
 
 \subsection{Tarski's parallel axiom implies Euclid 5}
 Tarski proved%
 \footnote{The cited proof is in a book with two co-authors, but Tarski used this 
 axiom from the beginning of his work in geometry,
  and it seems certain that he had this proof before Szmielew
 and Schwabh\"auser were involved.}
  that his parallel axiom implies Playfair's axiom
 (see \cite{schwabhauser}, Satz~12.11, p.~123).  Here we give a constructive 
 proof that Tarski's parallel axiom implies the points-only version of Euclid 5.
 See Fig.~\ref{figure:TarskiImpliesEuclidFigure}. 
\medskip

\begin{figure}[h]   
\caption{Tarski's parallel axiom implies Euclid 5.   Points $x$ and $y$ are produced
by Tarski's parallel axiom because $q$ is in the interior of angle $upa$.
Then apply Pasch to line $L$ and triangle $xpy$ to get $e$.}
\label{figure:TarskiImpliesEuclidFigure}  
\TarskiImpliesEuclidFigure
\end{figure}

\medskip
 
 \begin{Theorem} \label{theorem:8.3}
  Tarski's parallel axiom implies Euclid 5 in neutral 
 intuitionistic Tarski geometry.
 \end{Theorem}
 
\noindent{\em Proof}.  
Let $L$ be a line, $p$ a point not on $L$,  
$M$ be another line through $p$, and suppose points $p$, $r$, $s$,
and $q$ are as in the hypothesis of Euclid 5 (see Fig.~\ref{figure:TarskiImpliesEuclidFigure}).
In particular, because the shaded triangles are congruent, $K$ is parallel to $L$ and 
makes alternate interior angles equal with transversal $pq$. 

Let $u$ be a point 
to the left of $p$ on $K$, for example $u = ext(r,p,\alpha,\beta)$.  We 
can apply inner Pasch to the configuration $uprqa$, producing a point $v$ such 
that $\B(p,v,q)$ and $\B(u,v,a)$. (This is where we use the hypothesis $\B(r,a,q)$.)
 Then $v$ witnesses that $q$ is in the 
interior of angle $upa$. By 
Tarski's parallel axiom, there exist points $x$ and $y$ with $\B(x,q,y)$, $\B(x,u,p)$ and 
$\B(p,a,y)$.  Then line $L$ meets side $xy$ of triangle $xpy$ at $q$,
and does not meet the closed side $xp$, since $K$ is parallel to $L$.
Then $x$ and $y$ are on opposite sides of $L$, and $x$ and $p$ are on the same
side of $L$.  By the Plane Separation Theorem (Theorem~\ref{theorem:planeseparation}),
 $p$ and $y$ are on 
opposite sides of $L$.   
 
We next wish to prove $\B(p,a,e)$. By the stability of betweenness, we may 
prove it by contradiction.   By hypothesis, $r$ does not lie on $L$, so $a$ does not lie on $L$.
Hence $a \neq e$.  Since $L$ and $K$ are parallel, $e \neq p$.  Suppose $\B(p,e,a)$.
Then we can apply outer Pasch to $parqe$ (that is, to triangle $par$ with sequent $qe$),
obtaining a point lying both on $Line(q,e)$ (which is $L$) and $pr$ (which is part of $K$),
contradiction.    Hence $\neg \B(p,e,a)$.
Now suppose $\B(e,p,a)$. Then we can apply outer Pasch to triangle $eaq$ with secant $rp$,
obtaining a point on segment $eq$ (part of $L$) and also on $Line(p,r)$, which is $K$;
contradiction.  Also $e$ cannot be equal to $a$, since then $r$ would lie on $L$ and $K$;
and $e$ cannot be equal to $p$ since then $e$ would lie on both $K$ and $L$.  The only 
remaining possibility is $\B(p,a,e)$, which we set out to prove by contradiction.
That completes the proof of $\B(p,a,e)$.

 We still must show $\B(s,q,e)$.  Since we now have $\B(p,a,e)$, 
we can apply outer Pasch to $paest$ to construct a point $q^\prime$ such that 
$\B(p,t,q^\prime)$ and $\B(s,q^\prime,e)$.  Then $q^\prime$ lies on both $Line(p,t)$
and $L$ (which are distinct lines since $p$ does not lie on $L$); but $q$ also lies
on both those lines.  Hence $q^\prime = q$.  Hence $\B(s,q,e)$ as desired.
That completes the proof of the theorem.
 
\section{Uniform perpendicular and uniform rotation}

In classical geometry there are 
two different constructions,  one for ``dropping a perpendicular'' to line $L$ from a point 
$p$ not on $L$, and the other for ``erecting a perpendicular'' to $L$ at a point $p$ on $L$.
A ``uniform perpendicular'' construction is a method of constructing, given three points
$a$, $b$, and $x$,  with $a \neq b$,   a line perpendicular to $\Line(a,b)$ and passing through $x$, without a 
case distinction as to whether $x$ lies on $L$ or not. 

In constructive geometry, it is not sufficient to have dropped and erected perpendiculars;
we need uniform perpendiculars.  Similarly, we need uniform rotations: to rotate a given point
$x$ on $\Line(c,a)$ clockwise about center $c$ until it lies on a given line through $c$, without a case 
distinction as to whether $\B(x,c,a)$ or $x=c$ or $\B(c,x,a)$.  We also need 
uniform reflections:  to be able to reflect a point $x$ in a line $L$ without a case
distinction whether $x$ is on $L$ or not.  It will turn out (not surprisingly) 
that the three problems 
of perpendiculars,  rotations, and reflections are closely related.

\subsection{Uniform perpendicular using line-circle}

In this section, we take up the construction of the uniform perpendicular.  
 
\begin{Theorem} Uniform perpendiculars can be constructed in intuitionistic Tarski
geometry, using two-point line-circle continuity.
\end{Theorem}

 \begin{figure}[ht]
\caption{The simplest uniform perpendicular construction.  $M=Perp(x,L)$ is constructed 
perpendicular to $L$ without a case distinction whether $x$ is on $L$ or not.  Draw
a large enough circle $C$ about $x$.  Then bisect $pq$ and erect $K$ at the midpoint.
To draw $C$ we use radius $ab + ax > ax$.
}
\label{figure:UniformPerpendicularFigure}
\hskip 2cm    
\psset {unit=2cm}
\PerpFigure
\psset {unit=3cm}
\end{figure}

\noindent{\em Proof}.  Given distinct points $a$ and $b$ (defining $L = Line(a,b)$),
and point $x$ (without being told whether $x$ is or is not on $L$),  we desire to 
construct a line $K$ passing through $x$ and perpendicular to $L$.
  The idea is simple:  Draw a 
 circle $C$ around $x$ whose radius exceeds $ax$.  Then by two-point line-circle continuity
 there are distinct points $p$ and $q$ on $L = \Line(a,b)$ that lie on $C$.
 Then the perpendicular bisector of segment $pq$ is the desired perpendicular.
 The matter is, however,  trickier than the similar verification of Euclid~I.2, because
of the requirement to find a point $r$ through which to draw the circle $C$.
We must have $xr > ax$,  but we are not allowed to make a case distinction whether 
$x = a$ or not, and since we are not allowed to extend a null segment, we cannot 
find $r$ by extending the (possibly null) segment $ax$.   First we construct
point $c$ by extending segment $ab$ by $ax$.
Since $ab$ is not null, and we are allowed to extend a non-null segment by a possibly null 
segment, this is legal.  Then we make use of the construction $e_2$ corresponding to 
Euclid~I.2,  developed in Lemma~\ref{lemma:EuclidI.2}.  We define $r = e_2(x,a,c)$.
Then $xr = ac$ and since $ac > ax$ we have $ar > ax$.  Now draw the circle $C$
with center $x$ through $r$.  Point $r$ is not shown in the figure, since its 
exact location depends on the construction given by $e_2$, which depends on the constant $\alpha$.
All we know about $r$ is that $xr=ac$, so we can use it to draw $C$ with the desired radius.

Then $a$ is inside $C$, so $p$ and $q$ exist
by two-point line-circle continuity, and we can complete the proof using any construction 
of the perpendicular bisector,  for example the one developed using two-point line circle
continuity only, or Gupta's more complicated one.

It is not difficult to construct a script, and hence in principle a term, describing 
this construction.  This term does require an ``extra'' parameter $s$,  for a point
assumed to be not on the line $L$,  in order to erect the required perpendicular.

\subsection{Uniform perpendicular without any continuity}  
The above construction has one disadvantage:  it relies on line-circle continuity.
Although our main interest is in intuitionistic Tarski geometry,  nevertheless 
it is of some interest to study the theory that results from deleting line-circle 
continuity as an axiom, i.e.,the intuitionistic counterpart of Tarski's (A1)-(A10).
It turns out, one can also construct a uniform perpendicular in that theory,
although of course, one must use Gupta's perpendiculars.  It turns out that we 
also need the parallel axiom to construct uniform perpendiculars, although
it is not needed for Gupta's perpendiculars.

We begin with some lemmas.%
\footnote{Though these are fairly routine exercises in Euclidean geometry,
we need to verify that they are provable in intuitionistic Tarski geometry without
any continuity.  These theorems do not occur in \cite{schwabhauser}, and even if they 
had occurred, we would still need to check that they are provable with triangle circumscription
instead of Tarski's parallel axiom.  In fact we will use Euclid 5, which is 
provable from triangle circumscription, as shown in \cite{beeson-bsl}.}  
 
\begin{Lemma}\label{lemma:saccheri} Let $abcd$ be a quadrilateral with two 
adjacent right angles at $a$ and $d$, and two opposite sides equal, namely $ab = cd$.
Suppose $a$ and $d$ are on the same side of $Line(b,c)$.  
Then also the other pair of opposite sides are equal and $abcd$ is a rectangle.
\end{Lemma}

\noindent{\em Remarks}.  The proof necessarily will require the parallel axiom, as 
the lemma is false in hyperbolic geometry.  The assumption that $a$ and $d$ are 
on the same side of $Line(b,c)$  
 ensures that the quadrilateral lies in a plane (we do not use any dimension
axiom in the proof). 
\medskip

\noindent{\em Proof}.  We first claim that $c$ is on the same side of $ad$ as $b$.
If not then by the plane separation theorem, they are on opposite sides, so 
$bc$ meets $Line(a,d)$ in a point $m$.  By Lemma~\ref{lemma:saccheri-helper},
$\B(a,m,d)$, contradicting the hypothesis that $a$ and $d$ are on opposite sides of 
$Line(b,c)$.  That establishes that $c$ is on the same side of $ad$ as $b$.

Let $J$ 
be the perpendicular bisector of $ad$, meeting $ad$ at its
midpoint $q$.  Then $q$ is the foot of the perpendicular from $p$ to $J$.
Since reflection in $J$ preserves congruence, collinearity, and right angles,
the reflection $b^\prime$ of $b$ lies on the perpendicular to $ad$ at $d$,
which is $cd$, and either $b^\prime = c$ or $b^\prime$ is the 
reflection of $c$ in $d$.   The latter, however, contradicts
the fact that   $c$ is on the same side of $ad$ 
as $b$.  Therefore $c$ is the reflection of $b$ in $J$.  Hence $bc$
is perpendicular to $J$. 

 We next claim that $ad$ and $bc$ are parallel; that is,
 $Line(a,d)$ and $Line(b,c)$ do not meet.
  Suppose they  meet at point $p$.   Then point $p$ is not between 
   $a$ and $d$, since $a$ and $d$ would then be on opposite sides
of $bc$.   Interchanging $(a,b)$ and $(d,c)$ if necessary, we can assume
without loss of generality that $\B(p,a,d)$.   Now $pbc$ and $pad$ are
two distinct lines through $p$, both perpendicular to $J$.  That is 
a contradiction (as proved in
the first part of Satz~8.18).  That   
establishes the claim that $ad$ and $bc$ are parallel.

Also $ab$ and $cd$ are parallel, since if they meet at a point $u$ then 
there would be two perpendiculars from $u$ to $\Line(a,d)$.  Hence 
$abcd$ is a parallelogram.    

Using the parallel postulate, we can prove that the angles at $b$ and $c$
are right angles:  since $bc$ is parallel to $ad$, if the angle at $b$
is not a right angle, then the perpendicular to $ab$ at $b$ is a 
second, distinct, parallel to $ad$ through $b$, contradicting the parallel 
postulate (see \cite{beeson-bsl}).   Then triangles $abd$ and $cbd$ are 
right triangles with one leg and the hypotenuse congruent; hence 
their other legs $ad$ and $bc$ are congruent.  That completes the proof
of the theorem.

\begin{figure}[ht]
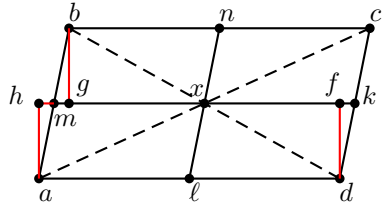

\caption{Given that $ac$ and $bd$ bisect each other at $x$.  Then the lines that 
look parallel, are parallel.}
\label{figure:parallelograms}
\ParallelogramsFigure
\end{figure}

\begin{Lemma} \label{lemma:parallelograms}
(in intuitionistic Tarski geometry 
without continuity) 
Let $abcd$ be a quadrilateral whose diagonals $ac$ and $bd$ bisect each other at $x$.
Then (i) opposite sides of $abcd$ are parallel, and (ii) the lines connecting midpoints of 
opposite sides pass through $x$, and (iii) they are parallel to the other sides, as 
shown in Fig.~\ref{figure:parallelograms}.  
\end{Lemma}

\noindent{\em Proof}.  
  Since reflection (in point $x$) 
preserves betweenness and congruence (see Satz~7.15 and Satz~7.16 of \cite{schwabhauser}),
we have $ab = cd$ and $bc = ad$, i.e.,opposite sides are equal.  Let $m$
be the midpoint of $ab$ and $k$ the reflection of $m$ in $x$.  Since 
reflection preserves congruence, $ck = kd$.  Since 
reflection in a point also preserves betweenness,
$k$ lies on segment $cd$.  Hence $k$ is the midpoint of $cd$.
Let $n$ be the midpoint of $bc$ and $\ell$ the midpoint of $ad$. Then 
similarly $n\ell$ passes through $x$.
That proves that the lines connecting midpoints of opposite sides pass through $x$
as claimed in (ii) of the theorem. 

Suppose point $u$ lies on $\Line(a,d)$ and also on $\Line(b,c)$.  Let $v$
be the reflection of $u$ in $x$.  Then since  reflection preserves collinearity, 
by Lemma~\ref{lemma:reflectionB}, $v$ also lies on both lines.  But $v \neq u$,
since if $v=u$ then $v=x=u$, but $m$ does not lie on $\Line(a,d)$.  Hence there
are two distinct points $u$ and $v$ on $\Line(a,d)$ and $\Line(b,c)$, contradiction.
Hence those two lines are parallel, as claimed.  Similarly, the other two sides
of the quadrilateral are parallel.  That proves part (i) of the theorem.

By the inner five-segment lemma (Lemma~\ref{lemma:interior5})  applied to $ambx$ and $ckdx$,
we have $mx = xk$.  Therefore triangle $ckx$ is congruent to triangle $amx$.
Then triangle $ncx$ is congruent to triangle $\ell ax$,
since they have all three pairs of corresponding sides equal.

We next claim that $mk=ad$.  By the stability of equality, we may argue 
by contradiction and cases.   Drop a perpendicular from $d$ 
to $mk$; let $f$ be the foot.  Case 1, $f=k$. Then $dk \perp mk$, and by 
reflection in $x$, $bm \perp mx$, so $ma \perp mk$ and $mkda$ has two adjacent right 
angles. Hence its opposite sides $mk$ and $ad$ are equal, by Lemma~\ref{lemma:saccheri},
since $a$ and $d$ are on the same side of $mk$ because $ab$ and $bd$ both meet $\Line(m,k)$
(in $m$ and $x$ respectively).
Case 2, $\B(m,f,k)$.  Then let $g$ be the reflection of $f$ in $x$, and $h$
the reflection of $g$ in $m$.  Then the triangles $dfk$, $bgm$, and $ahm$
are congruent right triangles (by reflection), and $\B(h,m,g)$.  Then 
$hm=mg=fk$.  Since $hadf$ has two adjacent right angles, its opposite sides $hf$ 
and $ad$ are equal, since again $a$ and $d$ are on the same side of $\Line(h,f)$,
so Lemma~\ref{lemma:saccheri} applies.
  But $hf = mk$, since they differ by adjoining and removing
equal segments $hm$ and $fk$.  Hence $mk=hf=ad$, so $mk=ad$ as claimed.  
Case 3, $\B(m,k,f)$.  Reflection in the line $mk$ reduces this case to Case 2.
That completes the proof that $mk=ad$.  

Now $mkda$ is a parallelogram with opposite sides equal.  Therefore, part (i)
of the theorem, which has already been proved, can be applied to it.  Hence 
$mk$ is parallel to $ad$ as claimed.  That completes the proof.

\begin{Lemma}\label{lemma:twoperps}  (in intuitionistic Tarski geometry 
without continuity) Let line $J$ be parallel to line $M$, and 
suppose $M \perp L$.   Then $J$ meets $L$ in a point $x$ and $J$ is perpendicular 
to $L$ at $x$.
\end{Lemma}

\begin{figure}[ht]
\caption{Given $M \perp L$ and $J$ parallel to $M$, show $J \perp L$ and 
construct the intersection point $x$ of $J$ and $L$ by dropping a perpendicular
from $p$ to $J$. Construct $r$, $t$, and $s$ by reflection. Then $ts$ is parallel to $M$
and hence lies on $J$, and $ts \perp L$.}
\label{figure:TwoPerpsFigure}
\TwoPerpsFigure
\end{figure}

\noindent{\em Proof}. See Fig.~\ref{figure:TwoPerpsFigure}.
Let $p$ be the point of intersection of $M$ and $L$.
Drop a perpendicular from $p$ to $J$; let $x$ be the foot of that perpendicular.
Let $m$ be the midpoint of $px$ (using Gupta's midpoint to avoid any need for 
continuity). Let $q$ be any point other than $p$ on $M$ and let $r$ be 
the reflection of $q$ in $p$.  Then $mq = mr$ since $M \perp L$.  Let $s$ and $t$
be the reflections of $q$ and $r$ in $m$, respectively.
Since reflection preserves congruence, $sx = pq$.  Then 
$xs$ is parallel to $M$, since if point $u$ on $M$ lies on $\Line(x,s)$,
then the reflection of $u$ in $m$ is another point on both lines, contradiction.
Now $J$ is by hypothesis parallel to $M$, and $J$ passes through $x$.
  But by the parallel axiom, there cannot be two 
parallels to $M$ through $x$.    Hence  the two coincide:  $J = \Line(x,s)$.
We have $ms = mq$ by reflection, $mq = mr$ since $M \perp L$, 
and $mr = mt$ by reflection; hence $ms = mt$.  Then $Line(p,x) \perp J$ by the 
definition of perpendicular.   We now claim that $Line(p,x) = L$; for 
that we shall need the dimension axiom, since otherwise $L$ might be another perpendicular
to $M$ in the plane perpendicular to $M$ at $p$.  We have $mq = mr$, so $Line(p,x) \perp M$.
Then $xq = xr$ also.   In order to show that $L$ coincides with $Line(p,x)$, suppose 
$w$ lies on $L$.  Then since $L \perp M$, we have $wq = wr$.  Then $w$, $x$, and $m$ 
are three points, each equidistant from $q$ and $r$.  By the upper dimension axiom 
these three points are collinear.  Hence $w$ lies on $Line(m,x)$, which coincides with
$Line(p,x)$.   Hence $L$ coincides with $Line(p,x)$.  Hence $L \perp J$.  Hence $J \perp L$.
That completes the proof.

\begin{Theorem}\label{theorem:uniformperp}
 
[Uniform perpendicular] Uniform perpendiculars can be constructed in 
intuitionistic Tarski geometry without any continuity axioms.  Explicitly:
there is a term $Project(x,a,b,w)$ (``the projection of $x$ on $Line(a,b)$'') 
in intuitionistic Tarski with Skolem functions,
such that if $a \neq b$ and $w$ is any point not on $L = \Line(a,b)$, and $f = Project(x,a,b,w)$,
then $Col(a,b,f)$, and the erected perpendicular to $L$ at $f$ contains $x$.
 
\end{Theorem}

\noindent{\em Remark}. The two main points of the lemma are that we do not need a case 
distinction whether $x$ is on $L$ or not, and we do not use any continuity axiom.
But we do use the parallel axiom.  
\medskip

\begin{figure}[ht]
\caption{Uniform perpendicular.  Given $L$ and $x$ construct perpendicular $J$ 
to $L$ passing through $x$, without a case distinction whether $x$ is on $L$ or not.}
\label{figure:UniformPerpFigureTwo}
\UniformPerpFigureTwo
\end{figure}

\noindent{\em Proof}. See Fig.~\ref{figure:UniformPerpFigureTwo}.
 We begin by constructing a point $p$ on $L$ such that $px$
is not perpendicular to $L$ (unless of course $x$ lies on $L$, which is a 
degenerate case). The existence of such a point $x$  is trivial classically,
but constructively, there is something to prove.)
 The line $L$ is given by two points $a$ and $b$;
let $p$ be the result of extending $ab$ by $ab$ to point $r$ and 
then extending $ar$ by $xb$ to point $p$.  Then $\B(b,r,p)$ and $rp = xb$.
(The idea is that, if $a$ is to the left of $b$, then $p$ is far enough to 
the right that $x$ cannot lie on the perpendicular to $L$ at $p$.) Now suppose
that $px \perp L$; we will show that $x$ lies on $L$.  By the stability of incidence, 
we may prove that by contradiction. 
If $x$ does 
not lie on $L$, then 
$bpx$ is a triangle, and  it has leg $bp$ greater than the hypotenuse $bx = rp$,
contradiction.  Therefore $x$ is on $L$ as claimed.
Then $px$ lies on $L$ too, and $p \neq x$;  
hence $px$ is not perpendicular to $x$.  That completes the preliminary construction of $p$.
Note that this part of the proof is not shown in the diagram; but it is needed to 
show that the diagram does not degenerate to a single vertical line instead of a 
parallelogram.

Erect a perpendicular $qp$ to $L$ at $p$. (For that we need the point $w$ not on $L$.)
  Let $s$ be the reflection of $q$ in $x$.
Let $t$ be the reflection of $p$ in $x$. Then $tx = xp$ and $sx = xq$.
By Lemma~\ref{lemma:parallelograms}, 
$stqp$ is a parallelogram whose diagonals bisect each other. 
Let $m$ be the midpoint of $tq$ and $k$ the midpoint of $sp$.  By 
Lemma~\ref{lemma:parallelograms}, $tq$ is parallel to $sp$, so $m \neq k$.
 We claim $J = \Line(k,m)$ is the desired perpendicular to $L$.  By Lemma~\ref{lemma:parallelograms},
 $J$ is parallel to $qp$.  But
 $qp \perp L$ by construction.
By the parallel axiom, lines parallel to $qp$ are also 
perpendicular to $L$; so $J \perp L$.  By Lemma~\ref{lemma:twoperps},
the intersection point $f$ of $J$ and $L$ exists and $J \perp L$.  The point $f$
is the value of $Project(x,a,b,s)$.   That completes the proof.  
\medskip

\noindent{\em Remark}.  We do not know how to construct uniform perpendiculars
without using either the parallel axiom or two-point line-circle, although 
either one suffices, and any form of the parallel axiom suffices (because we just 
need a few simple lemmas about parallelograms). 
\medskip

\subsection{Uniform reflection}
Another construction from \cite{beeson-bsl} that we need to check can be done with 
Tarski's Skolem symbols is the ``uniform reflection''.   The construction $$\Reflect(x,a,b,s)$$
gives the reflection of point $x$ in $L = \Line(a,b)$, without a case distinction as to whether 
$x$ is on $L$ or not.  (The parameter $s$ is some point not on $L$.)
 First we note a difficulty:  even though we can define 
$f = Project(x,a,b,s)$,  we cannot just extend segment $xf$ by $xf$, since 
$xf$ might be a null segment, and in constructive geometry, we can 
only extend non-null segments.  

The solution given in \cite{beeson-bsl} is to first define
rotations, and then use the fact that the reflection of $x$ in $\Line(a,b)$ is the same as
the result of two ninety-degree rotations of $x$ about $f = Project(x,a,b,s)$.
The construction given for rotations in \cite{beeson-bsl} only involves bisecting the angle 
in question, and dropping two perpendiculars to the sides, none of which is problematic in 
Tarski's theory.
\medskip

\subsection{Equivalence of line-circle and segment-circle continuity}

\begin{Lemma} \label{lemma:line-segment-equiv}
Two-point line-circle, one-point line-circle, and segment-circle continuity are equivalent
in (A1)-(A10).  
\end{Lemma}

\noindent{\em Remark}. This proof depends heavily on Gupta's 1965 dissertation. 
\medskip

\noindent{\em Proof}.  Because of Gupta, we have dropped perpendiculars, and we 
have shown above that from dropped perpendiculars, erected perpendiculars, midpoints,
and uniform perpendiculars follow.  Then one-point line-circle implies two-point line-circle, 
by reflection in the uniform perpendicular from line $L$ through the center of the circle.
Note that classically, the case when $L$ passes through the center is trivial, so a dropped
perpendicular is enough, and the parallel axiom is not required.  

Two-point line-circle implies segment-circle immediately.   

Next we will show that segment-circle implies one-point line-circle.  Assume 
the hypotheses of one-point line-circle:  
Let $C$ be a circle with center $a$ through $b$, and let $L$ be a line meeting
closed segment $ab$ at $p$.  The desired conclusion is the existence of a point 
$x$ on $L$ with $ax=ab$.  To obtain that from segment-circle, it suffices to 
construct a point $q$ on $L$ outside $C$ (see Fig.~\ref{figure:LineCircleContinuityFigureTwo}
and Fig.~\ref{figure:SegmentCircleContinuityFigure}). (Then segment-circle will give
us a point on $C$ between $p$ and $q$; since $p$ and $q$ are on $L$, that point will be on $L$ too.)

Inequality is defined between segments as follows:
$$uv < xy \mbox{\quad means \quad} \exists z\,(\T(x,z,y) \land uv = xy \land z \neq y$$
and
$$uv \le xy \mbox{\quad means \quad} \exists z\,(\T(x,z,y) \land uv = xy. $$
Both strict and non-strict inequality can be shown to be stable, e.g.,
$$ \neg \neg  uv < xy \implies uv < xy, $$
so inequalities can be proved by contradiction.

Here is how to construct $q$.   Let $z$ be any point on $L$
different from $p$, and define $q$ by extending segment $zp$ by three times the radius $ab$.
Then $qa  > ab-ap \ge qp -ab$, by the triangle inequality, Lemma~\ref{lemma:triangleinequality}.
Then we have 
$$ qa > qp - ab > 3ab - 2ab > ab.$$
This apparently algebraic calculation represents a geometric proof
of $qa > ab$. 
In view of the definition of segment inequality, this implies the existence of a 
``witness'' point $y$ such that $\B(y,b,a)$ and $ay = aq$.  This point and $q$ then 
satisfy the hypothesis of segment-circle continuity (see 
 Fig.~\ref{figure:SegmentCircleContinuityFigure}).
Segment-circle continuity then 
yields a point on $C$ between $q$ and $p$, which is therefore on $L$.

We have now proved a circular chain of implications between the three assertions of the lemma.
  That completes the proof.
\medskip

\noindent{\em Remark.}  One may wonder why Tarski chose segment-circle rather 
than line-circle continuity as an axiom.  It might be because segment-circle
continuity asserts the existence of something that turns out to be unique; but 
that consideration did not bother Tarski when he chose (A10) as his parallel axiom.
More likely it was just shorter.  
 Although the diagram for (either form of) line-circle continuity is simpler, the formal
expression as an axiom is longer, especially if collinearity is written out 
rather than abbreviated; and Tarski placed importance on the fact that his 
axioms could be written intelligibly without abbreviations.  
\medskip

\subsection{Representation theorems}
The following important theorem was stated in 1959 by Tarski  \cite{tarski1959}:

\begin{Theorem} \label{theorem:representation}
(i) The models of classical Tarski geometry with segment-circle continuity
are the planes $F^2$ over a Euclidean field $F$. 
\smallskip

(ii) The models of classical Tarski geometry with no continuity axioms are the 
planes $F^2$ over a Pythagorean field $F$.%
\footnote{A Pythagorean field is one in which sums of two squares always have
square roots, or equivalently, 
$\sqrt{1+x^2}$ always exists, as opposed to Euclidean in which all positive elements 
have square roots.}
\end{Theorem}

\noindent{\em Remark}.   Tarski wrote ``ordered field'' instead of ``Pythagorean field''
in (ii),  but one needs to be able to take $\sqrt{1+x^2}$ in $\F$ to verify the 
segment extension axiom in $\F^2$,  as Tarski surely knew.
\medskip

There is also a version of the representation theorem for intuitionistic
Tarski geometry.  For that to make sense, we must define ordered fields and 
Euclidean fields constructively.  That is done in \cite{beeson-bsl}, as follows:   We 
take as axioms the stability of equality and ``Markov's principle'' that 
$\neg x \le 0$ implies $x>0$, which is similar to the stability of betweenness
in geometry.  Then we require that nonero elements have multiplicative 
inverses; just as classically, a Euclidean field is an ordered field in which positive elements
have square roots (and a Pythagorean field is one in which sums of two squares have square roots).
Then we have a constructive version of 
Tarski's representation theorem:  

\begin{Theorem} \label{theorem:i-representation} 
(i) The models of intuitionistic Tarski geometry with two-point line-circle continuity
are the planes $F^2$ over a Euclidean field $F$.  
\smallskip

(ii) The models of intuitionistic Tarski geometry with no continuity axioms are the 
planes $F^2$ over a Pythagorean field $F$.
\smallskip

(iii) Given a model of geometry, the field $F$ and its operations can be 
explicitly and constructively defined. 
\end{Theorem}

\noindent{\em Proof (of both theorems)}.  Once we have (uniform) perpendiculars and midpoints, the classical 
constructions of Descartes and Hilbert that define addition and multiplication can 
be carried out.  In this paper we have proved the existence of (uniform) perpendiculars and 
midpoints in intuitionistic Tarski geometry, and the definitions of (signed uniform) addition and 
(signed uniform) multiplication are given in \cite{beeson-bsl}.
The field laws are proved for these definitions in \cite{schwabhauser}; since these are 
quantifier free when expressed with Skolem functions, they are provable in intuitionistic
Tarski geometry (without any continuity axiom) 
by the double-negation interpretation.  The simple construction 
of the uniform perpendicular suffices for (i), with two-point line-circle, 
 but for (ii) we need the more complicated
construction given above, based on Gupta's perpendiculars.  Moreover, even for 
Tarski's (i), with segment-circle in place of two-point line-circle, we need Gupta.
 That completes the proof.
 \medskip
 
With classical logic, the representation theorem gives a complete characterization of the 
consequences of the axioms, because according to G\"odel's completeness theorem, a sentence
of geometry true in all models $F^2$ is provable in the corresponding geometry.  
With intuitionistic logic, the completeness theorem is not valid.  The ``correct'' way to 
obtain a complete characterization of the theorems of constructive geometry in terms of 
field theory is to exhibit explicit interpretations mapping formulas of geometry to 
formulae of field theory,  and an ``inverse'' interpretation from field theory to geometry.
We have actually carried this program out, but it is highly technical due to the differences
in the two languages, and requires many pages, so we omit it here.
 
\subsection{Historical Note}   
Tarski claimed in 1959
\cite{tarski1959} (page 22, line 5)
that (i) he could define addition and multiplication geometrically,
and (ii) prove the field laws,  without using his continuity axiom; hence all models 
of the axioms (excluding continuity) are planes over ordered fields.  The first 
published proof of these claims was in 1983 \cite{schwabhauser}, and relies
heavily on Gupta's 1965 dissertation.  In this note we consider what Tarski
may have had in mind in 1959.   

To define multiplication, we need perpendiculars and midpoints, which 
were constructed in (A1)-(A10) by Gupta in 1965.   Tarski lectured in 1956
on geometry, but I have not found a copy of his lecture notes.  He lectured 
again on geometry in 1963, but according to Gupta (in a telephone conversation 
in September, 2014), Tarski did not base his lectures on his 
own axiom system (and again there are apparently no surviving notes).  The proofs in this section 
show that he could well have defined addition and multiplication using only two-point line-circle 
continuity, since  we showed here how to construct perpendiculars and midpoints from two-point
line-circle continuity.  But he made two claims that we do not see how to prove 
without Gupta:  that he could use segment-circle continuity in (i), and that he could get 
by without any continuity in (ii).

These claims could not have been 
valid in 1959,  since this was some years before Gupta's proofs.  In 1959, 
there was no way known to construct dropped perpendiculars using (A1)-(A10), 
even with the use of segment-circle continuity;  Tarski could not have justified 
Euclid I.12 on the basis of segment-circle continuity,
 since it requires the triangle inequality, which requires perpendiculars,
and Tarski was not using the triangle circumscription axiom but his own form of the 
parallel axiom, so he could not even prove that every segment has a midpoint, as far 
as I can see.    

Tarski did believe that line-circle continuity could be derived 
from a single instance of the continuity schema.  
He  explicitly claimed this in \cite{tarski1959}, page~26, line 8.  But   
the ``obvious'' derivation requires the triangle inequality, 
which in turn seems to require perpendiculars.%
\footnote{It is straightforward to derive line-circle continuity from A11 together with 
the facts that the interior and exterior of a circle are open, and the interior is convex.
These facts in turn are easy to derive from inner Pasch,  Euclid III.2 (chord lies inside circle),
the density lemma (Lemma~\ref{lemma:density}), and the triangle inequality.}
After Gupta, perpendiculars exist, and the proof of the 
triangle inequality follows easily, so indeed line-circle continuity follows from A11.
But in 1959,  there was no known way to get perpendiculars, and so, no way 
to derive the triangle inequality, and hence, no way to derive
line-circle continuity from the continuity schema, or to justify Euclid~I.12 to get 
dropped perpendiculars from line-circle continuity.  All these difficulties disappeared
once Gupta proved the existence of dropped perpendiculars in (A1)-(A8).  Thus Gupta's
work has a much more central place than is made apparent in \cite{schwabhauser}.  Tarski
desperately needed perpendiculars.  

It is possible that Tarski had in mind using two-point line-circle continuity to 
justify dropped perpendiculars,  
and overlooked the difficulty of proving two-point line-circle from segment-circle.  Even so, to 
complete a proof of his representation theorem about Euclidean fields, he would have had to duplicate 
the results in this paper about getting erected perpendiculars and midpoints from 
dropped perpendiculars.  Part (ii) of his theorem, about what happens with no continuity at all,
flat-out requires Gupta's work,  which put everything right, substantiating 
the claims of 1959.   A discovery of Tarski's missing 
1956 lecture notes seem to be the only way to resolve the question of ``what Tarski knew
and when he knew it.''
\medskip

\section{Geometry with terms for the intersections of lines}
It seems more natural, when thinking of straightedge-and-compass constructions,
to include a symbol $\il(a,b,c,d)$ for the (unique) intersection point of 
$\Line(a,b)$ and $\Line(c,d)$.  We say ``unique'' because we want the intersection 
point of two coincident lines to be undefined;  otherwise $\il(a,b,c,d)$ will not 
be continuous on its domain.

The difficulty with using this Skolem symbol
is that the definedness condition for $\il(a,b,c,d)$ is not easily expressible
in quantifier-free form.   Of course we need $a \neq b \land c \neq d$, and 
we want $\neg (Col(a,b,c) \land Col(a,b,d)) $ as just explained.  But in addition 
there are parallel lines that do not meet.  Using the strong parallel postulate,
  one can indeed express the definedness condition for $\il(a,b,c,d)$ in a 
  quantifier-free way, namely, $\il(a,b,c,d)$ is defined if and only if there
  is a point $p$  collinear with $a$ and $b$ but not $c$ and $d$, and a point 
  $q$ collinear with $c$ and $d$,  such that the transversal $pq$ of the $\Line(a,b)$
  and $\Line(c,d)$ makes alternate interior angles unequal.  This condition can 
  be expressed using points only, as shown in 
  Fig.~\ref{figure:AlternateInteriorAnglesFigure} above. We can use the 
  strong parallel axiom to prove stability:

  $$ \neg \neg \il(a,b,c,d) \defined \implies \il(a,b,c,d) \defined.$$
 
But one cannot do this for subtheories with no parallel postulate or other 
versions of the parallel postulate.  Therefore we prefer, when working with 
$\il(a,b,c,d)$ to use the Logic of Partial Terms (described below),  in which 
$t \defined$ is made into an official atomic formula for each term $t$, instead of 
an abbreviation at the meta-level.

\subsection{Logic of Partial Terms (\LPT)} \label{section:lpt}
This is a modification of first-order logic, in which the formation rules for 
formulas are extended by adding the following rule:  If $t$ is a term then $t \defined$ is a formula.  Then 
in addition the quantifier rules are modified so instead of $\forall x A(x) \implies A(t)$ we have 
$\forall x (t \defined \land A(x)) \implies A(t)$,  and instead of $A(t) \implies \exists x\, A(x)$ we have
$A(t) \land t\defined \implies \exists x\, A(x)$.  Details of  \LPT\ can be found in \cite{beeson-book}, p.~97.
 
\LPT\ includes axioms $c \defined$ for all constants $c$ of any theory formulated in \LPT; 
this is in accordance with the philosophy that terms denote things, and while terms may fail to denote (as in ``the King of France''),
there is no such thing as a non-existent thing.   Thus $1/0$ can be undefined, i.e.,fail to denote,  but if a constant $\infty$ is 
used in \LPT, it must denote something.  

The meaning of $t=s$ is that $t$ and $s$ are both defined and they are equal.   We write $t \cong s$ to 
express that if one of $t$ or $s$ is defined, then so is the other, and they are equal.  
\begin{Definition} \label{defn:equalifdefined}  For terms in any theory using the logic of partial terms, 
$t \cong q$ means
$$(t \defined \implies t=q )\land (q \defined \implies t=q).$$
This is read {\bf $t$ and $q$ are equal if defined}.
\end{Definition}
\medskip
Thus ``$\cong$'' 
is an abbreviation at the meta-level, rather than a symbol of the language. 

  \LPT\ contains the axioms
of ``strictness'',  which are as follows (for each function symbol $f$ and relation symbol $R$ in the language):
\begin{eqnarray*}
&&f(t_1,\ldots,t_n)\defined \implies t_1\defined \land \ldots \land t_n\defined \\
&&R(t_1,\ldots,t_n) \implies t_1\defined \land \ldots \land t_n\defined 
\end{eqnarray*}

Note that in \LPT, under a given ``valuation'' (assignment of elements of a structure to variables), each formula 
has a definite truth value, i.e.,we do not use three-valued logic in the semantics.  For example, if $P$ is a formula
of field theory with a reciprocal operation $1/x$, then $P(1/0)$ is false, since $1/0$ is undefined.   For the 
same reason $\neg P(1/0)$ is false.   Hence $P(1/0) \lor \neg P(1/0)$ is false too;  but that does not contradict the
classical validity of $\forall x(P(x) \lor \neg P(x))$ since we are required to prove $t\defined$ before deducing
an instance $P(t) \lor \neg P(t)$.  

As an example of the use of \LPT, we reformulate the theory of Euclidean 
fields \cite{beeson-bsl} using the logic of partial terms. 
 The existential quantifiers associated with the reciprocal axioms, with the axiom of additive inverse, 
 and with the  square-root axiom of 
Euclidean field theory are replaced by a function symbol $\sqrt{}$, a unary minus $-$,  and a function symbol for ``reciprocal'', which 
we write as $1/x$ instead of $\mbox{\it reciprocal}(x)$.    The  changed axioms are 
\medskip

\axioms
$x + (-x) = 0$ &(additive inverse) \\
$x \neq 0 \implies x \cdot (1/x) = 1 $  &(EF1${}^\prime$) \\
$P(x) \implies x \cdot (1/x) = 1 $  &(EF7${}^\prime$) \\
$x + y = 0 \land \neg P(y) \implies \sqrt x \cdot \sqrt x  = x $& (EF5${}^\prime$)
\endaxioms 
\medskip

\subsection{A version of Tarski's theory with terms for intersections of lines}
This version of Tarski's theory we call ruler and compass Tarski.   It is formulated as follows:
\begin{itemize}
\item It uses a function symbol $\il(a,b,p,q)$ for the intersection point of $\Line(a,b)$ and 
$\Line(p,q)$.  
\item It uses the logic of partial terms.
\item If $\il(a,b,p,q)$ is defined, then it is a point on both lines.
\item If there is a point on $\Line(a,b)$ and $\Line(p,q)$, and those lines do not coincide,
then $\il(a,b,p,q)$ is such a point.
\item Formally, the axioms involving $\il$ are
$$
Col(a,b,x) \land Col(p,q,x) \land \neg(Col(a,b,p) \land Col(a,b,q)) \implies \il(a,b,p,q) \defined 
$$
$$ \il(a,b,p,q) \defined \implies a \neq b \land p \neq q \land Col(a,b,\il(a,b,p,q)) \land Col(p,q,\il(a,b,p,q))
$$
\item $\il$ is used instead of separate Skolem functions for $ip$.  Specifically, the 
term $ip(a,p,c,b,q)$ in the Skolemized inner Pasch axiom become $\il(a,q,b,p)$.  The point $c$
does not occur in this term.
\item The Skolem functions $ext$ (for segment extension) is not changed.
\item The Skolem functions for intersections of lines and circles are not changed.
\item Stability for equality, betweenness, and congruence, as before.
\end{itemize}
We could consider replacing Skolem terms   $center(a,b,c)$ with terms built up from $\il$.
The two lines to be intersected are 
the perpendicular bisectors of $ab$ and $bc$, where $a$, $b$, and $c$ are three non-collinear
points.   We can define the perpendicular bisector by the erected
perpendicular at the midpoint, so it is indeed possible to eliminate the symbol $center$;
but there seems to be no special reason to do so.

We did not include the stability of definedness; that is because it can be proved.
The following lemma is proved in \cite{beeson-bsl}; here we give a different proof,
based on the triangle-circumscription form of the strong parallel axiom.

\begin{Lemma}\label{lemma:il-stability}
[Stability of $\il(a,b,c,d)$]
   The strong parallel postulate is 
equivalent (in ruler and compass Tarski minus the parallel postulate) to the 
stability of $\il(a,b,c,d) \defined$:
$$ \neg \neg il(a,b,c,d)\defined \implies \il(a,b,c,d) \defined.$$
\end{Lemma}

\noindent{\em Proof}.  (i)  First suppose the strong parallel postulate and $\neg \neg \il(a,b,c,d)\defined$.  We will show $\il(a,b,c,d)\defined$.   
Let $L = \Line(a,b)$ and $K = \Line(c,d)$.  Then lines $K$ and $L$ do not coincide, for 
then $\il(a,b,c,d)$ would be undefined.   Hence, by the strong parallel postulate,
 we can find a point on $L$ that is not on $K$.
We may assume without loss of generality that $b$ is such a point.  Construct point $f$ so that $bf$ is parallel 
to $K$; more explicitly, $K$ and $bf$ and the transversal $bc$ make alternate interior
angles equal.  If $a$, $b$, and $f$ are collinear, then $ab$ and $cd$ are parallel, so 
$\il(a,b,c,d)$ is undefined, contradiction.  Hence $a$, $b$, and $f$ are not collinear.
Then line $M = \Line(b,f)$ passes through point $b$ and is parallel to $K$, and line $L$
also passes through $b$, and has a point $a$ not on $M$.  Then by the strong parallel axiom,
$L$ meets $K$.  In that case $\il(a,b,c,d)$ is defined, as claimed.

(ii)  Conversely, suppose the stability of $\il(a,b,c,d)\defined$, and suppose $a$, $b$, and $c$ are 
not collinear.  Let $m$ be the midpoint of $ab$ and $n$ the midpoint of $cd$, 
with $pm$ the perpendicular bisector of $ab$ and $qn$ the perpendicular bisector of $cd$.
We must prove $\il(m,p,n,q)\defined$.  By stability it suffices to derive a contradiction 
from the assumption that it is not defined. If it is not defined then $mp$ is parallel to $nq$
(as not meeting is the definition of parallel).  But $\Line(a,b)$ and $\Line(b,c)$ are perpendicular to $mp$ 
and $nq$ respectively; hence they cannot fail to be parallel or coincident.  But since they
both contain point $b$, they are not parallel;  hence they are coincident.  Hence $a$, $b$, and 
$c$ are collinear, contradiction.   That completes the proof.

\begin{Theorem}[Stability of definedness] For each term $t$ of ruler and compass Tarski,
$\neg \neg t \defined \implies t \defined$ is provable.
\end{Theorem}

\noindent{\em Proof}.  By induction on the complexity of the term $t$.  If $t$ is a 
compound term $ts$, and $\neg \neg ts \defined$, then $\neg \neg t \defined$ and 
$\neg \neg s\defined$, so by induction hypothesis, $t \defined$ and $s \defined$.
Hence $ts \defined$.   We may therefore suppose $t$ is not a compound term.  If the 
functor is $\icone$, $\ictwo$, $\ilcone$, or $\ilctwo$, then it is easy to prove that 
the conditions for $t$ to be defined are given geometrically, by the same formulas
that were used to define $t \defined$ in Tarski with Skolem functions (and without \LPT).
Hence stability follows by the stability of equality, congruence, and betweenness.  The
stability of $\il(a,b,c,d)$ is equivalent to the strong parallel postulate, by 
the previous lemma.  That completes the proof.
\medskip

\subsection{Intersections of lines and the parallel axiom}
In the proof of the first part of
 Lemma~\ref{lemma:il-stability}, we showed that if lines $L$ and $M$ meet in a 
point $x$, then $x$ can be made to appear as the center of a circle circumscribed about 
suitably chosen points $a$, $b$, and $c$.  In this section, we will refine this construction to show 
that there is a single term $t(a,b,c,d)$ in the language of Tarski with Skolem functions
that gives the intersection point of $\Line(a,b)$ and $\Line(c,d)$, when it exists.

\begin{Lemma} \label{lemma:pointnotonL}
  Given two lines $L$ and $K$ that are neither coincident nor parallel,
one can construct a point $p$ that lies on $K$ but not on $L$.  More precisely, interpreting
$L$ as $\Line(a,b)$ and $K$ as $\Line(q,r)$,  there
is a single term $t(a,b,q,r)$ such that if $\neg(Col(a,b,q) \land Col(a,b,r))$ and 
for some $x$, $Col(a,b,x) \land Col(q,r,x)$,  then $e=t(a,b,q,r)$ satisfies
$Col(q,r,e) \land \neg Col(a,b,e)$.   
\end{Lemma}

\noindent{\em Remarks}.  The point $x$ cannot be used to construct $e$, which must 
depend only on $a$, $b$, $q$, and $r$, and must  be constructed by a single term, 
and hence depend continuously on the four parameters.  We will use the parallel 
postulate to construct $e$;  we do not know a construction that does not use the 
parallel postulate.
\medskip

\noindent{\em Proof}.   Let $M$ be the perpendicular to $K$ passing through $q$.
We are supposed to construct $M$ from $a$, $b$, $q$, and $r$ alone.  To construct
$M$, we need not just $p$ and $q$, but also a point not on $K$; and $a$ and $b$ 
are useless here as they might lie on $K$.   We must appeal to Lemma~\ref{lemma:pointoffline}
for the construction of some point not on line $K$; thus this apparently innocent
lemma requires the geometric definition of arithmetic and the introduction of coordinates,
and hence the parallel postulate.

The construction of $M$ gives us one point $u$ on $M$ different from $q$.  Let
$v$ be the reflection of $u$ in $q$.  Then $u$ and $v$ are equidistant from $q$.
Now, using the uniform perpendicular construction we construct the line $J$ through $u$ 
perpendicular to $L$. See Fig.~\ref{figure:pointnotonL}.%
\footnote{
In a theory with Skolem functions for the intersection points
of two circles,  the construction of $M$ and $u$ becomes trivial (just use the 
Euclidean construction of a perpendicular), but axiomatizing the Skolem functions 
for the circle in such a way as to distinguish the points of intersection (which is 
necessary to construct perpendiculars) requires the introduction of coordinates.  
So using circle-circle does not obviate the need for coordinates in this lemma.
}

\begin{figure}
\caption{Uniform construction of a point $c = center(d,u,v)$ on $K$ that is not on $L$.  The construction works whether or not $q$ is on $L$, or $u$ is on $L$. The dotted line bisects 
$ud$ and does not coincide with $L$.}
\label{figure:pointnotonL}
\PointnotonLFigure
\end{figure}

 While we do not know whether $u$ lies on $L$, the uniform 
perpendicular construction (Theorem~\ref{theorem:uniformperp})  provides two points
determining $J$, namely 
$$f = Project(u,a,b,c)$$ and   $$h= head(u,a,b,c),$$ where $f$ is on $L$
and $c$ is not on $L$.   Possibly $u$ is equal to $f$ or to $h$; we need a point $d$ on $J$
that is definitely not equal to $u$ or to the reflection of $u$ in $L$.
  To get one, our plan is to draw a circle of sufficiently large radius
about $u$ and intersect it with $J=\Line(f,h)$.  We use the uniform reflection
construction to define $w=Reflect(u,a,b)$, the reflection of $w$ in $L$.
Then we extend the non-null segment $\alpha\beta$ by the (possibly null) segment 
$uw$ to get a point $z$ such that $\alpha z > uw$.  Then we use $\alpha z$ as the 
``sufficiently large'' radius.   Here is the construction:
\begin{eqnarray*}
z &=& ext(\alpha,\beta,u,w) \\
d &=& \ilcone(f,c,u,e_2(u,\alpha,z))
\end{eqnarray*}
Now $d$ lies on $J$ and is different from $u$, and it is also 
different from $w$ since $w$ lies inside the circle centered at $u$ of 
radius $\alpha z$.   Finally define
$$ c = center(u,v,d).$$
The three points $u$, $v$, and $d$ are not collinear, since then $J$ and $M$ would 
coincide, and $L$ and $K$ would both be perpendicular to $J$, and hence parallel; but 
$L$ and $K$ are by hypothesis not parallel.
Since $u$, $v$, and $d$ are not collinear,
 $c$ is equidistant from $u$, $v$, and $d$.  Therefore $c$ lies on the perpendicular 
bisector of $uv$, which is $K$.  Also $c$ lies on the perpendicular bisector of 
$ud$, which is parallel to $L$, since both are perpendicular to $J$. This 
perpendicular bisector does not coincide with  $L$, since
 $d$ is not the reflection $w$ of $u$ in $L$.   Therefore $c$ does not lie on $L$.
 Then $c$ lies on $K$ but not on $L$, as desired.  That completes the proof of the lemma.

\begin{Theorem}\label{theorem:il-elimination}[Elimination of $\il$]
There is a term $t(a,b,p,r)$ of intuitionistic Tarski with Skolem functions (so $t$ contains
$ip$ and $center$ but not $\il$)  such that the following is provable:
$$ Col(a,b,x)  \land Col(p,r,x) \land \neg(Col(a,b,p) \land Col(a,b,r)) \land p \neq r
 \implies x = t(a,b,p,r).$$
Less formally, $t(a,b,p,r)$ gives the intersection point of $\Line(a,b)$ and $\Line(p,r)$.
\end{Theorem}

\noindent
{\em Remark}.  The problem here is to explicitly produce the term $t$ that is 
implicit in the proof of Lemma~\ref{lemma:pointnotonL}.  This is also closely
related to  the proof given in 
 \cite{beeson-bsl}  that in constructive neutral 
 geometry, the triangle circumscription principle implies the strong parallel axiom.
But here we have to check that this construction can be carried out in Tarski's 
geometry, i.e., all the lines that are required to intersect are proved to intersect
using only inner Pasch; so there is definitely something additional to check.  
\medskip

\noindent{\em Proof}.  By Lemma~\ref{lemma:pointnotonL}, we may assume without 
loss of generality that $p$ does 
not lie on $L$.  More explicitly, if we now produce a term $t$ as in the lemma, that 
works under the additional assumption that $p$ is not on $L$, then we can compose that 
term with the term given in Lemma~\ref{lemma:pointnotonL}, which produces a point 
on $\Line(p,r)$ that is not on $L$, and the composed term will 
work without the assumption that $p$ is not on $L$.  
 
Recall that $Project(p,a,b)$ is the point 
$w$ on $\Line(a,b)$ such that $pw \perp ab$, and $Project(p,a,b)$ is given by a term 
in Tarski with Skolem functions.   There is also a term $erect(p,r)$ that produces
a point $j$ such that $j$ is not on $\Line(p,r)$ and $jp \perp pr$. (Since $p$ is 
on $\Line(p,r)$, the uniform perpendicular construction is not needed; the simple 
Euclidean construction will be enough.)   Finally, there is a term $\Reflect(x,a,b)$
that produces the reflection of $x$ in $\Line(a,b)$; if we assume $x$ is not 
on $\Line(a,b)$ this is particularly easy:  let $q = Project(x,a,b)$ and take
$\Reflect(x,a,b) = ext(x,q,x,q)$. The requirement that $x$ is not on $\Line(a,b)$
means that we are extending a non-trivial segment, which is constructively allowed.

\begin{figure}[h]
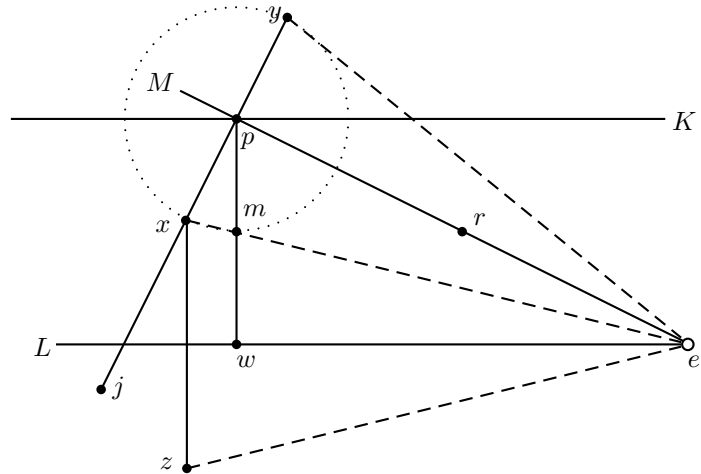
   
\caption{Triangle circumscription implies the  strong parallel axiom. Given lines $L$
and $M$, to construct their intersection point as the center $e$ of an appropriate circle.
$y$ and $z$ are reflections of $x$ in $M$ and $L$.
}
\label{figure:CenterImpliesStrongParallelFigure}  
\hskip 1.9cm   
\CenterImpliesStrongParallelFigure
\end{figure}
\medskip

In Fig.~\ref{figure:CenterImpliesStrongParallelFigure}, $L = \Line(a,b)$ and $M=\Line(p,r)$,
and $p$ does not lie on $L$.   First we claim that $x$ does not lie 
on $L$.  If $x$ does lie on $L$, then $pwx$ is a right triangle, so the 
hypotenuse $px$ is greater than the leg $pw$; but $px=mp$, which is less 
than $pw$ since $m = midpoint(p,w)$, and $p \neq w$.  Therefore $\Reflect(x,a,b)$
can be defined using the easy construction given above.

Now we give a construction script corresponding to the figure:
\begin{eqnarray*}
w &=& Project(p,a,b)\\
m &=& midpoint(p,w)\\
j &=& erect(p,r)\\
x &=& \ilcone(j,p,p,m)\\
y &=& \ilctwo(j,p,p,m)\\
z &=& \Reflect(x,a,b)\\
e &=& center(x,y,z)
\end{eqnarray*}
Composing the terms listed above we find a (rather long) term that produces
$e$ from $a$, $b$, $p$, and $r$.  We claim that $e$ is the intersection 
point of $\Line(a,b)$ and $\Line(p,r)$.
By the stability of collinearity, we can argue by cases on whether $x$, $y$, 
and $z$ are collinear or not.  (Here it is important that we are not proving
a statement with an existential quantifier, but a quantifier-free statement 
involving a single term that constructs the desired point $e$.)
In case $x$, $y$, and $z$ are not collinear, then $center(x,y,z)$ produces
a point $e$ that is the center of a circle containing $x$, $y$, and $z$.
Then Euclid III.1  implies that $L$ and $M$ both pass through $e$, and we 
are done.  On the other hand, if
 $x$, $y$, and $z$ are collinear, then $M$ and $L$ are both perpendicular to 
$\Line(x,y)$, so $M$ and $L$ are parallel; then there is no $x$ as in the hypothesis 
of the formula that is alleged, so there is nothing more to prove.
That completes the proof of the theorem.

\subsection{Interpreting ruler and compass Tarski in intuitionistic   Tarski}
\label{section:meta}

Ruler-and-compass Tarski clearly suffices to interpret intuitionistic Tarski (with or without 
Skolem functions), because the points asserted to exist by inner Pasch and the triangle 
circumscription principle are given as intersections of lines. 
  Conversely we may ask, whether ruler and compass Tarski can be 
interpreted in intuitionistic Tarski with Skolem functions.  That is, can terms in $\il$ 
be effectively replaced by terms built up from $ip$ and $center$?  The answer is ``yes''.

\begin{Theorem} \label{theorem:il-elimination2}
 Suppose ruler and compass Tarski geometry (with $\il$ and other 
Skolem functions) proves a theorem $\phi$ that does not contain $\il$.  Then 
Tarski geometry with Skolem functions proves $\phi$.   Moreover if $\phi$ contains
no Skolem functions, then Tarski geometry proves $\phi$.  These claims hold both 
for the theories
with intuitionistic logic and those with classical logic.
\end{Theorem}

\noindent{\em Proof}.  We assign to each term $t$ of ruler and compass Tarski, 
a corresponding term $t^\circ$ of   Tarski with Skolem functions.
Let $\il^\circ(a,b,q,r)$ be the term given by Theorem~\ref{theorem:il-elimination}.
The term $t^\circ$ is defined inductively by 
\begin{eqnarray*}
x^\circ &=& x \mbox{\qquad where $x$ is a variable or constant} \\
\il(a,b,c,d)^\circ &=& \il^\circ(a^\circ,b^\circ,c^\circ,d^\circ) \\
f(a,b,c,d) &=& f(a^\circ,b^\circ,c^\circ,d^\circ) \mbox{\qquad where $f$ is $\ilcone$, $\ilctwo$, $\icone$, $\ictwo$, or $ext$}\\
e(a)^\circ &=& e(a^\circ)
\end{eqnarray*}
Then we assign to each formula $A$ of ruler and compass Tarski, a corresponding
formula $A^\circ$ of   Tarski with Skolem functions.  Namely, 
the map $A \mapsto A^\circ$ commutes with logical operations and quantifiers, and
for atomic $A$ not of the form $t\defined$, we define
$$A(t_1,\ldots,t_n)^\circ = A(t_1^\circ,\ldots,t_n^\circ).$$
For the case when $A$ is $t\defined$, we define $(t \defined)^\circ$ to be $t=t$ 
when $t$ is a variable or constant, and when it is a compound term, we use
Definition~\ref{definition:definedness}.
By induction on the complexity of $A$, 
$A^\circ[x:=t^\circ]$ is provably equivalent to 
$(A[x:=t])^\circ$.

Then by induction on the length of proofs in ruler and compass Tarski, we show that
if ruler and compass Tarski proves $\phi$, then   Tarski with Skolem 
functions proves $\phi^\circ$.  A propositional axiom or inference remains one under
the interpretation, so it is not even vital to specify exactly which propositional axioms
we are using.  In this direction (from  \LPT\ to ordinary logic), the quantifier rules 
and axioms need no verification, as the extra conditions of definedness needed in
\LPT\ are superfluous in ordinary logic.
For example, one of those axioms is $\forall x A \land t \defined \implies A[x:=t]$.
That becomes $\forall x A^\circ \land (t\defined)^\circ \implies A[x:=t^\circ]$,
in which the $t\defined $ can just be dropped.  There are some special axioms
in \LPT, for example $c \defined$ for $c$ a constant and $x \defined$ for $x$ a 
variable.  

We check the  basic axioms for $\ilcone$.  These say that (i) if 
$\ilcone(a,b,c,d)$ is defined, then it is a point on $\Line(a,b)$ and also on 
the circle with center $c$ passing through $d$,   and  (ii) if there is an point $x$
on both the line and circle, then $\ilcone(a,b,c,d)$ is defined.
According to the definition of $(\ilcone(a,b,c,d)\defined)^\circ$, the interpretation of 
``$\ilcone \defined$'' is ``there is a point on the line inside the circle'', where 
``inside'' means not strictly inside.  Since $\ilcone$ is not affected by the interpretation
(except in the atomic formula $\ilcone(a,b,c,d)\defined$), the interpretations of 
the basic axioms for $\ilcone$ are equivalent to those same axioms.  Similarly for 
$\ilctwo$, $\icone$, and $\ictwo$.  

Now consider the axioms for $\il(a,b,p,r)$.  These   say that if $\il(a,b,p,r)\defined$
then $e=\il(a,b,p,r)$ is a point on $\Line(a, b)$ and also on $\Line(p,r)$, and if 
$x$ is any point on both lines and not both $p$ and $r$ lie on $\Line(a,b)$, then
then $x=\il(a,b,p,r)$.   Recall that $(\il(a,b,p,r)\defined)^\circ$  says there exists an $x$
on both lines, and not both $p$ and $r$ are on $\Line(a,b)$.  Then the interpretation
of these axioms says that if $\Line(a,b)$ and $\Line(p,r)$ meet and not both $p$ and $r$
lie on $\Line(a,b)$, then $\il(a,b,p,r)^0$ is the intersection point.  But that 
is exactly Theorem~\ref{theorem:il-elimination}.  That completes the proof.
\medskip

\begin{Corollary} \label{theorem:classical-il-elimination}
Suppose (classical or intuitionistic) ruler and compass Tarski geometry proves a theorem
of the form $P \implies   Q(t)$, with $P$ and $Q$ negative.  Then (classical or 
intuitionistic)  Tarski geometry proves $P \implies \exists x\, Q(x)$.
\end{Corollary}

\noindent{\em Proof}.  Suppose  
 $P \implies   Q(t)$ is  provable in   ruler and compass
 geometry.  Then $P \implies \exists x\, Q(x)$ is provable, with
 one more inference.  
  But that formula contains no occurrences of $\il$.  
 Then by Theorem~\ref{theorem:il-elimination2}, it is provable in classical
 Tarski geometry.    That completes the proof.

\noindent{\em Remark}.  If we drop the strong parallel axiom (or triangle circumscription principle),
we obtain ``neutral 
geometry with $\il$\,'').   It is an open question whether neutral geometry with $\il$ can 
be interpreted in neutral Tarski with Skolem functions.  In other words,  can all terms for 
intersection points of lines that are needed in proofs of theorems not mentioning $\il$
be replaced by terms built up from $ip$ and $ext$?  We used $center$ in an essential way
in the proof of Theorem~\ref{theorem:il-elimination}, but did we have to do so?  

\section{Relations between classical and constructive geometry}

Our intuition about constructive geometry is this:  You may argue classically for 
the equality or inequality of points, for the betweenness of points, for collinearity,
for the congruence of segments.  But if you assert that something exists, it must be
constructed by a single, uniform construction, not by different constructions applying 
in different cases.   If you can give a uniform construction, you may argue by cases
that it works, but the construction itself cannot make a case distinction.   Thus the 
uniform perpendicular construction of a line through $x$ perpendicular to $L$ works 
whether or not $x$ is on $L$; if we wished, we could argue for its {\em correctness}
by cases, as we could always push a double negation through the entire argument and 
use stability to eliminate it.

There is in Szmielew's Part I of \cite{schwabhauser} an extensive development from 
Tarski's classical axioms, essentially deriving Hilbert's axioms and the definitions and
key properties of addition and multiplication.  We would like to be able to import 
arguments and results wholesale from this development into constructive geometry.  
In this section we investigate to what extent this is possible.   

It is certainly not completely possible to import results without modifying Tarski's axioms, 
since constructive proofs
will produce points that depend continuously   on parameters, while as we have 
discussed above, Tarski's version of inner Pasch and segment extension axioms do not have
this property.  Those defects have been remedied above by formulating ``continuous Tarski
geometry'',  a theory classically equivalent to Tarski's geometry.  

\subsection{The double-negation interpretation} \label{section:doublenegation}
The difference between intuitionistic and classical mathematics shows up
in the interpretations of existence ($\exists$) and disjunction ($\lor$).
 
\begin{Definition} A formula is called {\bf negative} if it does not contain 
$\exists$ or $\lor$ (existential quantifiers or disjunction).
\end{Definition}
Intuitively, negative formulas make no existential claims, and hence have the 
same meaning classically as constructively, at least in theories with stable atomic
formulas.  
Classically we can always express ``exists'' as ``not for all not'' and express $P \lor Q$
as $\neg(\neg P \land \neg Q)$.  If we do this, then every classical theorem should become
intuitionistically meaningful:  classical mathematics is, somewhat surprisingly,
contained in intuitionistic mathematics,  although often  intuitionistic 
mathematics is considered as a restriction of classical mathematics.

G\"odel made these considerations precise by introducing his double-negation 
interpretation  \cite{godel} , which  
assigns a negative formula $A^-$ to every formula $A$,  by replacing $\exists $ by $\neg \forall \neg$  
   and replacing $A \lor B$ by $\neg (\neg A \land \neg B)$.   For atomic formulae, $A^-$ is defined to be $\neg \neg A$.
The rules of intuitionistic logic are such that if $A$ is classically provable (in predicate logic) then 
$A^-$ is intuitionistically provable.   Hence, if we have a theory $T$ with classical logic, and 
another theory $S$ with intuitionistic logic, whose language includes that of $T$,  and for every
axiom $A$ of $T$, $S$ proves $A^-$,  then $S$ also proves $A^-$ for every theorem $A$ of $T$.
In case the atomic formulas in the language of $T$ are stable in $S$, i.e.,equivalent to their double negations,
then of course we can drop the double negations on atomic formulas in $A^-$.

In \cite{beeson-kobe} we applied this theorem to a version of constructive geometry
 based on Hilbert's axioms.  Given
the extensive almost-formal development of geometry from Tarski's axioms in \cite{schwabhauser}, one
might like to use the double-negation interpretation 
with  $T$ taken to be Tarski's theory, and 
$S$ taken to be  some suitable constructive version of Tarski's theory.  We now investigate
this possibility.

A double-negation interpretation from a classical theory to a constructive version of that
theory becomes a better theorem if it applies to the Skolemized versions of the theories, because
in the un-Skolemized version,  an existential quantifier is double-negated, while the corresponding
formula of the Skolemized theory may replace the existentially quantified variable by a term,
so no double negated quantifier is involved, and no constructive content is lost.
But if we Skolemize Tarski's version of inner Pasch, we get an essentially non-constructive
axiom, as shown above. Hence there is no double-negation interpretation for that theory.
However, it works fine if we replace Tarski's axioms by the (classically equivalent) 
axioms of continuous Tarski with Skolem functions:

\begin{Theorem}\label{theorem:doublenegation}
Let $T$ be intuitionistic Tarski geometry with Skolem functions.  If $T$ plus classical
logic proves $\phi$, then $T$ proves the double-negation interpretation $\phi^-$.   
\end{Theorem}

\noindent{\em Proof}.  It suffices to verify that the double-negation interpretations 
of the axioms are provable.  But the axioms are negative and quantifier-free, so they 
are their own double-negation interpretations.  
That completes the proof.
\medskip

\begin{Corollary} \label{corollary:doublenegation} If $\phi$ is negative, and classical Tarski geometry without Skolem
functions proves $\phi$,
then intuitionistic Tarski geometry proves $\phi$.
\end{Corollary}

\noindent{\em Proof}.  Suppose
 $\phi$ is provable in classical Tarski geometry (with or without Skolem functions).
Then since $\phi$ itself has no Skolem functions, 
$\phi$ is provable in classical Tarski geometry without Skolem functions, and hence 
by Theorem~\ref{theorem:intuitionisticTarskiImpliesClassical}, it is provable
in intuitionistic Tarski geometry with Skolem functions. Hence, 
by Theorem~\ref{theorem:doublenegation},
$\phi^-$ is provable.  Since $\phi$ is negative, it is equivalent to $\phi^-$.
That completes the proof. 

\subsection{Applications of the double-negation interpretation}

We illustrate the use of Theorem~\ref{theorem:doublenegation} by importing the work of 
Eva Kallin, Scott Taylor, H.~N.~Gupta, and Tarski mentioned in Section~\ref{section:TarskiProvesHilbert}. 

\begin{Corollary} The formulas (A16) through (A18),  which were once axioms of Tarski's
theory, but were shown classically provable from the remaining axioms,  are also provable
in intuitionistic Tarski  without Skolem functions.
\end{Corollary}

\noindent{\em Proof}.  By Corollary~\ref{corollary:doublenegation}. 
\medskip

We would like to emphasize
something has been achieved with the double-negation theorem even for negative theorems,
as it would be quite laborious to check the long proofs of (A16)-(A18) directly to verify 
that they are constructive.  For example, (A18) is Satz 5.1 in \cite{schwabhauser}.  Let us 
consider trying to check directly if this proof is constructive.   You can 
see that the proof proceeds by contradiction, which is permissible by stability;
in the crucial part of the proof, inner Pasch is applied to a triangle which 
ultimately must collapse (as the contradiction is reached) to a single point.  Therefore we 
can constructivize this part provided the non-collinearity hypothesis is satisfied for the 
application of Pasch.  By stability, we may assume that the vertices of this triangle are 
actually collinear.  But can we finish the proof in that case?  It looks plausible that
(A15) or similar propositions might apply, but it is far from clear. Yet
the double-negation interpretation applies, and we do not need to settle the issue by hand.
We had to assume (A15), but we do not have to assume (A18), because it is already provable.

\begin{Corollary} \label{corollary:gupta}
The correctness of Gupta's perpendicular construction given in \S\ref{figure:GuptaPerpendicularFigure} is 
provable in intuitionistic Tarski geometry.
\end{Corollary}

\noindent{\em Proof}.  Please refer to Fig.~\ref{figure:GuptaPerpendicularFigure} and the discussion 
in \S~\ref{section:gupta}.  Two conclusions are to be proved:  (i) that $x$ lies on $Line(a,b)$,
and (ii) that $cx \perp ax$.  We will show that both of these are equivalent to 
negative formulas; hence the double-negation interpretation applies.

Ad (i):  The statement that $x$ lies on $Line(a,b)$ is 
informal; its exact formal meaning is 
$$ \neg \neg\, (\T(x,a,b) \lor \T(a,x,b) \lor \T(a,b,x)). $$
(See the discussion after Lemma~\ref{lemma:stability2981}).
Pushing inner negation sign through the disjunction we see that this is equivalent 
to a negative sentence.

Ad (ii): By the definition of perpendicularity,
and the theorem that $L \perp K$ if and only if $K \perp L$,  and the fact that $x$ is by 
definition the midpoint of $cc^\prime$, $cx \perp ax$ is equivalent to $ac = ac^\prime$,
which is negative.  Hence the double-negation interpretation applies to that conclusion.
 
It follows from the double-negation interpretation that the correctness of Gupta's construction 
is provable in intuitionistic Tarski geometry.  That completes the proof.
\medskip

The following theorems (numbered as in \cite{schwabhauser}) have proofs
simple enough to check directly (as we did before developing the double negation interpretation),
but with the aid of the double negation interpretation, we do not need to check them directly.

\begin{Lemma} \label{lemma:Ch3}
 The following basic properties of betweenness are provable
in intuitionistic Tarski geometry.  Note that $\T(a,b,c)$ is a defined 
concept; $\B(a,b,c)$ is primitive.  The theorem numbers refer to \cite{schwabhauser}.
\smallskip

\axioms 
$\T(a,b,b)$ &  Satz 3.1\\
$\T(a,b,c) \implies \T(c,b,a)$ & Satz 3.2 \\
$\T(a,a,b)$ & Satz 3.3 \\
$\T(a,b,c) \land \T(b,a,c) \implies a=b$ & Satz 3.4\\
$\T(a,b,d) \land \T(b,c,d) \implies \T(a,b,c)$ & Satz 3.5a\\
$\T(a,b,c) \land \T(a,c,d) \implies \T(b,c,d)$ & Satz 3.6a \\
$\T(a,b,c) \land \T(b,c,d) \land b \neq c \implies \T(a,c,d)$& Satz 3.7a \\
$\T(a,b,d) \land \T(b,c,d) \implies \T(a,c,d)$ & Satz 3.5b \\
$\T(a,b,c) \land \T(a,c,d) \implies \T(a,b,d)$ & Satz 3.6b \\
$\T(a,b,c) \land \T(b,c,d) \land b \neq c \implies \T(a,b,d)$ & Satz 3.7b
\endaxioms
\end{Lemma}
    
\noindent{\em Proof}. By the double negation interpretation, since each of these 
theorems is negative.
\medskip

Does the double negation interpretation help us to be able to ``import'' proofs of existential theorems 
from \cite{schwabhauser}
to intuitionistic Tarski?  
It gives us the following recipe:  Given an existential theorem 
proved in classical Tarski,  we  examine the proof to see if we can 
construct a Skolem term (or terms)
for the point(s) asserted to exist.   If the proof constructs points using inner Pasch,
we need to verify whether degenerate cases or a possibly collinear case are used.   If they are not used
 then the strict inner Pasch axiom (A7-i) suffices.   
The crucial question is whether the point alleged to exist can be 
constructed by a single term,  or whether the proof is an argument by cases in which
different terms are used for different cases.  In the latter case, the proof is not
constructive (though the theorem might still be, with a different proof).  But in the 
former case, the double-negation interpretation will apply.

Thus the double-negation interpretation fully justifies the claim that the essence of constructive
geometry is the avoidance of arguments by cases,  providing instead uniform constructions
depending continuously on parameters.

\subsection{Euclid I.2  revisited}

Consider the first axiom of Tarski's geometry, which says 
any segment (null or not) can be extended:  $\exists d\,(\T(a,b,d) \land bd = bc)$.
Clearly $d$ cannot depend continuously on $a$ as $a$ approaches $b$ while $b$ and $c$ 
remain fixed, since as $a$ spirals in towards $b$,  $d$ circles around $b$ outside
a fixed circle.   Therefore Axiom (A4) of Tarski's (classical) geometry is essentially
non-constructive; the modification to (A4-i) that we made in order to pass to a 
constructive version was essential.

Euclid I.2 says that given three points $a,b,c$, there exists a point $d$ 
such that $ad = bc$.  Euclid gave a clever proof that works when the three points are 
distinct, and classically a simple argument by cases completes the proof.  Constructively,
that does not work, since
 when $b$ and $c$ remain fixed and $a$ approaches $b$,  $d$ from Euclid's
construction does not depend continuously on $a$.  We will show  in this section that it is only 
Euclid's proof that is non-constructive;  the theorem itself is provable in 
intuitionistic Tarski geometry, by a different proof.

\begin{Lemma} \label{lemma:subsegment}
Intuitionistic Tarski geometry proves 
\smallskip

(i) $\T(a,b,c)$ and $\T(p,q,r)$ and $ac = pr$ and $bc = qr$
implies $ab = pq$

(ii) a segment $ac$
cannot be congruent to (a proper subsegment) $bc$  with $\B(a,b,c)$.
\end{Lemma}

\noindent{\em Proof}.   We first show (ii) follows from (i).  Suppose $ac = bc$
and $\B(a,b,c)$.  Then  in (i) take $p=q=a$ and $r=c$.  Then (i) implies
$ab=aa$ , contrary to axiom (A3).  Hence (i) implies (ii) as claimed.

Now (i) is Satz 4.3 in \cite{schwabhauser}, and since it is negative, 
we can conclude from the double negation interpretation that it is 
constructively provable.  That completes the proof.

\begin{Lemma} \label{lemma:EuclidI.2}
In intuitionistic Tarski geometry, null segments can be extended, and Euclid I.2
is provable.  Indeed, there is a term (using Skolem functions) $e(x)$ such that 
$e(x) \neq x$ is provable, and a term $e_2$ corresponding to Euclid I.2, such that 
if $d=e_2(a,b,c)$, then $ad=bc$.
\end{Lemma}

\noindent{\em Remarks}.
Thus, it is only Euclid's proof of I.2 that is non-constructive, as discussed in \cite{beeson-bsl},
not the theorem itself.  Note that a constructive proof of the theorem should produce a continuous
vector field on the plane, so the constructive content of $\forall x \exists y\, (y \neq x)$ 
is nontrivial.  Notice how the proof fulfills this prediction. 
\medskip

\noindent{\em Proof}.
Let $\alpha$ and $\beta$ be two of the three constants used in the dimension axioms,
and define
$$ e(x) = ext(\alpha,\beta,\alpha,x)$$
Since $\alpha \neq \beta$, axiom (A4) applies.   Let $d=e(x)$; we claim $d \neq x$. 
By (A4) we have $\T(\alpha,\beta,d)$ and $\beta d = \alpha d$.  Then the subsegment 
$\beta d$ is congruent to the whole segment $\alpha d$, contrary to Lemma~\ref{lemma:subsegment}.
That completes the proof that $e(x) \neq x$.
 
Define
$$ e_2(a,b,c) := ext(e(a),a,b,c).$$
Then the segment with endpoints $e(a)$ and $a$ is not a null segment, so $e_2(a,b,c)$ is 
everywhere defined, and if $d=e_2(a,b,c)$, we have $ad = bc$ by (A4). 
That completes the proof.

\subsection{Constructing a point not on a given line}
Consider the proposition that for every line $L$ there exists a point $c$ not on $L$.
In Tarski's language that becomes 
$$\forall a,b (a \neq b \implies \exists c\,(\neg\, Col(a,b,c))).$$
Classically, the theorem is a trivial consequence of the lower dimension axiom (A8),
which gives us three non-collinear points $\alpha$, $\beta$, and $\gamma$.  One of 
those points will do for $c$.  But that argument is not constructively valid, since 
it uses a case distinction to consider whether $a$ and $b$ both lie on $\Line(\alpha,\beta)$
or not.  It is an interesting example, because it illustrates in a simple situation exactly what 
more is required for a constructive proof than for a non-constructive proof.  For a 
constructive proof,we would need to find a uniform ruler and compass 
construction 
that applies to any  two  points $a$ and $b$ (determining a line $L$),
and produces a point not on $L$.

If we could erect a perpendicular to line $L$ at $a$, then (since lines are given by two points)
we would already have constructed a point off $L$.  Gupta constructs perpendiculars without 
circles: maybe he has solved the problem?  No, as it turns out.  Gupta's construction has to 
{\em start} with a given point $p$ not on $L$.  He shows how to construct a perpendicular to 
$L$ at $a$, but the first step is to draw the line $ap$.   The same is true for all the 
constructions of erected perpendiculars discussed above.

\begin{Lemma}\label{lemma:pointoffline} Given $a \neq b$, there exists
a point $c$ not collinear with $a$ and $b$.
\end{Lemma}

\noindent{\em Remarks}.  We do not know a direct construction; the proof 
we give uses the introduction of coordinates.  It does in principle provide 
a geometric construction, but it will be very complicated and not visualizable.
If we assume 
circle-circle continuity, we have an easy solution:  by the method of Euclid I.1 we produce
an equilateral triangle $abc$, whose vertex $c$ can be shown not to lie on line $L$. 
But we note that the  proof of circle-circle continuity that we give below 
via the radical axis requires a point not on the line connecting the centers to get 
started,  so cannot be used to prove this lemma.  In order to prove circle-circle
continuity without assuming an ``extra'' point, we also need to introduce coordinates.

\medskip

\noindent{\em Proof}.  Let $\alpha$, $\beta$, and $\gamma$ be the three 
pairwise non-collinear points guaranteed by (A8).  Let $\Line(\alpha,\beta)$ 
be called the $x$-axis.   Let $0$ be another name for $\alpha$ and 1 another name
for $\beta$.  Since $\gamma$ does not lie on the $x$-axis, we can use it to erect 
a perpendicular $0p$ to the $x$-axis at 0, on the other side of the $x$-axis from $\gamma$.
Call that line the $y$-axis.  Let $i = ext(p,0,0,1)$.   (Then $i$ is on the same side
of the $x$-axis as $\gamma$.)  Using the uniform perpendicular construction, and 
the point $i$ not lying on the $x$-axis, we define $X(p)$ to be the point on the 
$x$ axis such that the perpendicular to the $x$-axis at $X(p)$ passes through $p$.
Similarly we define $Y(p)$ using the point 1 not on the $y$-axis.  As shown 
in \cite{beeson-bsl}, using the parallel axiom we can define a point $p=(x,y)$
given points $x$ and $y$ on the $x$-axis and $y$-axis, such that $X(p)=x$ and 
$Y(p)=y$, and define addition and multiplication on the $X$-axis and prove their 
field properties.  Then coordinate algebra can be used in geometry.  Given 
distinct points with coordinates 
$(a,b)$  and $(p,q)$ determining line $L$, we can calculate the coordinates of 
a point not on $L$,  for example $(a,b) + (b-q,p-a)$.  That completes the proof.
\medskip

Applying our metatheorems below to this lemma, we see that there  
is a term $t(a,b)$ such that $A(a,b,t(a,b))$ is provable.
Of course, since the theorem is classically provable, by Herbrand's theorem there must be a finite
number of terms, such that in each case one of those terms will work, and indeed, the 
three constants $\alpha$, $\beta$, and $\gamma$ illustrate Herbrand's theorem in this case:
$$ a \neq b \implies \neg Col(a,b,\alpha) \ \lor \ \neg Col(a,b,\beta) \ \lor \ \neg Col(a,b,\gamma).$$
 But constructively, the matter is much more delicate.

\subsection{Hilbert planes and constructive geometry}
\label{section:HilbertPlanes}
In this paper, we have considered line-circle as an axiom.
Classically, there is a tradition of studying the consequences of (A1)-(A9) alone, 
which is known as the theory of Hilbert planes;  this theory corresponds to Hilbert's 
axioms without any form of continuity and without the parallel axiom.  The question to 
be considered here is whether there is an interesting constructive geometry of Hilbert planes.
There should be such a theory, with ruler and compass replaced by ``Hilbert's tools'', 
which permit one to extend line segments and ``transport angles'', i.e.,to construct a 
copy of a given angle with specified vertex $b$ on a specified side of a given line $L$.
That ``tool'' corresponds to a Skolem function for Hilbert's axiom C3.  

Indeed, most of the development in \cite{schwabhauser} from A1-A9 is perfectly constructive.
In particular, we can prove Hilbert's C3 and the related ``triangle construction theorem''
enabling us to copy a triangle on a specified side of a line, where the side is 
specified by giving a point not on the line. Hence, there is a viable 
constructive theory of Hilbert planes.  Developing that theory, however, is 
beyond the scope of this paper; even the classical theory of Hilbert planes is 
a difficult subject.

 It is curious that if Hilbert's C3 is 
classically weakened by removing the (classically) superfluous hypothesis about ``a specified
side of the line'',  so that it just requires being able to copy a triangle $abc$ to
a congruent triangle $ABC$ with $AB$ on $L$ and $A$ given,  then it becomes much more 
difficult to prove constructively, because we first have to construct a point not on $L$,
which (apparently) requires the introduction of coordinates, and hence the parallel axiom.
We do not know whether A1-A9 can prove that for every line $L$, we can construct a point
not on $L$; we do know that coordinates cannot be introduced on the basis of A1-A9 alone,
so the proof that works for A1-A10 fails for A1-A9.

\section{Metatheorems}

In this section, we prove some metatheorems about 
the two 
Skolemized constructive theories of Tarski geometry, i.e., either intuitionistic Tarski
with Skolem functions, or ruler and compass Tarski.  Both theories have line-circle and 
circle-continuity with terms for the intersections, and a Skolem function symbol $center$ for 
the triangle circumscription principle;  ruler and compass Tarski has  the logic of partial 
terms and a symbol $\il(a,b,c,d)$ for the intersection point of two lines, while
intuitionistic Tarski with Skolem functions has a Skolem function symbol  $ip$ for inner Pasch.  
Straightedge and compass constructions correspond to terms of ruler and compass Tarski;
we have shown that these can all be imitated by terms of intuitionistic Tarski with Skolem 
functions, i.e.,$\il$ is eliminable.

\subsection{Things proved  to exist can be constructed} 
In this section we take up our plan of doing for constructive geometry
 what cut-elimination and recursive realizability
did for intuitionistic arithmetic and analysis, namely, to show that existence proofs lead to programs (or terms)
producing the object whose existence is proved. In the case of constructive geometry, we want to produce geometrical constructions,
not just recursive constructions (which could already be produced by known techniques, since 
geometry is interpretable in 
Heyting's arithmetic of finite types, using pairs of Cauchy sequences of rational numbers as points).

\begin{Theorem}[Constructions extracted from proofs] \label{theorem:extraction} 
Suppose intuitionistic Tarski geometry  with  Skolem functions
 proves $$P(x) \implies \exists y\, \phi(x,y)$$ where $P$ is negative.
 Then there is a term $t(x)$ of intuitionistic
Tarski geometry with Skolem functions  
such that $$P(x) \implies \phi(x,t(x))$$ is provable.
\end{Theorem}

\noindent{\em Remark}. We emphasize that there is a {\em single} term $t(x)$.  That
corresponds to a {\em uniform construction}, that applies without case distinctions on $x$.
We shall see in the next theorem that things proved to exist classically can also be 
constructed (under appropriate conditions), but that several different constructions may 
be needed, for different cases on $x$.
\medskip

\noindent{\em Proof}.
We use cut-elimination.  Readers unfamiliar with cut-elimination, or desiring a 
specific axiomatization, are referred to Chapter XV of \cite{kleene1952}. Cut-elimination
works with ``sequents'',  which we write $\Gamma \seq \Delta$.  
(Kleene wrote $\Gamma \rightarrow \Delta$, but we use $\rightarrow$ for implication.)
$\Gamma$ and $\Delta$ are finite lists of formulas; for intuitionistic systems, $\Delta$
contains at most one formula $\phi$, so we also write $\Gamma \seq \phi$.

Since our axiomatization is quantifier-free,
if $\psi \seq \exists y \, \phi$ is provable, then there is 
a list $\Gamma$ of quantifier-free axioms such that $\Gamma , \psi \seq \exists y \, \phi$ is provable by a cut-free (hence 
quantifier-free) proof.  If the existential quantifier is introduced at the last step
of that proof,  then we obtain the desired proof just by omitting the 
last step of the proof.   Hence (as always in such applications of cut-elimination) 
we need to be able to permute the inferences until the existential quantifier is indeed 
introduced at the last step.  This issue was studied in a fundamental paper by Kleene
\cite{kleene1951}, who verified that the desired permutations are possible except for 
a few combinations of connectives, all involving disjunction or $\exists$.  Hence,
if the axioms are purely universal and disjunction-free, and $P$ is negative, we do not have
a problem on the left side of the sequent.   Certain occurrences of $\exists$ also do not 
make trouble, and in first position on the right side of the proved sequent is one of 
the harmless cases.  Hence,  by \cite{kleene1951}, we can permute the inferences as desired.
That completes the proof.  The hard part of the
 work was in arranging the axiom system to be quantifier-free
and disjunction-free. 
\medskip

{\em Remark}. To see that inferences cannot {\em always} be permuted, consider the example 
of proving that there is a perpendicular to line $L$ through point $p$,  with the left 
side $\Gamma$ of the sequent saying
that $p$ is either on $L$ or not.  Using the separate dropped-perpendicular and erected-perpendicular
constructions we can prove both cases, and then finish the proof, introducing $\lor$ on the 
left at the last step.  Explicitly, if $\psi$ says that $p$ is not on $L$, and $\phi$ says
$Line(p,q)$ is perpendicular to $L$ and contains $p$, then both sequents 
$\psi \seq \exists p,q\, \phi$
and 
$\neg \psi \seq \exists p,q\, \phi$ 
are provable, and hence 
$\psi \lor \neg \psi \seq \exists p,q\, \phi$
is provable.  
 But this proof obviously does not produce a single construction that
works in either case, so we would not expect to be able to permute the introduction of $\exists$
on the right with the introduction of $\lor$ on the left.
\medskip

The term $t(x)$ in the preceding theorem represents a geometrical construction, but the 
points constructed by intersecting lines are always given either by $center$ or $ip$ terms,
so the construction contains a ``justification''  for the fact that the lines intersect.  
On the other hand, the construction cannot be read literally as a construction script, but 
requires extra steps to construct the lines implicit in the $center$ and $ip$ constructions.
Moreover, there is nothing in the theorem itself to guarantee that the ``definedness conditions''
for $t(x)$ are met, since the Skolem functions are total.  The following theorem about 
ruler and compass Tarski geometry does not have that defect, since that theory uses the logic 
of partial terms.

\begin{Theorem}[Constructions extracted from proofs] 
 \label{theorem:extraction2}  
Suppose intuitionistic ruler-and-compass Tarski geometry 
 proves $$P(x) \implies \exists y\, \phi(x,y)$$
where $P$ is negative
(does not contain $\exists$ or $\lor$).  Then there is a term $t(x)$ of intuitionistic
ruler and compass Tarski geometry 
such that $$P(x) \implies \phi(x,t(x))$$ is also provable.  Moreover, 
if the proof of $P(x) \implies \exists y\, \phi(x,y)$ does not use certain axioms, 
then the term $t(x)$ does not involve the Skolem symbols for the unused axioms.
\end{Theorem}

\noindent{\em Remark}. If the formula $P(x) \implies \exists y\, \phi(x,y)$ does not 
contain any Skolem function symbols at all, i.e., it lies in Tarski's original language,
then it lies in both   Tarski geometry with Skolem symbols and 
ruler-and-compass Tarski geometry.  Thus we have a choice whether to realize the 
existential quantifiers using the symbols $ip$ and $center$, as in the former theory,
or using the symbol $il$.
\medskip

\noindent{\em Proof}.  
We have a choice of two proofs.  We could use
cut-elimination directly, but then we need it for the logic of partial terms and not just for 
ordinary intuitionistic predicate calculus.   The details of the 
cut-elimination theorem for such logics have not been published, but they are not significantly different from Gentzen's
formulation for first-order logic. While we are explaining this point, it is no more 
complicated to explain it for multi-sorted theories with LPT, which were used in \cite{beeson-kobe}
with axioms for Hilbert-style geometry.  Specifically, we reduce such theories to ordinary predicate calculus
as follows:  introduce a unary predicate for each sort, and then if $t$ is a term of sort $P$,
interpret $t \defined$ as $P(t)$. Now we have a theory in first-order one-sorted predicate calculus,
which is quantifier-free and disjunction-free if the original theory was, and we can apply ordinary 
cut-elimination, as in the proof of Theorem~\ref{theorem:extraction}.

Alternately, we can avoid using cut-elimination for LPT, by first translating the original 
formula from ruler and compass Tarski to intuitionistic Tarski with Skolem functions using 
Theorem \ref{theorem:il-elimination}. 
Then the resulting construction term $t$ involves the function symbol $ip$; but it is easy
to express $ip$ in terms if $\il$ if that is desired, i.e., to interpret intuitionistic Tarski
with Skolem functions into ruler and compass Tarski.  That completes the proof.
\medskip

\subsection{Extracting constructions from classical proofs}
The following theorem illustrates the essential difference between constructive
and classical (non-constructive) geometry:  in a constructive existence theorem, we must supply
a single (uniform) construction of the point(s) whose existence is asserted, but in 
a classical theorem,  there can be several cases, with a different construction in each case. 

\begin{Theorem} [Constructions extracted from classical proofs]
\label{theorem:extraction-classical}
Suppose classical Tarski geometry with Skolem functions proves $$P(x) \implies \exists y\, \phi(x,y)$$  where $P$ is quantifier-free and 
disjunction-free.
Then there are terms $t_i(x)$  
such that $$P(x) \implies \phi(x,t_1(x)) \lor \ldots \lor \phi(x,t_n(x))$$ is also provable.
\end{Theorem}
\noindent{\em Proof}.  This is a special case of Herbrand's theorem. 
\smallskip

\noindent{\em Example 1}. There exists a perpendicular to line $L$ through point $p$.  Classically,
one argues by cases: if $p$ is on $L$, then we can ``erect'' the perpendicular, and if $p$ is 
not on $L$ then we can ``drop'' the perpendicular.  So the proof provides two constructions, 
$t_1$ and $t_2$.  This is not a constructive proof.  In this paper, we have given a construction
(in fact two different constructions) of the ``uniform perpendicular''.   This constructive 
proof of the existence of a perpendicular to $L$ through $p$ provides a single term for the 
construction, instead of two terms (one of which works in each case). 
\smallskip

\noindent{\em Example 2}.  Euclid's proof of Book I, Proposition 2 provides us with two such constructions, $t_1(a,b,c) = c$ 
and $t_2(a,b,c)$ the result of Euclid's construction of a point $d$ with $ad=bc$, valid if $a \neq b$.   Classically 
we have $\forall a,b,c \, \exists d (ad=bc)$, but we need two terms $t_1$ and $t_2$ to cover all cases. 
\smallskip

{\em Example 3}.  Let $p$ and $q$ be distinct points and $L$ a given line, and $a$, $b$, and $c$ points on $L$, with $a$ and $b$ on the 
same side of $L$ as $c$.  Then there exists a point $d$ which is 
equal to $p$ if $b$ is between $a$ and $c$ and equal to $q$ if $a$ is between $b$ and $c$.    The two terms $t_1$ and $t_2$
for this example can be taken to be the
variables $p$ and $q$.  One term will not suffice, since $d$ cannot depend continuously
on $a$ and $b$, but all constructed points do depend continuously on their parameters.   This classical theorem is therefore not 
constructively provable.  

\subsection{Disjunction properties}
We mentioned above that intuitionistic Tarski geometry cannot prove any non-trivial disjunctive theorem.  That is a simple consequence of 
the fact that its axioms contain no disjunction.  We now spell this out:

\begin{Theorem}[No nontrivial disjunctive theorems]
\label{theorem:disjunctionpropertyECG} $\mbox{}$
Suppose intuitionistic Tarski geometry proves $H(x) \implies P(x) \lor Q(x)$, where $H$ is negative.  Then either 
$H(x) \implies P(x)$ or $H(x) \implies Q(x)$ is also provable.   (This result depends only on the 
lack of disjunction in the axioms.)
\end{Theorem}

\noindent{\em Proof}.   
Consider a cut-free proof of $\Gamma, H(x) \implies P(x) \lor Q(x)$,
where $\Gamma$ is a list of some axioms.
Tracing the disjunction upwards in the proof, if we reach a place where the disjunction was introduced on the right before
reaching a leaf of the proof tree, then we can erase the other disjunct below that introduction, obtaining a proof of one 
disjunct as required.  If we reach a leaf of the proof tree with $P(x) \lor Q(x)$ still present on the right, then it 
occurs on the left, where it appears positively.   Its descendants will also be positive, so it cannot participate in application of the rule for proof by cases (which introduces $\lor$ in the left side of a sequent); and it 
cannot reach left side of the bottom sequent, namely $\Gamma, H(x)$,  as these formulas contain no disjunction.  But 
a glance at the rules of cut-free proof, e.g. on p. 442 of \cite{kleene1951},  will show that these are the only possibilities.
That completes the proof.

We note that order on a fixed line $L$ can be defined using betweenness, so it makes sense
to discuss the provability of statements about order.

\begin{Corollary} Intuitionistic Tarski geometry does not prove apartness 
$$a < b \implies x < b \lor a < x.$$
\end{Corollary}

\noindent{\em Proof}.   The statement in question is a disjunctive theorem, so 
Theorem~\ref{theorem:disjunctionpropertyECG} applies.
\medskip

\begin{Corollary} Intuitionistic Tarski geometry
does not prove the principle $x \neq 0 \implies x < 0 \ \lor \ x > 0$
or the equivalent principle that if point $p$ does not lie on line $L$, then any other point $x$
is either on the same side of $L$ as $p$ or the other side.
\end{Corollary}

\noindent{\em Proof}.  The statement in question is a disjunctive theorem, so 
Theorem~\ref{theorem:disjunctionpropertyECG} applies.
\medskip

\subsection{Interpretation of Euclidean field theory}
A Euclidean field is defined constructively as an ordered ring in which nonzero elements 
have reciprocals.  The relation $a < b$ is primitive;  $a \le b$ abbreviates $\neg b < a$.
The axioms of Euclidean field theory include stability of equality and order.  Stability of order,
that is $\neg b \le a \implies a < b$, is also known as Markov's principle.
Classically, the models of ruler and compass geometry are planes over Euclidean fields.
We showed in \cite{beeson-bsl} that a plane over a Euclidean field is a model of ruler-and-
compass geometry, when ruler and compass geometry is defined in any sensible way;  constructively,
this theorem takes the form of an interpretation $\phi \mapsto \bar \phi$ from some geometric
formal theory to 
the theory \EF\ of Euclidean fields.  

The converse direction is much more difficult;  we have to show that any model of geometry 
is a plane over a Euclidean field $F$.  To do that, we fix a line $F$ to serve as the
$x$-axis  (and the domain of the field);  fix a point $0$ on that line, erect a perpendicular $Y$
to $F$ at 0 to serve as the $y$-axis.   Given any pair of points $(x,y)$ on $F$, we 
rotate $y$ by ninety degrees to a point $y^\prime$ on the $y$-axis, and then erect perpendiculars
at $x$ to $F$ and at $y^\prime$ to $Y$.   These perpendiculars should meet at a point 
$\MakePoint(x,y)$.   It is possible to show by the strong parallel axiom that they do meet.
This construction is the starting point for the following theorem:

\begin{Theorem}  Every model of intuitionistic Tarski geometry
 is a plane over a Euclidean field.  Moreover,
there is an interpretation $\phi \mapsto \phi^\circ$ from the theory 
of Euclidean fields to intuitionistic Tarski geometry. 
\end{Theorem}

\noindent{\em Proof}.  In addition to introducing coordinates as discussed above, one also
has to define addition and multiplication geometrically in order to interpret the 
addition and multiplication symbols of Euclidean field theory.  It has been shown in 
\cite{beeson-bsl} how to do this; the proofs there can be formalized in intuitionistic
Tarski geometry,  so we 
obtain a model-theoretic characterization of the models of that theory. 

Moreover, our work with the double-negation interpretation above can now be put to good use.
For example, the definition of multiplication can be given directly following Hilbert's 
definition, which is based on the triangle circumscription principle.  It is easy
to give a term $\HilbertMultiply(a,b)$  that takes two points $a$ and $b$ on a fixed 
line (the ``$x$-axis'') and produces their product (also a point on the $x$-axis), 
using $center$ and the uniform rotation construction. (See \cite{beeson-bsl} for details.)
But once that term is given, the assertions that it satisfies the associative and 
commutative laws are quantifier-free, and hence, the proofs in \cite{schwabhauser} 
are ``importable.''   Technically, one must check that the degenerate cases of inner Pasch
are not used, but that is all that one has to check by hand.  In \cite{beeson-bsl},
there is a definition of ``uniform addition'', i.e., without a case distinction on the signs
of the addends.  A term $\Add(x,y)$ defining the sum of $x$ and $y$ is given in 
\cite{beeson-bsl}.  Again, once the term is given, we can be assured by the double-negation
interpretation that its properties are provable in intuitionistic Tarski with Skolem 
functions, if we just check \cite{schwabhauser} to make sure the degenerate cases of 
inner Pasch are not used.

The terms $\Add$ and $\HilbertMultiply$ can then be used to define a
syntactic interpretation $\phi \mapsto \phi^\circ$ from the theory of Euclidean fields 
to intuitionistic Tarski geometry.  That 
completes the proof of the theorem.

\section{Circle-circle continuity} 

In this section we show that circle-circle continuity is a theorem of
intuitionistic Tarski geometry; that is,  
 we can derive the existence of the intersection points
of two circles (under the appropriate hypotheses).  The similar 
theorem for classical Tarski geometry can be derived indirectly, using
the representation theorem (Theorem~\ref{theorem:representation}) and 
G\"odel's completeness theorem; but for intuitionistic Tarski geometry,
we must actually exhibit a construction for the intersection points of 
two circles, and prove constructively that it works.  This question 
relies on Euclid III.35, a theorem about how two chords of a circle 
divide each other into proportional segments, and III.36, a similar
theorem, and constructively it requires a combined version of those 
two propositions without a case distinction as to which applies (i.e.,
two lines cross inside or outside a circle).
  Aside from those theorems,
it uses only very straightforward geometry.  

\begin{Theorem} \label{theorem:circlecircle}
 In intuitionistic Tarski geometry with only line-circle 
continuity:  given two circles $C$ and $K$ satisfying the hypotheses of 
circle-circle continuity,  and a point not on the line $L$ connecting
their centers,  we can construct the point(s) of intersection of $C$ and 
$K$.  
\end{Theorem}

A point not on $L$ can of course be constructed, by Lemma~\ref{lemma:pointoffline},
but to do so we have to introduce coordinates, which we wish to avoid here
for esthetic reasons, so here we just leave that point as a parameter of the construction.

There are two approaches to proving this theorem, which we will discuss separately.
The first method proceeds by introducing coordinates and reducing the problem 
to algebraic calculations by analytic geometry.  While in principle this does
produce an ultimately purely geometric proof, one cannot visualize the lengthy
sequence of constructions required.   For esthetic reasons, therefore, we also
give a second proof, which avoids the use of coordinates by the use of a 
well-known construction called the ``radical axis.''  In this more geometric
proof, essential use is made of the ``extra point'' in the hypothesis.%
\footnote{
Both these proofs use the parallel axiom essentially.  There is a proof
in the literature that line-circle continuity implies circle-circle continuity
without the use of the parallel axiom \cite{strommer1973}.  It is based on 
Hilbert's axioms rather than Tarski's, and we do not know if it is constructively
justifiable or not.  But we have no reason to avoid the parallel axiom for
present purposes.} 
\medskip

\subsection{Circle-circle continuity via analytic geometry}
Given two circles $C$ and $K$ with distinct centers $s$ and $t$, let 
$L$ be the line through the centers.  Given a point not on $L$, we can 
erect a perpendicular to $L$ at $s$, and introduce coordinates in the 
manner of Descartes and Hilbert,  with constructive extension to negative
arguments as developed in \cite{beeson-bsl}.  We can choose the point $t$
as $1$.  Now the tools of analytic geometry are available.  Let $r$
be the radius of circle $C$ and $R$ the radius of circle $K$, and 
calculate the equations of $C$ and $K$ and solve for a point $(x,y)$ 
lying on both circles.  It turns out that some crucial terms cancel, and 
we can solve the equations using only square roots,  which means that we 
can solve them geometrically using the methods of Descartes and Hilbert,
with the constructive modifications {\em op.~cit}.   To derive circle-
circle continuity, we must show that the hypothesis that circle $C$ has 
a point inside circle $K$ makes the quantities under the square root non-negative,
and the extra hypothesis that the circle $C$ has a point strictly inside $K$
makes the two solutions of the equation distinct.  

 A similar theorem is proved classically in \cite{hartshorne}, p. 144,
but a few details are missing there.    
 The issue is that it is not enough to observe that the equations for the intersection points are 
quadratic.  One has to translate the hypothesis that one circle has a point inside and a point outside the 
other circle into algebra and show algebraically that this implies the equations for the intersection are solvable. 

 This is a fairly routine exercise and is all perfectly 
constructive; but it is a bit long, and besides, it
 is somewhat unsatisfying to have to resort to coordinates. One would like
to see a direct geometric construction of the points of intersection of two 
circles, using only a few steps, rather than the dozens or perhaps hundreds 
of not-visualizable steps required to geometrize an algebraic calculation.
There is indeed such a geometric construction, using the ``radical axis'' of 
the two circles.  Below we verify that the radical axis construction can 
be carried out constructively (i.e., does not require any case distinctions in its definition),
and that the correctness proof can be carried out in intuitionistic Tarski geometry.
Although the construction itself is easy to visualize (it is only a few steps with 
ruler and compass),  the correctness proof in Tarski geometry is more complicated.

\subsection{Euclid Book III in Tarski geometry}

The correctness proof of the radical axis construction 
requires the last two propositions of Euclid Book III; and 
moreover the formalization of Book III in Tarski geometry is of independent interest.

Book III of Euclid can be formalized in intuitionistic Tarski geometry, 
but since most of the theorems mention angles, we need to use the developments
of Chapter~11 of \cite{schwabhauser}, where angles, angle congruence, etc. are 
developed.  However, some propositions can be proved quite simply, for 
example III.31 (an angle inscribed in a semicircle is right), which 
goes back to I.29 and I.11 and hence to the construction of perpendiculars
(and not to Book II at all).   We also will need a related proposition 
  that might have been (but does not seem to be) in Euclid:

\begin{Lemma} \label{lemma:31b}
 (in Tarski geometry with segment-circle continuity)
If $axb$ is a right angle and $ab$ is a diameter of a circle $C$ then $x$ lies on $C$.
\end{Lemma}

\noindent{\em Proof}.
By segment-circle continuity, we can 
find a point $y$ on $C$ and on the ray from $a$ through $x$.  Then by Euclid~III.30,
$ayb$ is a right angle, so if $y\neq x$ then
 $xb$ and $yb$ are two perpendiculars to $ay$ through $b$,
contradiction.  By the stability of equality we have $y=x$, so $x$ lies on $C$.
That completes the proof.

We present here another proposition from Euclid Book III that can be proved
directly from the Tarski axioms.

\begin{Lemma}[Euclid III.18] \label{lemma:euclidIII.18} Suppose line $L$ meets
circle $C$ with center $c$  in exactly one point $a$. Then $ca \perp L$.
\end{Lemma}

\noindent{\em Remark}. We need dropped perpendiculars to prove this lemma,
but we already derived the existence of dropped perpendiculars from line-circle
continuity, so that is not a problem.  The proof uses Euclid's idea, but 
Tarski's definition of perpendicular.
\medskip

\noindent{\em Proof.}  Drop a perpendicular $cb$ from $c$ to line $L$,
which can be done by Lemma~\ref{theorem:droppedperp}. 
We want to prove $b=a$.  By the stability of equality, we can proceed by 
contradiction, so suppose $b \neq a$.  
 Point $b$ is outside the circle (i.e.,there is a point $e$ 
on $C$ with $\B(c,e,b)$,  since otherwise by line-circle continuity, $L$
meets $C$ in a second point.  Let $e$ be the reflection of $a$ in $b$;
then $ab = be$, and since $ce \perp L$, we have $ca = ce$.  That is, $e$
lies on circle $C$ as well as line $L$, contradicting the hypothesis 
that $L$ meets $C$ only once.  That completes the proof.
\medskip

Despite these   examples, 
the radical axis construction that we use makes use of some developments
of Euclid Book III that are not so straightforward, because they
rest on the theory of proportionality 
in Euclid Book II.  We want to make sense of the phrase 
$$ab \cdot cd  = pq \cdot rs.$$ 
In Euclid, this is written  ``the rectangle 
contained by $ab$ and $cd$ is equal to the rectangle contained by $pq$ and $rs$.''
But Euclid has no concept of ``area'' as represented by (the length of)
a segment.  To define this relation 
geometrically,  we could use a definition of similar triangles, involving two right 
triangles with sides $ab$, $pq$ and $cd$, $rs$ respectively. 
 Book II of Euclid develops (quadratic)  algebra on that basis.   

There is another way to define the notion $ab \cdot cd  = pq \cdot rs$.  
Namely,  introduce coordinates,
define multiplication geometrically, and interpret $ab \cdot cd$ as   multiplication
of segments on the $x$-axis congruent to $ab$ and $cd$.   This is not actually 
so different from the first (Euclidean) interpretation, since similar triangles are
used in defining multiplication.    
 For the modern analysis of Euclid's notion of equality 
for ``figures'' see page~197 of \cite{hartshorne}, and for the connection 
to geometric multiplication see page~206. There it is proved that equality in 
Euclid's sense corresponds to algebraic equality using geometric arithmetic;
in particular, the two definitions of $ab \cdot cd  = pq \cdot rs$ are provably
equivalent.  

To carry out the radical axis construction constructively, we need to extend
the notion $ab \cdot cd$ to allow signed segments.  To do this directly using 
similar triangles would be to duplicate the effort of defining signed multiplication 
in \cite{beeson-bsl}.  Therefore we use the geometric-multiplication definition,
following \cite{hartshorne}.

Either definition can be expressed in a quantifier-free way using 
intuitionistic Tarski geometry with Skolem functions,  so the two notions are 
provably equivalent in intuitionistic Tarski geometry if and only if 
they are provably equivalent in classical Tarski geometry.  Since Hartshorne
{\em op.~cit.} proves them equivalent in classical Hilbert geometry (and 
without using circle-circle continuity), and since Hilbert's axioms (except continuity)
are provable
in classical Tarski geometry (without continuity) (as shown in \cite{schwabhauser}),
it follows that the two definitions are provably equivalent in intuitionistic
Tarski geometry.  
 
We need to define the notion of the power of a point in intuitionistic
Tarski geometry, and check that its principal properties can be proved there.
 That notion is usually defined as follows, when $b$ is not the 
center of $C$:  Let $c$ be the center of 
$C$, and $x$ the point of intersection of $\Line(b,c)$ with $C$ that is on the same 
side of $c$ as $b$, and $y$ the other point of intersection.
 Then the power of $b$ with respect to $C$ is $bx \cdot by$.
If we interpret the dot as 
 (signed)  multiplication (see \cite{beeson-bsl}) of (directed) segments
 on the $x$-axis of points congruent to the segments mentioned, then  
this definition makes sense in intuitionistic Tarski geometry.
It does, however, have the disadvantage that the power of the center of the 
circle is not defined; and we cannot just define it to be $-1$, as we could 
classically, because constructively we cannot make the case distinction whether $b$
is or is not the center.  This definition can be fixed constructively as follows.
Fix a diameter $pq$ of circle $C$, whose center is $c$.  Then given any point $b$,
extend segment $pc$ by $bc$ to produce point $B$.  (If $b=c$, this is still legal
and produces $c$.)  Then the power of $b$ with respect to $C$ is $Bq \cdot Bp$, 
where the dot is signed   multiplication (so that $Bq$ and $Bp$ have opposite
signs when $b$ is inside $C$.)
This gives the same answer as the usual definition when $b$ is not the center.

\begin{figure}[ht]
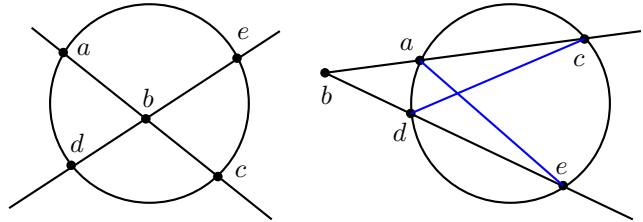

\caption{The power of $b$ with respect to $c$ can be computed from any chord,
because $ba\cdot bc = bd \cdot be$.}
\label{figure:powerfigure}
\PowerFigure
\end{figure}

The following lemma shows that the power of $p$ with respect to $C$ can 
be computed from any chord, not just from the diameter.  See Fig.~\ref{figure:powerfigure},
which shows separately the cases when $b$ is inside or outside $C$.

\begin{Lemma}\label{lemma:III.35} In Tarski geometry with only line-circle
continuity:  Let $C$ be a circle and 
let $b$ be any point.  Let $L$ be any line through $b$ meeting
$C$ in points $a$ and $c$.  Then the power of $C$ with respect to $C$ is 
$ba \cdot bc$.
\end{Lemma}

\noindent{\em Remark}.   For the case when $b$ is inside $C$ (first part of Fig.~\ref{figure:powerfigure}), this is Euclid III.35, and for $b$ outside it is III.36. 
The general case is mentioned in Heath's commentary on III.35 \cite{euclid1956} as a corollary of 
III.35 and III.36.  It also occurs as  Exercise~20.3 in \cite{hartshorne}. The proof
implicitly suggested there by Exercise~20.2 is the same one suggested in Heath's 
commentary.  
\medskip

\noindent{\em Proof}.  Hilbert multiplication can be defined by 
Skolem terms in intuitionistic Tarski geometry, without circle-circle continuity,
as the constructions in  \cite{beeson-bsl} show.   To recap:  multiplication is defined
by using the triangle circumscription axiom to draw a circle, and then the product 
is given by the intersection of that circle with a line.  In addition there are some 
rotations involved, which also can be defined by terms. 
  When formulated in intuitionistic Tarski geometry with Skolem functions,
the statement of the lemma is quantifier-free.   By the double-negation interpretation,
it is provable constructively if and only if it is classically provable.  And 
it is classically provable, cf. the exercise mentioned in the remark.  Of course,
the textbook containing the exercise is based on Hilbert's axioms, but \cite{schwabhauser} derives 
all of Hilbert's axioms from Tarski's (classically), so once you solve Exercise~20.3 
{\em op.~cit.},  
it follows that the result is provable in intuitionistic Tarski geometry with Skolem 
functions.   
 
\subsection{The radical axis}
In this section, we discuss the construction of the ``radical axis'' of two 
circles, with attention to constructivity.   In the next section we will 
use the radical axis to give a second proof that line-circle continuity
implies circle-circle continuity. 

\begin{figure}[ht]
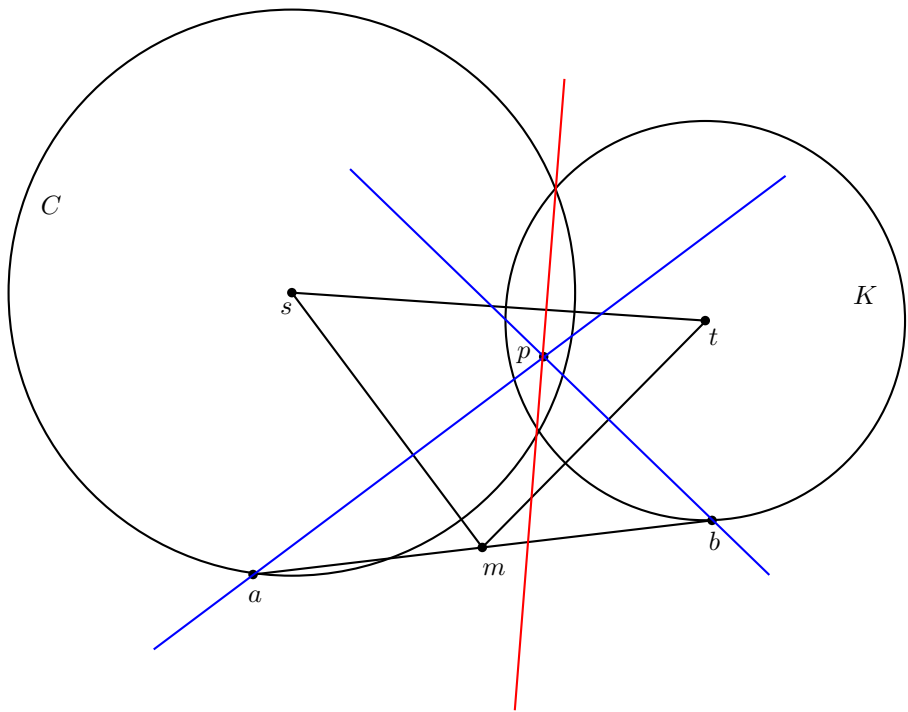

\caption{Construction of the radical axis.  \goodbreak $m = midpoint(a,b)$; draw
$sm$ and $mt$ and drop perpendiculars from $a$ to $sm$ and from $b$ to $mt$.
Their intersection is $p$ and the radical axis is perpendicular to $st$ through $p$.}
\label{figure:radicalaxis}
\RadicalAxisFigure
\end{figure}

 The ``radical axis'' of 
two circles is a line, defined whether or not the circles intersect, 
such that if they do intersect, the line passes through the points of intersection
(and if they are tangent, it is the common tangent line).  On page~182 of 
\cite{hartshorne}, a ruler and compass construction of the radical axis 
is given.  Fig.~\ref{figure:radicalaxis}
illustrates the construction, for the benefit of readers who do not have
\cite{hartshorne} at hand.%
\footnote{The radical axis was already old in 1826 \cite{steiner1826}, although
there it is constructed from the intersection points of circles, rather than 
the other way around.  I do not know the origin of the ruler and compass 
construction used here.}

\noindent
The initial data are the centers $s$ and $t$ of the two circles, with $s \neq t$, 
and two points $a$ and $b$ on the circles, such that $a \neq b$ and $ab$ does
not meet the line joining the centers.   (That 
hypothesis allows  one 
of the circles to be a null circle (zero radius), but not both). 
The construction is as follows:  First define $m$ as the midpoint of $ab$.
Then $\B(a,m,b)$.  If $s=m$ or $t=m$ then $ab$ meets $st$, contrary to hypothesis.
Hence we can
 construct the lines $sm$ and $tm$.  Then construct
 perpendiculars to those lines 
through $a$ and $b$ respectively (using the uniform perpendicular, so we do not 
need to worry if $a$ and $b$ are on those lines or not).    
Then $p$ is to be the intersection of these two perpendiculars, and the radical 
axis is the perpendicular to $st$ through $p$, again using the uniform perpendicular.  
\medskip

\begin{Lemma}\label{lemma:radicalaxisworks}  Given two circles $C$ and $K$,
the radical axis as constructed above does not depend on the particular points
$a$ and $b$ chosen,  and can be constructed from the two circles and one
additional point not lying on the line joining the centers.
\end{Lemma}

\noindent{\em Remark}.  The ``extra point'' is a necessary parameter.  The 
circles are presumed given by center and point, but the points giving the circles
might happen to lie on the center line.
\medskip

\noindent{\em Proof}.
 Let the two circles $C$ and $K$ have centers $s$ and $t$, 
and let $a$ and $b$ be distinct points on $C$ and $K$, 
respectively,  such that the midpoint $m$ of $ab$ does not lie on the line $L$
containing the centers $s$ and $t$.  We need a point $r$ not on $L$ to
be able to choose $a$ and $b$. For example, we can use $r$ to erect a perpendicular
to $L$ at $s$, and let $a$ be one of its intersections with $C$, and 
then construct a perpendicular to $L$ at $t$, on the same side of $L$ as $a$,
  and let $b$ be the intersection
of this perpendicular with $K$.
 
We wish to construct the radical 
axis $R$ of $C$ and $K$.  The first step is to construct the midpoint 
$m$ of $ab$, which can be done using only line-circle continuity by
Lemma~\ref{lemma:midpointsG}. Then $ms$ and $mt$ are not null segments, since
$m$ does not lie on $L$ (since $s$ and $t$ are on the same side of $L$,
so not on opposite sides of it).  Then we need to construct perpendiculars to $sm$ 
and $mt$ that pass through $p$.  By 
the uniform perpendicular construction, we can do that without
worrying about whether $p$ lies on $Line(s,m)$ or not, or whether $p$ 
lies on $Line(m,t)$ or not.
 But we do have to worry about whether the intersection
$p$ 
of those perpendiculars exists.   By the strong parallel axiom, it will exist 
if the two perpendiculars are not parallel or coincident.  That can only happen 
if $m$, $s$, and $t$ are collinear, but we have chosen $a$ and $b$ so that $b$
does not lie on $L$; hence indeed $p$ exists.  Then define line $R$ as the 
(uniform) perpendicular to $L$ through $p$.  

\begin{figure}[ht]
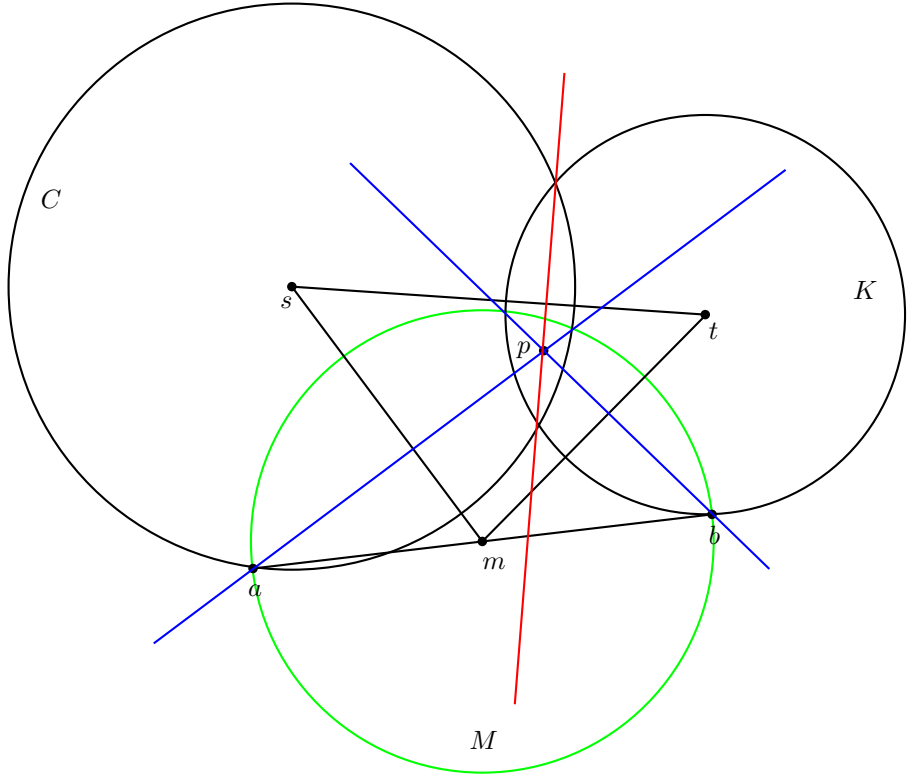

\caption{The power of  $p$ with respect to each of $C$ and $K$ is equal to the power
of $p$ with respect to $M$.}
\label{figure:radicalaxisprooffigure}
\RadicalAxisProofFigure
\end{figure}

Now we will prove (constructively and using only line-circle continuity) that 
every point $x$ on $R$ has equal powers with respect to $C$ and $Q$.  Suppose
$x$ is on the radical axis $R$.
Define circle $M$ to be the circle with center $m$ and passing through $a$ and $b$.
See Fig.~\ref{figure:radicalaxisprooffigure}.

\noindent
Let $z$ be the intersection of $M$ with $\Line(a,p)$, and let $y$ be the 
intersection of $M$ with $\Line(b,x)$.  Then the power of $x$ with respect 
to $C$, namely $xz \cdot ax$ (by Lemma~\ref{lemma:III.35})
  is equal to the power of $x$ with respect to $M$.
Similarly, the power of $x$ with respect to $K$ is $yx\cdot yb$, which is also the power
of $x$ with respect to $M$.   Since the powers of $x$ with respect to $C$ and to $K$
are both equal to the power of $x$ with respect to $M$, they are equal. 
That completes the proof.

\subsection{Circle-circle continuity via the radical axis}

We observe that the power of $x$ with respect to circle $C$ is negative when $x$ is 
inside $C$, zero when $x$ is on $C$, and positive when $x$ is outside $C$. 
Intuitively, the radical axis of two intersecting circles is the line joining 
the two intersection points.  In this section we use this idea to 
prove circle-circle continuity.  One more idea is necessary: the ``radical center.''

\noindent
{\em Second proof of Theorem~\ref{theorem:circlecircle}.}
Let circles $C$ and $K$ be given,
 satisfying the  hypothesis of the circle-circle continuity axiom.
 That is, $K$ has a point (say $b$) 
outside $C$ and a point (say $x$) inside $C$.   We will use $b$ as one of
the points to start the radical axis construction, and we will choose the 
other point $a$ very carefully.  
Construct the perpendicular to 
$xb$ at $x$.   Since $x$ is inside $C$,  by segment-circle continuity that perpendicular
meets $C$ in some point $a$ such that the interior of $ab$ does not meet $\Line(s,t)$.
In other words,   
$a$ and $b$ lie on the same side of $\Line(s,t)$ unless $b$ already lies on $\Line(s,t)$. 
   
 Now use $a$ and $b$ as starting points
for the radical axis construction; let the constructed point be $p$.   
By construction of $a$, if $x \neq b$, then 
angle $axb$ is a right angle with $a$ and $x$ at the ends of a diameter of circle 
$M$; that is the result of choosing
$a$ as we did.  See Fig.~\ref{figure:radicalaxisprooffiguretwo}. 
Therefore the vertex $x$  of this right angle lies on circle $M$, by Lemma~\ref{lemma:31b}.
On the other hand, if $x=b$ (that is $vb$ is tangent to $K$) then $x$ also lies on $M$;
by the stability of equality, $x$ lies on $M$ (that is $xt=bt$) whether or not $x=b$.  
 
\begin{figure}[ht]
\caption{$b$ and $x$ are given on $K$ with $b$ outside $C$ and $x$ inside $C$.
Then $a$ is chosen so $ax \perp xb$. Then
it turns out that $x$ lies on $M$, and the constructed point $p$ is inside both circles.}
\label{figure:radicalaxisprooffiguretwo}
\RadicalAxisProofFigureTwo
\end{figure}

By segment-circle continuity, there is a point $v$ on $C$ with $\B(v,x,p)$.
We claim $\T(v,p,b)$.  By the stability of betweenness we can prove this by 
contradiction.   Since $\neg \neg (A \lor B)$ is equivalent to $\neg (\neg A \land \neg B)$,
we can argue by cases for the contradiction.  There are two cases to consider:
$\B(p,v,b)$ and $\B(v,b,p)$.  Let $y$ be the other point of intersection (besides $v$)
of $Line(b,x)$ with $C$.  Then $\B(v,x,y)$ since $x$ is inside $C$.

Case 1: $\B(p,v,b)$.  We will show that the power of $p$ with respect to 
$C$ is less than  the power of $p$ with respect to $K$.   The former is 
$pv \cdot py$, the latter is $px \cdot pb$, and because $\T(v,x,y)$ we have $pv \le px$
and because $\T(p,y,b)$ we have $py \le pb$.  Therefore  
  $pv \cdot py \le px \cdot pb$.  For equality to hold, we would need 
$v=x$ and $y=b$, but we have arranged $x \neq v$.  Hence the power of $p$ with 
respect to $C$ is strictly less than the power of $p$ with respect to $K$,
contradiction.  

Case 2: $\B(v,b,p)$.  We similarly can show that the power of $p$ with respect to $K$
is less than the power of $p$ with respect to $C$.
Hence both cases are contradictory.  Hence $\T(v,p,b)$ as claimed.

Since $\T(v,p,b)$ and $\T(v,x,y)$ and $\T(x,p,b)$,  it follows classically 
that $\T(v,p,y)$ or $\T(y,p,b)$.  In the first case, $p$ is inside $C$; in 
the second case, $p$ is inside $K$.  But since the power of $p$ with respect to $C$
is equal to the power of $p$ with respect to $K$,  $p$ is inside $C$ if and only 
if it is inside $K$.  Double-negating each step of the argument, we find that 
$p$ is not not inside $C$;  but by the stability of ``inside'',   $p$ is inside $C$.

 Then by line-circle continuity, since 
$p$ lies on the radical axis $R$,  
$R$ meets $C$ in a point $x$.  Since $x$ lies on $C$,
the power of $x$ with respect to $C$ is zero.  Since $x$ lies on $R$, the power of 
$x$ with respect to $K$ is equal to the power of $x$ with respect to $C$, which is zero.
Hence $x$ lies on $K$ as well as on $C$.
That completes the proof. 
\medskip

We have shown that points on the radical axis have equal powers with respect to 
both circles.  The following lemma is the converse.  We do not need it, 
but the proof is short and pretty.  

\begin{Lemma} \label{lemma:radicalpower}
Let $C$ and $K$ be distinct circles.  If they meet, then the radical axis of $C$ and $K$ consists
of exactly those points whose powers with respect to $C$ and $K$ are equal.
\end{Lemma}

\noindent{\em Remark}.  The lemma is true even if the circles do not meet, but I 
do not know a simple geometric proof.
\medskip

\noindent{\em Proof}.  We have already proved that points on the radical axis
have equal powers.  It suffices to prove the converse.  Suppose that $u$ has
equal powers with respect to $C$ and $K$.  We must prove $u$ lies on the radical 
axis.  Let $v$ lie on both circles.  When we compute the powers of $u$ with respect
to both $C$ and $K$ using the line $uv$, we get different answers 
unless $u$ lies on the radical axis (so that the endpoints on $C$ and $K$ are the same).
That completes the proof.

\subsection{Skolem functions for circle-circle continuity}
Terms of Tarski geometry (intuitionistic or continuous, which has classical logic)
correspond to (certain) ruler and compass constructions; in effect, to constructions 
in which you can form the intersection of lines that must intersect by inner Pasch,
and intersections of lines and circles.  Since inner Pasch implies outer Pasch
(in the presence of other axioms) 
the points formed by outer Pasch are also given by terms; of course those terms
will involve the complicated constructions of Gupta's perpendiculars.  That is,
the Skolem function for outer Pasch, defined in terms of the Skolem function for 
inner Pasch, is very complicated.   Similarly, since circle-circle
continuity is implied by line-circle continuity, then there must be terms for 
constructing those intersection points as well.  

But there is an issue to consider, in that the terms for erected and uniform perpendiculars
have an extra parameter, a point not on the line;  and the radical axis construction 
also has an extra parameter for a point not on the center line.  Examination of those
constructions reveals that if the ``extra'' point is changed to the other side of the line,
then the ``head'' of the perpendicular changes sides too; and in the radical axis 
construction, the result is that the two intersection points of the two circles switch 
places.  If we then fix that choice once and for all by using Lemma~\ref{lemma:pointnotonL},
we can construct terms that give the two intersection points of two circles 
continuously.   However, those terms will be complicated, because the term 
for constructing a point not on a line involves coordinates and cross products. 

In previous work on constructive geometry, we had circle-circle continuity as an axiom,
and built-in function symbols for the intersection points.  It was a point of difficulty
to distinguish the two, which we wanted to do by saying whether the triple of points
from the center of one circle to the other center to the point of intersection was a 
``right turn'' or a ``left turn''.   The concepts $\Right(a,b,c)$ and $\Left(a,b,c)$
had to be defined, either by a complicated set of axioms, or by introducing coordinates
and using cross products as in \cite{beeson-bsl}.  Here we can recover those same 
terms, but now the coordinates and cross products are at least no longer in the axioms!
If we perform a complex conjugation (i.e.,reflect in the $x$-axis) then the two terms 
for the intersection points of two circles change places, exactly as in \cite{beeson-bsl}.

If one wishes (for example for connecting these theories to computer graphics) to have 
explicit function symbols for circle-circle continuity, of course they can be conservatively 
added.  

\section{Conclusion}
We have exhibited a constructive version of Tarski's 
Euclidean geometry.  Because of the double-negation interpretation, it 
can prove at least some version of each classical theorem. Using the 
uniform perpendicular, rotation, and reflection constructions given in this 
paper, it is possible (by the methods of \cite{beeson-bsl}) to give
  geometric 
definitions of addition and multiplication, without case distinctions as to the 
sign of the arguments,  and proofs of their properties,
so that coordinates in a Euclidean field provably exist.  Hence the theory 
has not omitted anything essential.
   To achieve these results, we had to modify Tarski's axioms to eliminate degenerate
cases, and add back some former axioms that Tarski had eliminated using those degenerate
cases.  Even with classical logic, this theory now connects nicely with ruler-and-
compass constructions, since each of the points asserted to exist can be constructed
with ruler and compass. 

By cut-elimination, things proved to 
exist (under a negative hypothesis, as is always the case in Euclid) can 
be constructed, by a uniform straightedge-and-compass construction.  Even
stronger, these constructions need not involve taking the intersections of 
arbitrary lines, but only those lines that have to intersect by the strong parallel
axiom or inner Pasch. 

By contrast, in Tarski's (classical) theory,  we obtain (by Herbrand's theorem) a similar result but without 
uniformity, i.e.,  there are several constructions (not necessarily just one),  such that 
for every  choice of the ``given points'', one of the constructions will work.  (The 
classical result (unlike the constructive one) holds only for formulas $\forall x \exists y\, A(x,y)$,
where $A$ is quantifier-free.)

These points-only axiom system have conservative extensions with variables for lines and circles,
and further conservative extensions with variables for angles, segments, and arcs,
which can serve for the direct constructive formalization of Euclidean geometry using
Hilbert's primitives (as in \cite{beeson-kobe}).
Therefore, this points-only theory, with its short list of axioms, can be said 
to provide the logical foundations of constructive Euclidean geometry.   In particular, 
it supplies one detailed example of a formalization of constructive geometry,
 to which the independence
results about the parallel postulate of \cite{beeson-bsl} apply.


\end{document}